\documentclass[leqno,11pt]{amsart}
\usepackage{amsmath}
\usepackage{amsfonts}
\usepackage{amssymb}
\usepackage{multicol}
\usepackage{amscd}
\usepackage{amsthm}
\usepackage{latexsym}
\usepackage{amsbsy}
\usepackage[dvips]{graphics}
\usepackage[dvips]{color}
\usepackage{lscape}
\newtheorem{theo}{Theorem}[section]
\newtheorem{lemm}[theo]{Lemma}
\newtheorem{prop}[theo]{Proposition}
\newtheorem{coro}[theo]{Corollary}
\theoremstyle{definition}
\newtheorem{defi}[theo]{Definition}
\newtheorem{example}[theo]{Example}
\newtheorem{nota}[theo]{Notation}

\theoremstyle{remark}
\newtheorem{rema}[theo]{Remark}
\newcommand{\bahram}{\textbf{(CR)}}
\begin{document}
\title[$C^{\star}$-algebras from Anzai flows]
{$C^{\star}$-algebras from Anzai flows and their $K$-groups}
\author[K. Reihani]{Kamran Reihani}
\address{Department of Mathematics\\
Tarbiat Modarres University\\
Tehran, P.O. Box: 14115-111\\
Iran}
\email{reyhan\_\hspace{2pt}k@modares.ac.ir}
\author[P. Milnes]{Paul Milnes}
\address{Department of Mathematics \\
University of Western Ontario \\
London, Ont. N6A 5B9\\
Canada}
\email{milnes@uwo.ca}
\date{}
\maketitle

\begin{abstract}
We study the $C^{*}$-algebra $\mathcal{A}_{n,\theta}$ generated by the Anzai 
flow on the 
\text{$n$-dimensional} torus $\mathbb{T}^n$. It is proved
that this algebra is a simple quotient of the group $C^{*}$-algebra of a 
lattice
subgroup $\mathfrak{D}_n$ of a \text{$(n+2)$-dimensional} connected simply 
connected nilpotent Lie
group $F_n$ whose \text{corresponding} Lie algebra is the generic filiform 
Lie algebra
$\mathfrak{f}_{n}$. Other simple infinite dimensional quotients of 
$C^{*}(\mathfrak{D}_n)$
are also characterized and \text{represented} as matrix algebras over simple 
affine Furstenberg
transformation group \text{$C^*$-algebras} of the lower dimensional tori.
The $K$-groups of the $\mathcal{A}_{n,\theta}$ and other simple quotients
of $C^{*}(\mathfrak{D}_n)$ are studied, the Pimsner-Voiculescu \text{6-term}
exact sequence being a useful tool. The rank of the \text{$K$-groups} 
of $\mathcal{A}_{n,\theta}$ is
studied as explicitly as possible, and is proved to be the same as for more 
general
transformation group \text{$C^*$-algebras} of $\mathbb{T}^n$ including the
Furstenberg transformation group \text{$C^*$-algebras} $A_{F_{f,\theta}}$. 
An error (about these \text{$K$-groups}) in the 
literature is addressed.
\end{abstract}
\section{Introduction}

In 3 dimensions there is  a unique (up to isomorphism) connected, simply 
connected,
nilpotent Lie group, which we call $G_3$ (following Nielsen  \cite{oN83}); 
$G_3$
$(\,=\Bbb R^3$ as a set) is the Heisenberg group with multiplication
$$(k,m,n)(k^{\prime},m^{\prime},n^{\prime})=(k+k^{\prime}+nm^{\prime}
,m+m^{\prime},n+n^{\prime}).$$
The faithful irreducible representations of the lattice subgroup $H_3$
$(\,=\Bbb Z^3$ as a set) of $G_3$ generate the irrational rotation algebras
$A_{\theta}$. In 4 dimensions there is also a unique such connected group 
$G_4$, in 5
dimensions there are 6 such groups $G_{5,i}$ for $1 \leq i \leq 6$, and in 6 
dimensions
there are 24 such groups $G_{6,j}$ for $1 \leq j \leq 24$.
The main thrust in \cite{pM93, pM97, pM02} was  to find
cocompact subgroups $H_4 \subset G_4$, $H_{5,i}\subset G_{5,i}$, and
$H_{6,10}\subset G_{6,10}$, that would be analogous to $H_3 \subset G_3$, 
and then for
these H's to identify the infinite dimensional simple quotients of 
$C^{*}(H)$,
both the faithful ones (generated by a faithful
representation of $H$) and the non-faithful ones, and also to give matrix
representations over lower dimensional algebras for as many of the 
non-faithful
quotients as possible.\\

The most attractive concrete representations on $L^2(\Bbb T)$ of the
irrational rotation algebras $A^3_{\theta}$ (for irrational $\theta$), the 
infinite
dimensional simple quotients of the group $C^{*}$-algebra $C^{*}(H_3)$, use 
the rotation
flow $(\Bbb Z,\Bbb T)$, $v \mapsto\lambda v$, on the circle $\Bbb T$
(with $\lambda =e^{2 \pi i \theta}$). The analogous 2, 3, and 4 dimensional 
(Anzai) flows,
generated by $(w,v)\mapsto (\lambda w,wv)$ on $\Bbb T^2$, $(x,w,v)\mapsto
(\lambda x,xw,wv)$ on $\Bbb T^3$, and $(y,x,w,v)\mapsto
(\lambda y,yx,xw,wv)$ on $\Bbb T^4$ give the analogous concrete 
representations
of the algebras $A^4_\theta$, $A^{5,5}_\theta$, and $A^{6,10}_\theta$ on
$L^2(\Bbb T^2)$, $L^2(\Bbb T^3)$, and $L^2(\Bbb T^4)$ of simple quotients of
$C^{*}(H_4)$ \cite{pM93}, $C^{*}(H_{5,5})$ \cite{pM97}, and 
$C^{*}(H_{6,10})$ \cite{pM02}.\\

The $K$-groups of $A^3_{\theta}=A_{\theta}$ were computed in \cite{mP80} and 
\cite{mR81},
and those of the Heisenberg $C^*$-algebra $A^4_{\theta}$ were studied by J. 
Packer in
\cite{jP88} (where it is referred to as \emph{class 2}). Since 
$A^{5,5}_\theta$
is isomorphic to a crossed product 
$\mathcal{C}(\Bbb{T}^3)\rtimes_{\sigma}\Bbb{Z}$
as in \cite{pM97}, S. Walters in \cite{sW02} used the Pimsner-Voiculescu six 
term exact sequence
\cite{mP80} to compute its $K$-groups, and in order to calculate the action 
of the
underlying automorphism $\sigma$ of the crossed product on 
$K_*(\mathcal{C}(\Bbb{T}^3))$
he made use of Connes' non-commutative geometry involving cyclic cocycles 
and the
Connes-Chern character. Application of this method to compute the $K$-groups 
of
$A^{6,10}_\theta$ and its higher dimensional analogues is much more 
difficult,
and no computations for these cases seem to have been attempted. 
In Section 6 this problem will be
discussed in a general form via another approach (looking at
$K_*(\mathcal{C}(\Bbb{T}^n))$ as an exterior algebra over $\Bbb{Z}^n$);
we have calculated the $K$-groups of $\mathcal{A}_{n,\theta}$ for 
$1 \leq n \leq 11$ (see Table 1).
\\

In the present paper, we consider the $n$-dimensional Anzai flow
$\mathcal{F}=(\mathbb{Z},\Bbb{T}^n)$ generated by the Anzai transformation
$$\sigma:(v_1,v_2,\ldots,v_n)\longmapsto(\lambda v_1,v_1 v_2,
\ldots,v_{n-1}v_n)$$
with $\lambda=e^{2 \pi i \theta}$ for an irrational $\theta$, and the 
associated operator
equations
\begin{equation*}
[U,V_1]=\lambda,\,\,\,[U,V_2]=V_1,\,\,\,\ldots,\,\,\,[U,V_n]=V_{n-1}. 
\tag*{${{\bahram}_n}$}
\end{equation*}
We show that the group $\mathfrak{D}_n$ ($\mathfrak{D}_1=H_3$, 
$\mathfrak{D}_2=H_4$,
$\mathfrak{D}_3=H_{5,5}$, and $\mathfrak{D}_4=H_{6,10}$), to which these are 
related, is a cocompact
subgroup of a connected \text{$(n+2)$-dimensional} group $F_n$ ($F_1=G_3$, 
$F_2=G_4$,
$F_3=G_{5,5}$, and $F_4=G_{6,10}$). In Theorems \ref{theo:1}
and \ref{theo:2} the infinite dimensional simple quotients of 
$C^*(\mathfrak{D}_n)$
are identified and displayed as $C^*$-crossed products generated by minimal 
actions.
The faithful ones using the flow $\mathcal{F}$ are called 
$\mathcal{A}_{n,\theta}$
($\mathcal{A}_{1,\theta}=A^3_\theta$, $\mathcal{A}_{2,\theta}=A^4_\theta$,
$\mathcal{A}_{3,\theta}=A^{5,5}_\theta$, and 
$\mathcal{A}_{4,\theta}=A^{6,10}_\theta$),
and the others using suitable modifications
of $\mathcal{F}$ are called $A^{(n)}_1,A^{(n)}_2,\ldots,A^{(n)}_{n-1}$. The 
non-faithful
ones $A^{(n)}_i$ are displayed also as matrix algebras over simple
\text{$C^*$-algebras} from groups of lower dimension (Theorem \ref{theo:4}); 
these groups
are certain subgroups of $\mathfrak{D}_{n-i}$ for $i=1,2,\ldots,n-1$.
These \text{$C^*$-algebras} are completely classified in Corollary 
\ref{coro:4} by $K$-theoretic invariants, namely the
ranks of the $K$-groups and the range of the unique tracial state. 
In Section 6, we make use of the algebraic properties of 
$K_*(\mathcal{C}(\Bbb{T}^n))$. More precisely, $K_*(\mathcal{C}(\Bbb{T}^n))$
is an exterior algebra over $\Bbb{Z}^n$ with a certain natural basis
, and the induced automorphism $\sigma_*$
is in fact a ring automorphism, which makes computations much easier. In 
fact, it is shown in Proposition \ref{prop:1} that the problem of finding 
the $K$-groups of any transformation group $C^*$-algebras of the tori 
is completely computable in the sense that one only needs to calculate 
the kernels and cokernels of a finite number of integer matrices.\\

We have learned that a similar method was used by R. Ji in 
his Ph.D. thesis \cite[unpublished]{rJ86} to study 
the $K$-groups of the $C^*$-algebras $A_{F_{f,\theta}}$
associated with the descending Furstenberg transformations $F_{f,\theta}$ 
on tori. He has computed a 2-dimensional case \cite[Corollary 2.20]{rJ86}
(which includes $A^4_\theta$), and has made a claim 
\cite[Proposition 2.17]{rJ86} 
about the form of the torsion subgroup of 
$K_*(A_{F_{f,\theta}})$ that is not true in general 
(see Remark \ref{rema:2} below), and does not deal with the rank of the
$K$-groups. Thus the $K$-groups of $A_{F_{f,\theta}}$ (including 
$\mathcal{A}_{n,\theta}$) have not yet been calculated in general.
In Corollary \ref{coro:3}, the $K$-groups of any transformation 
group $C^*$-algebra of an $n$-torus are proved to be finitely 
generated with the same rank. In the case of $\mathcal{A}_{n,\theta}$, 
this common rank is called $a_n$, and is studied in detail in 
Section 6, and given as explicitly as possible via generating 
functions (Theorem \ref{theo:6}). It is proved in 
Theorem \ref{theo:5} that $a_n$ is the common rank of the $K$-groups
of many other 
transformation group $C^*$-algebra of $\Bbb{T}^n$ as well, including 
the $C^*$-algebras from Furstenberg transformations on $\Bbb{T}^n$. An 
explicit formula for the torsion parts of the $K$-groups of these 
algebras seems much more challenging to find.
\\

To present the results and proofs of the paper, especially for Section 6, we 
need some
definitions about transformations on the tori and the corresponding 
\text{$C^*$-crossed products.}\\
Let $\mathbb{T}^n$ denote the $n$-dimensional torus with coordinates 
$(v_1,v_2,\ldots,v_n)$.
We consider homotopy classes of continuous functions from 
$\mathbb{T}^n$ to $\Bbb{T}$.
It is known that in each class there is a unique ``linear" function 
$f(v_1,\ldots,v_n)=
v_1^{b_1} v_2^{b_2}\ldots v_n^{b_n}$, $b_1,b_2,\ldots,b_n \in\mathbb{Z}$. 
Following
\cite[p. 35]{hF81}, we denote the exponent $b_i$, which is uniquely determined by 
the
homotopy class of $f$, as $b_i=A_i[f]$.
\begin{defi}
An \textbf{affine transformation} on $\Bbb{T}^n$ is given by
\begin{multline*}
\alpha(v_1,v_2,\ldots,v_n)=\\(e^{2 \pi i t_1} v_1^{b_{11}}\ldots 
v_n^{b_{n1}},
e^{2 \pi i t_2} v_1^{b_{12}}\ldots v_n^{b_{n2}},\ldots,
e^{2 \pi i t_n} v_1^{b_{1n}}\ldots v_n^{b_{nn}}),
\end{multline*}
where $\boldsymbol{t}:=(t_1,t_2,\ldots,t_n)\in\Bbb{R}^n$ and 
$\mathsf{A}:=[b_{ij}]_{n \times n}\in
\mathrm{GL}(n,\Bbb{Z})$. We identify the pair $(\boldsymbol{t},\mathsf{A})$ 
with
$\alpha$.
\end{defi}

Note that any automorphism of $\Bbb{T}^n$ followed by a rotation can be 
expressed in such
a fashion. The set of affine transformations on $\Bbb{T}^n$ form a group 
$\mathrm{Aff}(\Bbb{T}^n)$,
which can
be identified with the semidirect product $\Bbb{R}^n 
\rtimes\mathrm{GL}(n,\Bbb{Z})$.
More precisely, for two affine transformations 
$\alpha=(\boldsymbol{t},\mathsf{A})$
and $\alpha'=(\boldsymbol{t}',\mathsf{A}')$ on $\Bbb{T}^n$, we have
$$\alpha\circ\alpha'=(\boldsymbol{t}+\mathsf{A}\boldsymbol{t}',
\mathsf{A}\mathsf{A}')
~~\text{and}~~\alpha^{-1}=(-\mathsf{A}^{-1}\boldsymbol{t},\mathsf{A}^{-1}).$$
(In the expression $\mathsf{A}\boldsymbol{t}$, $\boldsymbol{t}$ is a  column 
vector, but for convenience we write it as a row vector.)\\

We have the following definitions in accordance with \cite{rJ86}. These 
transformations
of the tori appear in applications of ergodic theory to number theory 
\cite{hF81}, and sometimes
are called \emph{skew product transformations of the tori}.

\begin{defi}
\begin{itemize}
\item [\textnormal{(a)}]
A \textnormal{\textbf{Furstenberg transformation}} $F_{f,\theta}$ on 
$\Bbb{T}^n$ is given by
\begin{multline*}
F_{f,\theta}(v_1,v_2,\ldots,v_n)=(e^{2 \pi i \theta} 
v_1,f_1(v_1)v_2,f_2(v_1,v_2)v_3,
\ldots,\\ f_{n-1}(v_1,\ldots,v_{n-1})v_n),
\end{multline*}
where $\theta$ is a real number and each
$f_i:\Bbb{T}^i \rightarrow\Bbb{T}$ is a continuous function with
$A_i[f_i]\neq 0$ for $i=1,\ldots,n-1$.
\item [\textnormal{(b)}]
An \textnormal{\textbf{affine Furstenberg transformation}} on $\Bbb{T}^n$ is 
given by
$$\alpha(v_1,v_2,\ldots,v_n)=(e^{2 \pi i \theta} v_1,v_1^{b_{12}}v_2,
v_1^{b_{13}}v_2^{b_{23}} v_3,\ldots
,v_1^{b_{1n}}v_2^{b_{2n}}\ldots v_{n-1}^{b_{n-1,n}}v_n),$$
where $\theta$ is a real number and the exponents $b_{ij}$
are integers and $b_{i,i+1}\neq 0$ for $i=1,\ldots,n-1$.
\item [\textnormal{(c)}]
An \textnormal{\textbf{ascending Furstenberg transformation}} on $\Bbb{T}^n$ 
is given by
$$\alpha(v_1,v_2,\ldots,v_n)=(e^{2 \pi i \theta} v_1,v_1^{k_1} v_2,
v_2^{k_2} v_3,\ldots
, v_{n-1}^{k_{n-1}} v_n),$$
where $\theta$ is a real number and the exponents $k_i$ are nonzero
integers and $k_i \mid k_{i+1}$ for $i=1,\ldots,n-2$.
\item [\textnormal{(d)}]
In \textnormal{(c)}, if $k_i=1$ for $i=1,\ldots,n-1$, the transformation is 
called
the \textnormal{\textbf{Anzai transformation}} on $\Bbb{T}^n$. Thus it is 
given by
$$\sigma(v_1,v_2,\ldots,v_n)=(e^{2 \pi i \theta} v_1,v_1 v_2,\ldots,v_{n-1} 
v_n),$$
where $\theta$ is a real number.
\end{itemize}
\end{defi}

Note that one can easily verify that $F_{f,\theta}$ is a
homeomorphism. Also, in the above definition, we have converted
``descending", which is used in \cite[Definition 2.16]{rJ86},
to ``ascending" since the order of coordinates there is opposite to ours.
For certain Furstenberg transformations on $\Bbb{T}^n$, we have the 
following
theorem.

\begin{theo}
\textnormal{(\cite{hF61}, 2.3)}
If $\theta$ is irrational, then $F_{f,\theta}$ defines a minimal dynamical
system on $\Bbb{T}^n$. If in addition, each $f_i$ satisfies a uniform
Lipschitz condition in $v_i$ for $i=1,\ldots,n-1$, then $F_{f,\theta}$ is a 
uniquely ergodic
transformation and the unique invariant measure is the normalized Lebesgue
measure on $\Bbb{T}^n$. In particular, every affine Furstenberg 
transformation
defines a minimal and uniquely ergodic dynamical system if $\theta$ is 
irrational.
\end{theo}

As a conclusion, we have the following result for the Furstenberg
transformation group $C^*$-algebra
$A_{F_{f,\theta}}:=\mathcal{C}(\Bbb{T}^n)\rtimes_{F_{f,\theta}}\mathbb{Z}$
as introduced in \cite{rJ86}.

\begin{coro}\label{coro:1}
$A_{F_{f,\theta}}=\mathcal{C}(\Bbb{T}^n)\rtimes_{F_{f,\theta}}\mathbb{Z}$ is 
a simple
$C^*$-algebra for irrational $\theta$. If in addition, each $f_i$ satisfies 
a uniform
Lipschitz condition in $v_i$ for $i=1,\ldots,n-1$, then $A_{F_{f,\theta}}$ has 
a
unique tracial state.
\end{coro}
\begin{proof}
For the first part, the minimality of the action as stated in
the preceding theorem implies the simplicity of $A_{F_{f,\theta}}$
\cite{eE67,scP78}. For the second part, one can
easily check that since $\theta$ is irrational, the action of
$\mathbb{Z}$ on $\Bbb{T}^n$ generated by $F_{f,\theta}$ is
free. So, there are no periodic points in $\Bbb{T}^n$. This and
the unique ergodicity of $F_{f,\theta}$ yield the result
\cite[Corollary 3.3.10, p. 91]{jT87}.
\end{proof}
\begin{rema}\label{rema:1} Using the preceding corollary and 
much like the
proof of Theorem \ref{theo:1}, one can prove that for irrational $\theta$, 
$A_{F_{f,\theta}}$ is
in fact the unique $C^*$-algebra generated by unitaries $U,V_1,\ldots,V_n$ 
satisfying
the commutator relations
\begin{equation*}
[U,V_1]=e^{2 \pi i \theta},[U,V_2]=f_1(V_1),\ldots,
[U,V_n]=f_{n-1}(V_1,\ldots,V_{n-1}) \tag*{${{\bahram}_{f}}$}
\end{equation*}
(where $[a,b]:=aba^{-1}b^{-1}$ and all other pairs of operators
from $U,V_1,\ldots,V_n$ commute).
\end{rema}

\begin{rema}\label{rema:2} In \cite[Proposition 2.17]{rJ86}, R. 
Ji claims to have proved\\
\linebreak
\textnormal{($\ast$)} \emph{If $F_{f,\theta}$ is an ascending Furstenberg 
transformation
on $\Bbb{T}^n$ with the ascending sequence $\{k_1,k_2,\ldots,k_{n-1}\}$, 
then the
torsion subgroup of $K_*(A_{F_{f,\theta}})$ is isomorphic to
$\mathbb{Z}_{k_1}\oplus\mathbb{Z}_{k_2}^{(m_2)}\oplus\ldots\oplus
\mathbb{Z}_{k_{n-1}}^{(m_{n-1})}$ ,where the group 
$\mathbb{Z}_{k_i}^{(m_i)}$ is the direct
product of $m_i$ copies of the cyclic group $\mathbb{Z}_{k_i}=\mathbb{Z}/k_i 
\mathbb{Z}$.}\\
\linebreak
From this claim one can immediately deduce that the $K$-groups 
of the
\text{$C^*$-algebra} 
$\mathcal{A}_{n,\theta}:=\mathcal{C}(\Bbb{T}^n)\rtimes_\sigma \mathbb{Z}$
generated by the Anzai transformation $\sigma$ on $\Bbb{T}^n$ are 
torsion-free. We will
show in the last section that this is not true in general. As the first 
counterexample,
we will see that 
$K_1(\mathcal{A}_{6,\theta})\cong\mathbb{Z}^{13}\oplus\mathbb{Z}_2$
(Example 6.2).
In fact, the error in the proof of ($\ast$) is in \cite[p. 29, l. 2]{rJ86};
there it is ``clearly'' assumed that using a matrix $S$ in 
$\mathrm{GL}(2^n,\Bbb{Z})$, one can delete all entries
denoted by $\star$'s in $\mathsf{K_*-I}$, where $\mathsf{K_*}$ is the 
$2^n \times 2^n$ 
integer matrix corresponding to $A_{F_{f,\theta}}$ that acts 
on $K_*(\mathcal{C}(\Bbb{T}^n))=\Lambda^*\Bbb{Z}^n$ with respect to
a certain ordered basis. This error arose originally from the general 
form of the matrix $\mathsf{K}_*$ in \cite[p. 27]{rJ86}, which is not 
correct.
R. Ji went on to use the torsion subgroup in ($\ast$) as an invariant
to classify the $C^*$-algebras generated by ascending transformations and 
matrix
algebras over them \cite[Theorem 3.6]{rJ86}. We do not know whether that
classification holds.
\end{rema}

\section{The Anzai flow $\mathcal{F}=(\mathbb{Z},\Bbb{T}^{\hspace{1pt}n})$
and the group $\mathfrak{D}_n$}

Let $ \lambda:=e^{2 \pi i \theta}$ for an irrational number $\theta$,
and consider the Anzai transformation $$\sigma:(
v_1,v_2,\ldots,v_n)\longmapsto(\lambda v_1,v_1 v_2,\ldots,v_{n-1} v_n)$$
on $\Bbb{T}^n$, which generates
(by iteration)
the Anzai flow $\mathcal{F}=(\mathbb{Z},\Bbb{T}^n)$
\begin{multline*}
(v_1,\ldots,v_n)\overset{m}{\longmapsto}\sigma^{m}(v_1,\ldots,v_n)=\\
(\lambda^m v_1,\lambda^{\binom{m}{2}} 
v_1^{m}v_2,\lambda^{\binom{m}{3}}v_1^{\binom{m}{2}}v_2^m v_3
,\ldots,\lambda^{\binom{m}{n}}v_1^{\binom{m}{n-1}}\ldots v_{n-1}^m v_n).
\end{multline*}
With $v_1,\ldots,v_n$ denoting also the functions
in $\mathcal{C}(\Bbb{T}^n)$, $$(v_1,\ldots,v_n)
\longmapsto v_1,\ldots,v_n,~\text{respectively,}$$ we then get unitaries on
$L^2(\Bbb{T}^n)$, $$Uf=f \circ \sigma
,\,\,\,V_1 f
=v_1 f,\,\,\,\ldots,\,\,\,V_n f=v_n f.$$
These unitaries satisfy the commutator equations
\begin{equation*}
[U,V_1]=\lambda,\,\,\,[U,V_2]=V_1,\,\,\,\ldots,\,\,\,[U,V_n]=V_{n-1}, 
\tag*{${{\bahram}_n}$}
\end{equation*}
(where $[a,b]:=aba^{-1}b^{-1}$ and all other pairs of operators
from $U,V_1,\ldots,V_n$ commute).
A ``discrete group construction" \cite{pM93} shows how to construct a group 
from unitaries like this; use ${\bahram}_n$ 
to collect terms in the
product$$(\lambda^{k_1}V_1^{k_2}V_2^{k_3}\ldots
V_n^{k_{n+1}}U^k)(\lambda^{k^\prime_1}V_1^
{k^\prime_2}V_2^{k^\prime_3}\ldots
V_n^{k^\prime_{n+1}}U^{k^\prime})$$then, the exponents give the
multiplication for a group $\mathfrak{D}_n$$(=\mathbb{Z}^{n+1}
\times \mathbb{Z}$, as a set$)$. In fact
$\mathfrak{D}_n=\mathbb{Z}^{n+1} \rtimes_\eta \mathbb{Z}$ is a semidirect
product for which $\eta:\mathbb{Z}\rightarrow
\mathrm{GL}(n+1,\mathbb{Z})$ is given by
$\eta(k)=\eta_k=\mathsf{M}_n^k$, where $\mathsf{M}_n$ is a 
$(n+1)\times(n+1)$ matrix
defined as
\[\mathsf{M}_n = \left[\begin{array}{ccccc}
1 & 1 & 0 & \cdots& 0\\
0 & 1 & 1 & &\vdots\\
0 &\ddots&\ddots&\ddots&0\\
\vdots&   & 0 & 1 & 1\\
0 &\cdots  & 0 & 0 & 1
\end{array}
\right]_{(n+1)\times(n+1)},
\]
and (by induction) one can show that $\mathsf{M}_n^k=[m_{ij}^{(k)}]$, where
\[
m_{ij}^{(k)} =
\binom{k}{j-i}
\]
with the following notation.\\

\begin{nota}\label{nota:1}
\begin{equation*}
\binom{k}{r} :=
\begin{cases}
\frac{k(k-1)\ldots (k-r+1)}{r!},  &\text{if $0 \leq r < k$ or $(k<0$ and $ r 
 >0)$;}\\
1 , &\text{if $r = (k+|k|)/2$;}\\
0 , &\text{otherwise.}
\end{cases}
\end{equation*}
\end{nota}
~\\

Note that in this construction, we identify an element 
$((k_1,\ldots,k_{n+1}),k)$ of $\mathfrak{D}_n$
with $\lambda^{k_1}V_1^{k_2}V_2^{k_3}\ldots V_n^{k_{n+1}}U^k$ and the 
multiplication
of $\mathfrak{D}_n$ is given by
\begin{multline*}
((k_1,\ldots,k_{n+1}),k).((k^\prime_1,\ldots,k^\prime_{n+1}),k^\prime)=\\
((k_1,\ldots,k_{n+1})+\eta_k(k^\prime_1,\ldots,k^\prime_{n+1}),k+k^\prime).
\end{multline*}
One can easily check that $\mathfrak{D}_n$ is the discrete group generated 
by
$x,y_0,y_1,\ldots,\linebreak y_n$ such that $x y_0=y_0x$ and $y_i y_j=y_j y_i$ for $0 
\leq i,j \leq  n$ and
\begin{equation}
[x,y_1]=y_0,\,\,\,[x,y_2]=y_1,\,\,\,\ldots,\,\,\,[x,y_n]=y_{n-1}
\end{equation}
The group $\mathfrak{D}_n$ is discrete, nilpotent, finitely generated and 
torsion-free.
So by a result of Malcev \cite{aM49}, it follows that $\mathfrak{D}_n$ is 
isomorphic
to a discrete cocompact subgroup of a connected $(n+2)$-dimensional 
nilpotent Lie group
${F}_n$.
More precisely, one can verify that the corresponding Lie algebra is the so 
called
\textbf{generic filiform Lie algebra} \cite{vG97,mG96} $\mathfrak{f}_{n}$, 
which is spanned
by $X,Y_0,Y_1,\ldots,Y_n$ with non-zero brackets
\begin{equation}
[X,Y_1]=Y_0,\,\,\,[X,Y_2]=Y_1,\,\,\,\ldots,\,\,\,[X,Y_n]=Y_{n-1}.
\end{equation}
Let $\mathcal{A}_{n,\theta}$ denote the $C^{*}$-algebra generated by the 
unitaries $U,V_1,\ldots,V_n$. Note
that $U,V_n$ generate $\mathcal{A}_{n,\theta}$; the unitaries 
$V_1,\ldots,V_n$ have been introduced only
to control the notation. An obvious property of the construction is 
that$$\pi:((k_1,\ldots,k_{n+1}),k)
\longmapsto \lambda^{k_1}V_1^{k_2}V_2^{k_3}\ldots V_n^{k_{n+1}}U^k$$ is a 
unitary representation of
$\mathfrak{D}_n$ on $L^2(\Bbb{T}^n)$ that generates 
$\mathcal{A}_{n,\theta}$.

\section{The faithful simple quotients $\mathcal{A}_{n,\theta}$ of 
$C^{*}(\mathfrak{D}_n)$}

Let $ \lambda:=e^{2 \pi i \theta}$ for an irrational number $\theta$. Since
$\mathcal{A}_{n,\theta}$ is generated by a representation of 
$\mathfrak{D}_n$,
it is a quotient of $C^{*}(\mathfrak{D}_n)$.

\begin{theo}\label{theo:1}
Let $ \lambda:=e^{2 \pi i \theta}$ for an irrational number $\theta$.
  \begin{itemize}

\item [\textnormal{(a)}]

There is a unique (up to isomorphism) $C^{*}$-algebra  
$\mathcal{A}_{n,\theta}$ generated by unitaries
$U,V_1,\ldots,V_n$ satisfying \textnormal{${\bahram}_n$}; 
$\mathcal{A}_{n,\theta}$ is simple and is universal for the
equations \textnormal{${\bahram}_n$}. Let $\sigma$ be the homeomorphism 
used in the definition of
$\mathcal{A}_{n,\theta}$.
  Then $$\mathcal{A}_{n,\theta}\cong 
\mathcal{C}(\Bbb{T}^n)\rtimes_{\sigma}\mathbb{Z}.$$

\item[\textnormal{(b)}] Let $\pi^\prime$ be a representation of 
$\mathfrak{D}_n$ such that $\pi=\pi^\prime$
(as scalars)
on the
  center $((\mathbb{Z},0,\ldots,0),0)$ of $\mathfrak{D}_n$, and let 
$\mathcal{A}$ be a  $C^{*}$-algebra
  generated by $\pi^\prime$. Then $\mathcal{A}\cong \mathcal{A}_{n,\theta}$ 
via a unique isomorphism $
  \omega:\mathcal{A}_{n,\theta}\rightarrow \mathcal{A}$ such that $\omega\circ 
\pi=\pi^\prime.$

  \item[\textnormal{(c)}] The $C^{*}$-algebra  $\mathcal{A}_{n,\theta}$ has 
a unique tracial state.

\item[\textnormal{(d)}] There is an automorphism $\alpha$ of 
$\mathcal{A}_{n-1,\theta}$ such that
   $\mathcal{A}_{n,\theta}\cong 
\mathcal{A}_{n-1,\theta}\rtimes_\alpha\mathbb{Z}$.
   \end{itemize}
\end{theo}
\begin{proof}
(a) The flow $\mathcal{F}$ used in the definition of 
$\mathcal{A}_{n,\theta}$ is minimal \cite{hF61}
, so the $C^{*}$-crossed product 
$\mathcal{C}(\Bbb{T}^n)\rtimes_{\sigma}\mathbb{Z}$ is simple
\cite{eE67,scP78}.
On the other hand,  $\mathcal{A}_{n,\theta}$ provides a covariant 
representation of the dynamical system
$(\mathcal{C}(\Bbb{T}^n),\mathbb{Z} ,\sigma)$. More precisely, let 
$\varphi:k \mapsto U^k$; $\varphi(k)(g)=g
\circ \sigma^{k}$ be the unitary representation
of $\mathbb{Z}$ on $L^2(\Bbb{T}^n)$ 
and $\nu:f \mapsto M_f$; $M_f(g)=fg$ be the $*$-representation
of $\mathcal{C}(\Bbb{T}^n)$ on $L^2(\Bbb{T}^n)$ $(g \in L^2(\Bbb{T}^n))$.
Then it is easy to see that $(\nu,\varphi)$ is a covariant pair for 
$(\mathcal{C}(\Bbb{T}^n),\mathbb{Z} ,\sigma)$. Hence
by the universal property of 
$\mathcal{C}(\Bbb{T}^n)\rtimes_{\sigma}\mathbb{Z}$, there is a 
$*$-homomorphism
$\rho:\mathcal{C}(\Bbb{T}^n)\rtimes_{\sigma}\mathbb{Z}\rightarrow 
C^{*}(\nu(\mathcal{C}(\Bbb{T}^n)),\varphi(\mathbb{Z}))$
obtained by setting $\rho(\sum_k f_k u^k)=\sum_k M_{f_k}U^k$ ($u$ is the 
unitary implementing the automorphism
$\sigma$) on the dense $*$-subalgebra $\mathcal{C}(\Bbb{T}^n)
\mathbb{Z}$ and extending by continuity on 
$\mathcal{C}(\Bbb{T}^n)\rtimes_{\sigma}\mathbb{Z}$. $\rho$ is
surjective since $\mathcal{C}(\Bbb{T}^n)$ is unital and 
$C^{*}(\nu(\mathcal{C}(\Bbb{T}^n),\varphi(\mathbb{Z}))$
is the $C^{*}$-subalgebra
of $\mathcal{B}(L^2(\Bbb{T}^n))$ generated by $U$ and the multiplication 
operators $M_f$ coming from
$\mathcal{C}(\Bbb{T}^n)$. But, as $\mathcal{C}(\Bbb{T}^n)$ is generated by 
the coordinate functions $v_1,\ldots,v_n$,
$\nu(\mathcal{C}(\Bbb{T}^n))$ is generated by special unitary operators 
$\nu(v_i)=M_{v_i}=V_i$ for $i=1,\ldots,n$.
So, 
$C^{*}(\nu(\mathcal{C}(\Bbb{T}^n)),u(\mathbb{Z}))=\mathcal{A}_{n,\theta}$. 
We now have a surjective $*$-
homomorphism 
$\rho:\mathcal{C}(\Bbb{T}^n)\rtimes_{\sigma}\mathbb{Z}
\rightarrow\mathcal{A}_{n,\theta}$ 
with
$\mathcal{C}(\Bbb{T}^n)\rtimes_{\sigma}\mathbb{Z}$ simple. 
Therefore $\rho$ is an isomorphism and $\mathcal{A}_{n,\theta}$ 
is simple too.

Now, let $\mathcal{A^\prime}$ be another $C^{*}$-algebra generated by 
unitaries $U^\prime,V_1^\prime,\ldots,V_n^\prime$ satisfying  
${\bahram}_n$. Since $V_1^\prime,\ldots,V_n^\prime$ commute, there is a
$*$-homomorphism $\mu:\mathcal{C}(\Bbb{T}^n)\rightarrow\mathcal{A}^\prime$ 
such that $\mu(v_i)=V_i^\prime$ for
   $i=1,\ldots,n$. In fact $\mu(f)=f(V_1^\prime,\ldots,V_n^\prime)$. Let 
$\widetilde{\omega}:\mathbb{Z}
\rightarrow \mathcal{A}^\prime$ be the unitary representation 
$\widetilde{\omega}(k)={U^\prime}^k$. Noting
that $\mu(f \circ \sigma^{k})=\widetilde{\omega}(k) \mu(f) 
\widetilde{\omega}(k)^{*}$ holds for $f=v_1,\ldots,
v_n$ and hence for all $f \in \mathcal{C}(\Bbb{T}^n)$; by the universal 
property of $\mathcal{C}(\Bbb{T}^n)\rtimes_{\sigma}
\mathbb{Z}$, the covariant pair $(\mu,\widetilde{\omega})$ yields a 
homomorphism of
$\mathcal{C}(\Bbb{T}^n)\rtimes_{\sigma}\mathbb{Z}$ onto $\mathcal{A}^\prime$ 
mapping $v_i$ to $V_i^\prime$ for $i=1,\ldots,n$,
and $u$ to $U^\prime$. So, $\mathcal{A}_{n,\theta}$ is universal for 
equations  ${\bahram}_n$.

(b) The hypothesis imply that  ${\bahram}_n$ is satisfied by the unitaries 
$U^\prime,$ $V_1^\prime,$ $\ldots
,$ $V_n^\prime$ given by 
$$\pi^\prime((k_1,\ldots,k_{n+1}),k)=
\lambda^{k_1}{V_1^\prime}^{k_2}{V_2^\prime}^
{k_3}\ldots {V_n^\prime}^{k_{n+1}}{U^\prime}^k.$$Part (a) and its proof now 
yields the result.

(c) This flow is minimal and uniquely ergodic with respect to (normalized) 
Haar measure $\lambda
$ on $\Bbb{T}^n$ \cite{hF61}. So, $\mathcal{A}_{n,\theta}$ has a unique 
tracial state $\tau$ given by
$\tau(\sum_k f_k u^k)=\int f_{\hspace{-1pt}_{_0}} d \lambda$
\cite[Corollary VIII.3.8, p. 91]{krD96}.

(d) Let $\mathcal{A}_{n-1,\theta}$ be generated by operators 
$U^\prime,V_1^\prime,\ldots,V_{n-1}^\prime$
satisfying \linebreak ${\bahram}_{n-1}$. Define $\alpha$ as 
$\alpha(U^\prime)={V^\prime}^{-1}_{n-1}U^\prime$, $\alpha
(V_i^\prime)=V_i^\prime$ for $i=1,\ldots,n-1$. Since $\alpha(U^\prime)$, 
$\alpha(V^\prime_1)$,\ldots,
$\alpha(V^\prime_{n-1})$ also satisfy ${\bahram}_{n-1}$, $\alpha$ can be 
extended to an
automorphism of $\mathcal{A}_{n-1,\theta}$. Define a $*$-homomorphism
$p:\mathcal{A}_{n-1,\theta}\rightarrow \mathcal{A}_{n,\theta}$ by 
$p(U^\prime)=U$ and $p(V^\prime_i)=
V_i$ for $i=1,\ldots,n-1$, and a unitary representation 
$\psi:\mathbb{Z}\rightarrow \mathcal{A}_{n,\theta}$ by
$\psi(k)=V_n^k$. $(p,\psi)$ is a covariant pair for 
$(\mathcal{A}_{n-1,\theta},\mathbb{Z},\alpha)$ since
the equality $p(\alpha^{(k)})(a)={V^k_n}\pi(a){V^{-k}_n}$ holds for 
$a=U^\prime,V^\prime_1,\ldots,
V^\prime_{n-1}$; hence for all $a \in\mathcal{A}_{n-1,\theta}$. So, there is 
a $*$-homomorphism from
$\mathcal{A}_{n-1,\theta}\rtimes_\alpha\mathbb{Z}$  onto 
$\mathcal{A}_{n,\theta}$ mapping $U^\prime$
to $U$ and $V^\prime_i$ to $V_i$ for $i=1,\ldots,n-1$, and $V^\prime_n$ (the 
unitary implementing the
automorphism $\alpha$) to $V_n$.

Conversely, $\mathcal{A}_{n-1,\theta}\rtimes_\alpha\mathbb{Z}$ is generated 
by the unitaries $U^\prime,
V^\prime_1,\ldots,V^\prime_{n-1}$ and $V^\prime_n$ satisfying  
${\bahram}_n$. So, by universality of
$\mathcal{A}_{n,\theta}$, there is a $*$-homomorphism from  
$\mathcal{A}_{n,\theta}$ onto
$\mathcal{A}_{n-1,\theta}\rtimes_\alpha\mathbb{Z}$ mapping $U$ to $U^\prime$ 
and $V_i$ to $V^\prime_i$
for $i=1,\ldots,n$. Clearly, these two $*$-homomorphisms are inverses of 
each other.
  \end{proof}

  \section{Non-faithful simple quotients of $C^{*}(\mathfrak{D}_n)$}

  When $\lambda=e^{2 \pi i \theta}$ for an irrational number $\theta$, 
$\mathcal{A}_{n,\theta}$ is a
  simple quotient of $C^{*}(\mathfrak{D}_n)$ and the representation 
$$\pi:((k_1,\ldots,k_{n+1}),k)
\longmapsto \lambda^{k_1}V_1^{k_2}V_2^{k_3}\ldots 
V_n^{k_{n+1}}U^k;~~~\mathfrak{D}_n
\longrightarrow \mathcal{A}_{n,\theta}$$ is faithful. But there are other 
infinite dimensional simple
quotients of $C^{*}(\mathfrak{D}_n)$; for them $\pi$ is not faithful.

Suppose that $\lambda$ is a primitive $q_1$-th root of unity and that $A$ is 
a simple quotient of
$C^{*}(\mathfrak{D}_n)$ that is irreducibly represented and generated by 
unitaries $U,V_1,\ldots,V_n$
satisfying  ${\bahram}_n$. Then $V_1^{q_1}$ commutes with $U$ and $V_i$ 
for $i=1,\ldots,n$ and so by irreducibility
equals ${\tilde\mu}_1 I$, a multiple of identity. Take $V_1=\mu_1 
V^\prime_1$ for $\mu_1^{q_1}={\tilde\mu}_1$
, so that ${V_1^\prime}^{q_1}=1$, and substitute $V_1=\mu_1 V^\prime_1$ in  
${\bahram}_n$ to get
\begin{equation*}
\begin{cases}
[U,V^\prime_1]=\lambda,~~[U,V_2]=\mu_1 V^\prime_1,\\
[U,V_3]=V_2,~~\ldots,~~[U,V_n]=V_{n-1},\tag*{${{\bahram}_{n,1}}$}\\
{V^\prime_1}^{q_1}=I
\end{cases}
\end{equation*}

1. If $\mu_1$ is not a root of unity, then we can modify the presentation 
$\mathcal{C}(\Bbb{T}^n)\rtimes_{\sigma}\mathbb{Z}$
for $\mathcal{A}_{n,\theta}$ in Theorem \ref{theo:1} and present the 
operators $U,V^\prime_1,V_2,\ldots,\linebreak V_n$, and their
generated algebra $A^{(n)}_1$, with the flow 
$\mathcal{F}_1=(\mathbb{Z},\mathbb{Z}_{q_1}\times\Bbb{T}^{n-1})$ generated
by the homeomorphism $\phi_1$ of 
$\mathbf{X}_1:=\mathbb{Z}_{q_1}\times\Bbb{T}^{n-1}$,

$$\phi_1(v_1,v_2,\ldots,v_n)=(\lambda v_1,\mu_1 v_1 v_2,v_2 
v_3,\ldots,v_{n-1} v_n).$$

To see that $\mathcal{F}_1$ is minimal, we use a lemma on minimality of
skew products of dynamical systems with $\Bbb{T}$ \cite[2.1, 2.3]{hF61}.
In fact, let $(X,\phi)$ be a dynamical system, where $X$ is a compact metric
space and $\phi$ is a homeomorphism of $X$. Let $g:X \rightarrow\Bbb{T}$
be a continuous function and consider the dynamical system $(X 
\times\Bbb{T},\Phi)$
defined by $\Phi(x,\zeta)=(\phi(x),g(x)\zeta)$, which is called a \emph{skew 
product} of
$(X,\phi)$ with $\Bbb{T}$. Then we have the following lemma.

\begin{lemm}\textnormal{(\cite[2.1, 2.3]{hF61})}
Let $(X,\phi)$ be a minimal dynamical system and $g:X \rightarrow\Bbb{T}$ be 
a continuous
function. Consider the skew product dynamical system
$(X \times\Bbb{T},\Phi)$ as above. Then $(X \times\Bbb{T},\Phi)$ is minimal 
if, and only if,
  for any non-zero integer $k$, the functional equation
\begin{equation}
g^k=\frac{R \circ\phi}{R} \tag{$\triangle$}
  \end{equation}
  has no continuous solution $R:X \rightarrow\Bbb{T}$.
\end{lemm}

Using the preceding lemma, we can prove the minimality of $\mathcal{F}_1$ as 
follows
by induction on $n$. For $n=2$, we use the preceding lemma, although 
\cite[Theorem 3]{pM93}
yields the result. In this case, $X=\Bbb{Z}_{q_1}$, $\phi(v_1)=\lambda v_1$ 
and
$g(v_1)=\mu_1 v_1$ for $v_1 \in\Bbb{Z}_{q_1}$. $\phi$ is clearly minimal and 
suppose that
$(\triangle)$ has a solution $R:\Bbb{Z}_{q_1} \rightarrow\Bbb{T}$ for some 
$0 \ne k \in\Bbb{Z}$.
Then one can easily obtain $R(\lambda^j)=R(1)\mu_1^{kj}\lambda^{k 
\binom{j}{2}}$, which
yields the contradiction $\mu_1^{k q_1}\lambda^{k \binom{q_1}{2}}=1$ since 
$\mu_1$ is not
a root of unity but $\lambda$ is a root of unity. Thus $\mathcal{F}_1$ is 
minimal.\\
Now, suppose that $\mathcal{F}_1$ is minimal for $n=m-1$. In this case
$X=\mathbb{Z}_{q_1}\times\Bbb{T}^{m-2}$, $\phi(v_1,v_2,\ldots,v_{m-1})=
(\lambda v_1,\mu_1 v_1 v_2,v_2 v_3,\ldots,v_{m-2} v_{m-1})$ and
$g(v_1,\ldots,v_{m-1})=v_{m-1}$. If $(\triangle)$ has a solution
$R:X \rightarrow\Bbb{T}$ for some $0 \ne k \in\Bbb{Z}$, then we have
$$v_{m-1}^k=R(\lambda v_1,\mu_1 v_1 v_2,v_2 v_3,\ldots,v_{m-2} v_{m-1})/
R(v_1,v_2,\ldots,v_{m-1})$$ for all $(v_1,v_2,\ldots,v_{m-1})\in X$.
But this equality is impossible, for \linebreak $R(v_1,\ldots,v_{m-1})$ would 
have a
certain degree in $v_{m-1}$ and one verifies that
$R(\lambda v_1,\mu_1 v_1 v_2,v_2 v_3,\ldots,v_{m-2} v_{m-1})$ has the same 
degree in
$v_{m-1}$, hence the right side of the above equality would have degree 0 
but the
left side has degree $k \ne 0$. So the skew product system $(X 
\times\Bbb{T},\Phi)=
(\mathbb{Z}_{q_1}\times\Bbb{T}^{m-1},\phi_1)$ is minimal. \\

Now take $\mathbf{Y}_1=\mathbb{Z}_{q_1}$. So the $C^{*}$-crossed product
$\mathcal{C}(\mathbf{Y}_1 \times\Bbb{T}^{n-1})\rtimes_{\phi_1}\mathbb{Z}$ is 
simple and
isomorphic to $A^{(n)}_1$.\\

2. Suppose that $\mu_1$ is also a root of unity, say a primitive $p_1$-th 
root of unity, and let
$q_2=\mathrm{lcm}\{q_1,p_1 \}$, the least common multiple of $q_1$ and $p_1$. 
Then $V_2^{q_2}={\tilde\mu}_2 I$,
  a multiple of the identity. If ${\tilde\mu}_2$ is not a root of unity, 
substitute $V_2=\mu_2 V^\prime_2$
  (as well as $V_1=\mu_1 V^\prime_1$) in  ${\bahram}_n$, where 
$\mu_2^{q_2}={\tilde\mu}_2$, and get
\begin{equation*}
\begin{cases}
[U,V^\prime_1]=\lambda,~~[U,V^\prime_2]=\mu_1 V^\prime_1,~~[U,V_3]=\mu_2 
V^\prime_2,\\
[U,V_4]=V_3,~~\ldots,~~[U,V_n]=V_{n-1},  \tag*{${{\bahram}_{n,2}}$}\\
{V^\prime_1}^{q_1}={V^\prime_2}^{q_2}=I
\end{cases}
\end{equation*}
Then we can present the generated algebra $A^{(n)}_2$ using the 
homeomorphism $\phi_2$ on $
\mathbf{X}_2:=\mathbb{Z}_{q_1}\times\mathbb{Z}_{q_2} \times\Bbb{T}^{n-2},$ 
$$\phi_2(v_1,v_2,\ldots,v_n)
=(\lambda v_1,\mu_1 v_1 v_2,\mu_2 v_2 v_3,v_3 v_4,\ldots,v_{n-1}v_n).$$
The flow $(\mathbb{Z},\mathbf{X}_2)$ that $\phi_2$ generates is usually not 
minimal, but we can restrict $\phi_2$
to $\mathbf{Y}_2 \times \Bbb{T}^{n-2}\subset \mathbf{X}_2$, where 
$\mathbf{Y}_2 \subset \mathbb{Z}_{q_1}\times
\mathbb{Z}_{q_2}$ is the finite set
\begin{align*}
\mathbf{Y}_2 
&=\{(v_1,v_2)\mid(v_1,v_2,1,1,\ldots,1)
\in\phi_2^r(1,1,\Bbb{T}^{n-2})~~\text{for 
some}~r \in\mathbb{N}\}\\
             &=\{(\lambda^r,\lambda^{\binom{r}{2}}\mu_1^r)\mid r \in 
\mathbb{N}\}.
\end{align*}

Then the flow $\mathcal{F}_2=(\mathbb{Z},\mathbf{Y}_2 \times \Bbb{T}_{n-2})$ 
is minimal; the proof of this is similar
to the minimality proof in case 1 above. So $\mathcal{C}(\mathbf{Y}_2 
\times\Bbb{T}^{n-2})\rtimes_{\phi_2}\mathbb{Z}$
is simple and isomorphic to $A^{(n)}_2$.\\

Continue this process assuming $\mu_2$ is also a root of unity and so on
down to the following last cases\\

($n-1$). When $\mu_{n-2}$ is also a root of unity, say a primitive 
$p_{n-2}$-th root of unity, let$$
q_{n-1}=\mathrm{lcm}\{q_{n-2},p_{n-2}\}=
\mathrm{lcm}\{q_1,p_1,\ldots,p_{n-2}\}$$
Then $V_{n-1}^{q_{n-1}}={\tilde\mu}_{n-1} I$, 
a multiple of the identity. If 
${\tilde\mu}_{n-1}$ is not
a root of unity, substitute $V_{n-1}=\mu_{n-1} V^\prime_{n-1}$ (as well as 
$V_{i}=\mu_{i}V^\prime_{i}$ for
$i=1,\ldots,{n-2}$) in  ${\bahram}_n$, where 
$\mu_{n-1}^{q_{n-1}}={\tilde\mu}_{n-1}$, and get
\begin{equation*}
\begin{cases}
[U,V^\prime_1]=\lambda,~~[U,V^\prime_{i+1}]=\mu_i 
V^\prime_i,~~(i=1,2,\ldots,n-2),\\
[U,V_n]=\mu_{n-1}V^\prime_{n-1}, \tag*{${{\bahram}_{n,n-1}}$}\\
{V^\prime_1}^{q_1}={V^\prime_2}^{q_2}=\ldots={V'}^{q_{n-1}}_{n-1}=I
\end{cases}
\end{equation*}
Then we can present the generated algebra $A^{(n)}_{n-1}$ using the 
homeomorphism $\phi_{n-1}$ on
$\mathbf{X}_{n-1}:=\mathbb{Z}_{q_1}\times\mathbb{Z}_{q_2} 
\times\ldots\times\mathbb{Z}_{q_{n-1}}\times\Bbb{T},$
$$\phi_{n-1}(v_1,v_2,\ldots,v_n)
=(\lambda v_1,\mu_1 v_1 v_2,\mu_2 v_2 v_3,\ldots,\mu_{n-1}v_{n-1}v_n).$$
The flow $(\mathbb{Z},\mathbf{X}_{n-1})$ that $\phi_{n-1}$ generates is 
usually not minimal, but we can restrict $\phi_{n-1}$
to $\mathbf{Y}_{n-1} \times \Bbb{T}\subset \mathbf{X}_{n-1}$, where 
$\mathbf{Y}_{n-1} \subset \mathbb{Z}_{q_1}\times
\mathbb{Z}_{q_2}\times\ldots\times\mathbb{Z}_{n-1}$ is the finite set
\begin{multline*}
\mathbf{Y}_{n-1} 
=\{(v_1,v_2,\ldots,v_{n-1})\mid(v_1,v_2,\ldots,v_{n-1},1)\in\phi_
{n-1}^r(1,1,\ldots,1,\Bbb{T}),\\
\text{for some}~r \in\mathbb{N}\}
\end{multline*}
\begin{equation*}
~~~~~~~=\{(\lambda^r,\lambda^{\binom{r}{2}}\mu_1^r,
\lambda^{\binom{r}{3}}\mu_1
^{\binom{r}{2}}\mu_2^r,\ldots,
\lambda^{\binom{r}{n-1}}\mu_1^{\binom{r}{n-2}}
\mu_2^{\binom{r}{n-3}}\ldots\mu_{n-2}^r)\mid r \in \mathbb{N}\}.
\end{equation*}
Then the flow $\mathcal{F}_{n-1}=(\mathbb{Z},\mathbf{Y}_{n-1} \times 
\Bbb{T})$ is minimal; the proof of this is similar
to the minimality proof in case 1 above. So $\mathcal{C}(\mathbf{Y}_{n-1} 
\times\Bbb{T})\rtimes_{\phi_{n-1}}\mathbb{Z}$
is simple and isomorphic to $A^{(n)}_{n-1}$.\\

$n$. When $\mu_{n-1}$ is a root of unity (as well as $\mu_i$ for 
$i=1,\ldots,n-2$), all
the unitaries are of finite order, so the generated $C^*$-algebra $A$ is 
finite dimensional.\\

The preceding comments are summarized in the next theorem.

\begin{theo}\label{theo:2}
A $C^*$-algebra $A$ is isomorphic to a simple infinite dimensional quotient 
of $C^{*}(\mathfrak{D}_n)$
if, and only if, $A$ is isomorphic to $\mathcal{A}_{n,\theta}$ for some 
irrational number $\theta$, or to
an $A^{(n)}_i=\mathcal{C}(\mathbf{Y}_{i} 
\times\Bbb{T}^{n-i})\rtimes_{\phi_{i}}\mathbb{Z}$ for a suitable finite set
$\mathbf{Y}_i$ as in cases $i=1,2,\ldots,n-1$ above.
\end{theo}

A result that has been used implicitly above, and should be stated 
explicitly, is the analogue of Theorem
\ref{theo:1} that holds for $A^{(n)}_i$ for $i=1,\ldots,n-1$.

\begin{theo}\label{theo:3}
For $i=1,2,\ldots,n-1$, $A^{(n)}_i$ is the unique (up to isomorphism) 
$C^*$-algebra generated
by unitaries $U,V_1,V_2,\ldots,V_n$ satisfying \linebreak
\textnormal{${\bahram}_{n,i}$}; $A^{(n)}_i$ is simple and is universal
for the equations \textnormal{${\bahram}_{n,i}$}.
\end{theo}

As for Theorem \ref{theo:1}, the result is a consequence of the minimality 
of the flow involved.

\begin{rema}
There are concrete representations for the $A^{(n)}_i$'s that 
are analogous to the
concrete representation on $L^2(\Bbb{T}^n)$ used in the definition of
$\mathcal{A}_{n,\theta}$ (and indeed the $A^{(n)}_i$'s could have been 
defined in terms of these
concrete representations). For $i=1,2,\ldots,n-1$, the representation for 
$A^{(n)}_i$
uses the flow $\mathcal{F}_i$ and is on $L^2(\mathbf{Y}_i 
\times\Bbb{T}^{n-i})$.
\end{rema}

  \section{Matrix representations for non-faithful quotients of 
$C^*(\mathfrak{D}_n)$}
The algebras $A^{(n)}_i$ above have representations as matrix algebras with 
entries
in simple $C^*$-algebras from faithful representations of groups of lower 
dimension.
More precisely, let $B^{(n)}_i$ be the universal $C^*$-algebra generated by 
unitaries
  $\tilde U,\tilde V_1,\ldots,\tilde V_{n-i}$ satisfying the following 
commutator relations
  \begin{equation*}
\begin{cases}
[\tilde{U},\tilde V_1]=\zeta_i\\
[\tilde{U},\tilde{V_2}]={\tilde V_1}^{C_i}\\
[\tilde{U},\tilde{V_3}]={\tilde V_1}^{\binom{C_i}{2}}{\tilde V_2}^{C_i}\\
\vdots \tag*{${{\widetilde{\bahram}}_{n,i}}$}\\
[\tilde{U},\tilde{V_{n-i}}]={\tilde V_1}^{\binom{C_i}{n-i-1}}{\tilde V_2}^
{\binom{C_i}{n-i-2}}\ldots{\tilde V_{n-i-1}}^{C_i}\\
[\tilde V_r,\tilde V_s]=1~~(1 \leq r,s \leq n-i),\\
\end{cases}
\end{equation*}
\linebreak
where $\zeta_i$ is not a root of unity. Then we will see in Theorem 
\ref{theo:4} below that $A^{(n)}_i \cong M_{C_i}(B^{(n)}_i)$.
Note that according to Remark 
\ref{rema:1},
$B^{(n)}_i$ is a simple (affine) Furstenberg transformation group 
$C^*$-algebra of
$\Bbb{T}^{n-i}$ since $C_i \ne 0$ (note that $b_{rs}=\binom{C_i}{s-r}$).
The group generated by ${{\widetilde{\bahram}}_{n,i}}$ (as $\mathfrak{D}_n$ 
was generated
by ${\bahram}_{n}$) is $\mathfrak{D}'_{n,i}$ ($=\Bbb{Z}^{n-i+1} \rtimes 
\Bbb{Z}$).
We will see in Corollary \ref{coro:2} below that $\mathfrak{D}'_{n,i}$ is 
isomorphic to a
subgroup of $\mathfrak{D}_{n-i}$. \\

More generally, let $\alpha:\Bbb{T}^n \rightarrow \Bbb{T}^n$ be an affine 
Furstenberg
transformation given by
$$\alpha(v_1,v_2,\ldots,v_n)=(e^{2 \pi i \zeta} v_1,v_1^{b_{12}}v_2,
v_1^{b_{13}}v_2^{b_{23}} v_3,\ldots
,v_1^{b_{1n}}v_2^{b_{2n}}\ldots v_{n-1}^{b_{n-1,n}}v_n),$$
where $\zeta$ is an irrational number and $b_{i,i+1}\neq 0$ for $i=1 
,\ldots, n-1$.
Let $A:=\mathcal{C}(\Bbb{T}^n) \rtimes_\alpha \Bbb{Z}$. Then $A$ is a simple
$C^*$-algebra, which is the unique $C^*$-algebra generated by unitaries 
$U,V_1,\ldots,V_n$
satisfying the commutator relations
\begin{equation*}
[U,V_1]=e^{2 \pi i \zeta},\,\,\,[U,V_2]=V_1^{b_{12}},\,\,\ \ldots,
\,\,\,[U,V_n]=V_1^{b_{1n}}\ldots V_{n-1}^{b_{n-1,n}} 
\tag*{${{\bahram}_{\alpha}}$}
\end{equation*}
(all other pairs of operators from $U,V_1,\ldots,V_n$ commute).
Let $\Gamma_\alpha$ denote the group generated by
${{\bahram}_{\alpha}}$. Then $\Gamma_\alpha=\mathbb{Z}^{n+1}
\rtimes_\gamma \mathbb{Z}$, in which $\gamma:\mathbb{Z}\rightarrow
\mathrm{GL}(n+1,\mathbb{Z})$ is given by
$\gamma(k)=\gamma_k=\mathsf{G}_{\alpha}^k$, where $\mathsf{G}_{\alpha}$ is a 
matrix defined as
\begin{equation}
\mathsf{G}_{\alpha}=\left[\begin{array}{cccccc}
1 & 1 & 0 & 0 &\cdots & 0\\
0 & 1 & b_{12} & b_{13} & \cdots& b_{1n}\\
0 & 0 & 1 & b_{23} & \cdots&\vdots\\
0 &\ddots&\ddots&\ddots&\ddots& b_{n-2,n}\\
\vdots&   & & 0 & 1 & b_{n-1,n}\\
0 &\cdots  & & & 0 & 1
\end{array}
\right]_{(n+1) \times (n+1)}
\end{equation}
One can check that $\Gamma_\alpha$ is the discrete group generated by
$x',y'_0,y'_1,\ldots,y'_n$ such that $x' y'_0=y'_0 x'$ and $y'_i y'_j=y'_j 
y'_i$ for $0 \leq i,j \leq  n$ and
\begin{equation}
[x',y'_1]=y'_0,\,\,\,[x',y'_2]={y'_1}^{b_{12}},\,\,\,\ldots,\,\,\,
[x',y'_n]={y'_1}^{b_{1n}}\ldots {y'}_{n-1}^{b_{n-1,n}}.
\end{equation}
Since $A$ is generated by a representation of $\Gamma_\alpha$, it is a 
quotient of
$C^*(\Gamma_\alpha)$.
\begin{lemm}\label{lemm:1}
Let $\Gamma_\alpha$ denote the group generated by 
\textnormal{${{\bahram}_{\alpha}}$} as above.
Then there is a monomorphism $\iota:
\Gamma_\alpha\rightarrow\mathfrak{D}_n$.
\end{lemm}
\begin{proof}
We construct $\iota$ recursively. Let $x,y_0,\ldots,y_n$ be the
generators of $\mathfrak{D}_n$ as in Section 1 satisfying (1) and
$x',y'_0,y'_1,\ldots,y'_n$ be the generators of $\Gamma_\alpha$
satisfying (4), as above. Take $\iota(x')=x$ and
$\iota(y'_j)=y_j$ for $j=0,1$. Now to define $\iota(y'_k)$, note
that we should have $\iota([x',y'_k])=[\iota(x'),\iota(y'_k)]=
[x,\iota(y'_k)]=\iota({y'_1}^{b_{1k}}\ldots{y'^{b_{k-1,k}}_{k-1}})
=\iota(y'_1)^{b_{1k}}
\ldots\iota(y'_{k-1})^{b_{k-1,k}}$. For example,
$[x,\iota(y'_2)]=\linebreak\iota(y'_1)^{b_{12}}$ $=y_1^{b_{12}}$, so
$\iota(y'_2)$ must be defined as $y_2^{b_{12}}$ according to
equations (1). Similarly,
$\iota(y'_3)=y_2^{b_{13}}y_3^{b_{12}b_{23}}$ and $\iota(y'_4)=
y_2^{b_{14}}y_3^{b_{12}b_{24}+b_{13}b_{34}}y_3^{b_{12}b_{23}b_{34}}$
and so on. Therefore using induction, one can show that
$\iota(y'_k)$ is of the form $y_2^{c_{2k}}y_3^{c_{3k}}\ldots
y_k^{c_{kk}}$ $(k>2)$ and $c_{kk}=\prod_{i=1}^{k-1} b_{i,i+1}\ne
0$. This will also guarantee the injectivity of $\iota$, for the
matrix $(c_{ij})$ (for $2 \leq i,j \leq n; ~c_{ij}=0$ for $i>j$)
is upper triangular with non-zero diagonal entries.
\end{proof}

\begin{rema}
In a similar way, one can see that there is also a monomorphism
$\iota':\mathfrak{D}_n \rightarrow\Gamma_\alpha$.
\end{rema}

\begin{lemm}
Let $\Gamma:=\Bbb{Z}^m \rtimes_\gamma\Bbb{Z}$, where 
$\gamma(k)=\mathsf{G}^k$
for some $\mathsf{G} \in \mathrm{GL}(m,\Bbb{Z})$. Then $[\Gamma,\Gamma]=
(\mathsf{G}-\mathsf{I})\Bbb{Z}^m \rtimes_\gamma\{0\}$ and
$$\frac{\Gamma}{[\Gamma,\Gamma]}\cong
\frac{\mathbb{Z}^m}{(\mathsf{G}-\mathsf{I})\mathbb{Z}^m}\oplus\Bbb{Z}=
\mathrm{coker}(\mathsf{G}-\mathsf{I})\oplus\Bbb{Z}.$$
\end{lemm}
\begin{proof}
For arbitrary $\textbf{x}_i \in\Bbb{Z}^m$ and $k_i \in\Bbb{Z}$ ($i=1,2$),
one has
\begin{align*}
[(\textbf{x}_1,k_1),(\textbf{x}_2,k_2)]&=(\textbf{x}_1,k_1)(\textbf{x}_2,k_2)
(\textbf{x}_1,k_1)^{-1}(\textbf{x}_2,k_2)^{-1}\\
&=(\textbf{x}_1,k_1)(\textbf{x}_2,k_2)(\mathsf{G}^{-k_1}
(-\textbf{x}_1),-k_1)(\mathsf{G}^{-k_2}
(-\textbf{x}_2),-k_2)\\
&=(\textbf{x}_1+\mathsf{G}^{k_1}(\textbf{x}_2)+\mathsf{G}^{k_2}
(-\textbf{x}_1)-\textbf{x}_2,0)\\
&=((\mathsf{G}^{k_1}-\mathsf{I})\textbf{x}_2-(\mathsf{G}^{k_2}-\mathsf{I})
\textbf{x}_1,0).
\end{align*}
So it is easily seen that
$$[\Gamma,\Gamma]=\{((\mathsf{G}-\mathsf{I})\textbf{x},0)\mid \textbf{x} 
\in\Bbb{Z}^m \}
=(\mathsf{G}-\mathsf{I})\Bbb{Z}^m \rtimes_\gamma\{0\}.$$\\
For the next part, one can check that 
$\varphi:\Gamma\rightarrow\mathrm{coker}
(\mathsf{G}-\mathsf{I})\oplus\Bbb{Z}$
defined by 
$\varphi(\textbf{x},k)=(\textbf{x}+(\mathsf{G}-\mathsf{I})\Bbb{Z}^m,k)$
is an epimorphism with 
$\ker\varphi=((\mathsf{G}-\mathsf{I})\Bbb{Z}^m,0)=[\Gamma,\Gamma]$.
\end{proof}

\begin{coro}
If $|b_{i,i+1}|\ne 1$ for some $i \in\{1,\ldots,n-1\}$, then $\Gamma_\alpha$ 
is not isomorphic to $\mathfrak{D}_n$.
\end{coro}
\begin{proof}
Using the preceding lemma, 
$\mathfrak{D}_n/[\mathfrak{D}_n,\mathfrak{D}_n]
\cong\mathrm{coker}(\mathsf{M}_n-\mathsf{I})\oplus\Bbb{Z}\cong\Bbb{Z}^2$ where 
the matrix $\mathsf{M}_n$ was defined in Section 2,
and $\Gamma_\alpha/[\Gamma_\alpha,\Gamma_\alpha]\cong
\mathrm{coker}(\mathsf{G}_{\alpha}-\mathsf{I})
\oplus\Bbb{Z}$ and $\mathsf{G}_{\alpha}$ is the matrix defined in (3).
But it is easy to see that $\mathrm{coker}(\mathsf{G}_{\alpha}-\mathsf{I})=
\mathrm{coker}(\mathsf{B})\oplus\Bbb{Z}$, where
\begin{equation*}
\mathsf{B}=\left[\begin{array}{cccccc}
b_{12} & b_{13} & \cdots& b_{1n}\\
0 & b_{23} & \cdots&\vdots\\
\vdots &\ddots&\ddots& b_{n-2,n}\\
0 & \cdots& 0 & b_{n-1,n}\\
\end{array}
\right]_{(n-1) \times (n-1)}.
\end{equation*}
Thus $\Gamma_\alpha/[\Gamma_\alpha,\Gamma_\alpha]\cong
\mathrm{coker}(\mathsf{B})\oplus\Bbb{Z}^2$. So a necessary condition for
$\Gamma_\alpha\cong\mathfrak{D}_n$ is that $\mathrm{coker}(\mathsf{B})=0$
i.e. $|\det\mathsf{B}|=\prod_{i=1}^{n-1} |b_{i,i+1}|=1$, which ends the 
proof.
\end{proof}

\begin{coro}\label{coro:2}
$\mathfrak{D}'_{n,i}$ is isomorphic to a subgroup of
$\mathfrak{D}_{n-i}$, and is not isomorphic to
$\mathfrak{D}_{n-i}$ unless $C_i=1$.
\end{coro}
\begin{proof}
Note that in this special case, $b_{rs}=\binom{C_i}{s-r}$. Now use the 
preceding
corollary and Lemma \ref{lemm:1}.
\end{proof}

  \begin{lemm}\label{lemm:2}
  Let $C_i=|\mathbf{Y}_i|$ be the cardinality of $\mathbf{Y}_i$ for
$i=1,\ldots,n-1$. Then the flow
  $(\mathbb{Z},\mathbf{Y}_i \times \Bbb{T}^{n-i})$
  generated by $\phi_i$ is topologically conjugate to a flow
  $(\mathbb{Z},\mathbb{Z}_{C_i} \times \Bbb{T}^{n-i})$
  generated by
  $$\psi_i:(v_i,v_{i+1},\ldots,v_n)\longmapsto(\lambda_i v_i,\eta_i v_i 
v_{i+1},
  v_{i+1} v_{i+2},\ldots,v_{n-1} v_n),$$
  where $\lambda_i$ is a primitive $C_i$-th root of unity and $\eta_i \in 
\Bbb{T}$
  is chosen appropriately.
  \end{lemm}
  \begin{proof}
  We construct a homeomorphism $\tau_i:\mathbf{Y}_i \times 
\Bbb{T}^{n-i}\rightarrow
  \mathbb{Z}_{C_i} \times \Bbb{T}^{n-i}$ that commutes with the actions of
  $\Bbb{Z}$ i.e. $\tau_i \circ \phi_i=\psi_i \circ \tau_i$. For the moment, 
fix $\nu_{i+1},\nu_{i+2},\ldots,\nu_{n-1}
  \in \Bbb{T}$, and define $\tau_i$ as follows for $v_{i+1},\ldots,v_n \in 
\Bbb{T}$:
  \begin{equation*}
  \tau_i(\overbrace{1,\ldots,1}^{i},v_{i+1},\ldots,v_n)=(1,\nu_{i+1} 
v_{i+1},\nu_{i+2}v_{i+2},\ldots,\nu_{n-1}v_{n-1},v_n)
   \end{equation*}
  \underline{step (1)}:
  \begin{multline*}
   \tau_i \circ 
\phi_i(\overbrace{1,\ldots,1}^{i},v_{i+1},\ldots,v_n)=\\
\underline{\tau_i(\lambda,\mu_1,
  \ldots,\mu_{i-1},\mu_i v_{i+1},v_{i+1} v_{i+2},\ldots,v_{n-1} v_n)}\\
  =\psi_i \circ \tau_i(\overbrace{1,\ldots,1}^{i},v_{i+1},\ldots,v_n) 
  \\
  =\psi_i(1,\nu_{i+1} 
v_{i+1},\nu_{i+2}v_{i+2},\ldots,\nu_{n-1}v_{n-1},v_n)\\
  =\underline{(\lambda_i,\eta_i \nu_{i+1} v_{i+1},\nu_{i+1}\nu_{i+2}
   v_{i+1} v_{i+2},\ldots,\nu_{n-2}\nu_{n-1}v_{n-2}v_{n-1},
  \nu_{n-1}v_{n-1} v_n)},
  \end{multline*}
\underline{step (2)}:
  \begin{align*}
   &\tau_i \circ \phi_i^2 (\overbrace{1,\ldots,1}^{i},v_{i+1},\ldots,v_n)=\\
   &\underline{\tau_i(\lambda^2,\lambda\mu_1^2,\mu_1 \mu_2 
^2,\ldots,\mu_{i-2}\mu_{i-1}^2,
  \mu_{i-1}\mu_i^2 v_{i+1},\mu_i v_{i+1}^2 v_{i+2},v_{i+1} v_{i+2}^2 
v_{i+3},
  \ldots,}\\
&~~~~~~~~~~~~~~~~~~~~~~~~~~~~~~~~~~~~~~~~~~~~~~~~~
~~~~~~~~~~~~~~~~~~~~~~~~~~~~~~~~~
  \underline{v_{n-2} v_{n-1}^2 v_{n})}\\
  &=\psi_i^2 \circ \tau_i(\overbrace{1,\ldots,1}^{i},v_{i+1},\ldots,v_n) 
  \\
  &=\psi_i(\lambda_i,\eta_i \nu_{i+1} v_{i+1},\nu_{i+1}\nu_{i+2}
   v_{i+1} v_{i+2},\ldots,\nu_{n-2}\nu_{n-1}v_{n-2}v_{n-1},
  \nu_{n-1}v_{n-1} v_n)
  \end{align*}
  \begin{multline*}
  =\underline{(\lambda_i^2,\lambda_i \eta_i^2 \nu_{i+1} v_{i+1},\eta_i \nu_{i+1}^2
  \nu_{i+2}v_{i+1}^2 v_{i+2}
  , \nu_{i+2} \nu_{i+3}^2 \nu_{i+4}v_{i+2} v_{i+3}^2 v_{i+4},\ldots,}\\
\underline{
    \nu_{n-3} \nu_{n-2}^2 \nu_{n-1}v_{n-3} v_{n-2}^2 v_{n-1},
  \nu_{n-2} \nu_{n-1}^2 v_{n-2} v_{n-1}^2 v_n)},
  \end{multline*}
  and so on down to\\
  \underline{step ($C_i$)}:
  \begin{align*}
   &\tau_i \circ \phi_i^{C_i} 
(\overbrace{1,\ldots,1}^{i},v_{i+1},\ldots,v_n)\\
   &=\tau_i(\lambda^{C_i},\lambda^{\binom{{C_i}}{2}}\mu_1^{C_i},
   \lambda^{\binom{{C_i}}{3}}\mu_1
                 ^{\binom{{C_i}}{2}}\mu_2^{C_i},\ldots,
                 \lambda^{\binom{{C_i}}{i}}\mu_1^{\binom{{C_i}}{i-1}}
                 \mu_2^{\binom{{C_i}}{i-2}}\ldots\mu_{i-1}^{C_i},\\				                 
&~~~~~~~~~~~~~~~~~~~~~\lambda^{\binom{{C_i}}{i+1}}\mu_1^{\binom{{C_i}}{i}}
                 \mu_2^{\binom{{C_i}}{i-1}}\ldots\mu_{i}^{C_i} v_{i+1},
                 \ldots,\\
                 &~~~~~~~~~~~~~~~~~~~~~~~~~~~~~~~~~~~~~~
                 \lambda^{\binom{{C_i}}{n}}\mu_1^{\binom{{C_i}}{n-1}}
                 \mu_2^{\binom{{C_i}}{n-2}}\ldots\mu_{i}^{\binom{C_i}{n-i}} 
v_{i+1}^{
                 \binom{C_i}{n-i+1}} v_{i+2}^{\binom{C_i}{n-i+2}}\ldots 
v_n)\\
				 &=\underline{\tau_i(\overbrace{1,\ldots,1}^{i},
                 \lambda^{\binom{{C_i}}{i+1}}\mu_1^{\binom{{C_i}}{i}}
                 \mu_2^{\binom{{C_i}}{i-1}}\ldots\mu_{i}^{C_i} v_{i+1},
                 \ldots,}\\&~~~~~~~~~~~~~~~~~~~~~~~~~~~~~~~~~~~~~~                 
\underline{\lambda^{\binom{{C_i}}{n}}\mu_1^{\binom{{C_i}}{n-1}}
                 \mu_2^{\binom{{C_i}}{n-2}}\ldots\mu_{i}^{\binom{C_i}{n-i}} 
v_{i+1}^{
                 \binom{C_i}{n-i+1}} v_{i+2}^{\binom{C_i}{n-i+2}}\ldots 
v_n)}
				 \end{align*}
				 \begin{align*}
  &=\psi_i^{C_i} \circ \tau_i(\overbrace{1,\ldots,1}^{i},v_{i+1},\ldots,v_n) 
  =\psi_i^{C_i}(1,\nu_{i+1} v_{i+1},v_{i+2},\ldots,v_n)
  \end{align*}
  \begin{multline*}
  =(\lambda_i^{C_i},\lambda_i^{\binom{C_i}{2}}\eta_i^{C_i}\nu_{i+1} v_{i+1},
  \lambda_i^{\binom{C_i}{3}}\eta_i^{\binom{C_i}{2}}\nu_{i+1}^{C_i}\nu_{i+2} 
  v_{i+1}^{C_i}v_{i+2},
  \ldots,\\
\lambda_i^{\binom{C_i}{n-i+1}}\eta_i^{\binom{C_i}{n-i}}
\nu_{i+1}^{\binom{C_i}{n-i-1}}\nu_{i+2}^{\binom{C_i}{n-i-2}}\ldots\nu_{n-1}^{C_i}
  v_{i+1}^{\binom{C_i}{n-i-1}} v_{i+2}^{\binom{C_i}{n-i-2}}\ldots v_{n})
  \end{multline*}
  \begin{multline*}
  =\underline{(1,\lambda_i^{\binom{C_i}{2}}\eta_i^{C_i}\nu_{i+1} v_{i+1},
  \lambda_i^{\binom{C_i}{3}}\eta_i^{\binom{C_i}{2}}\nu_{i+1}^{C_i}\nu_{i+2} 
  v_{i+1}^{C_i}v_{i+2},\ldots,}\\
  \underline{
\lambda_i^{\binom{C_i}{n-i+1}}\eta_i^{\binom{C_i}{n-i}}
\nu_{i+1}^{\binom{C_i}{n-i-1}}\nu_{i+2}^{\binom{C_i}{n-i-2}}\ldots\nu_{n-1}^{C_i}
  v_{i+1}^{\binom{C_i}{n-i-1}} v_{i+2}^{\binom{C_i}{n-i-2}}\ldots v_{n})}.
  \end{multline*}
  The definition of $\tau_i$ on $(\overbrace{1,\ldots,1}^{i},\Bbb{T}^{n-i})$ 
at this last
  step must coincide with the definition at the first step, so $\eta_i$ and 
$\nu_{i+1},\nu_{i+2},\ldots,\nu_{n-1}$
  are chosen to satisfy the equations
  \begin{align*}
  \lambda^{\binom{{C_i}}{i+1}}\mu_1^{\binom{{C_i}}{i}}
                 \mu_2^{\binom{{C_i}}{i-1}}\ldots\mu_{i}^{C_i}
                 &=\lambda_i^{\binom{C_i}{2}}\eta_i^{C_i},\\
  \lambda^{\binom{{C_i}}{i+2}}\mu_1^{\binom{{C_i}}{i+1}}
                 \mu_2^{\binom{{C_i}}{i}}\ldots\mu_{i}^{\binom{C_i}{2}}&=                 
\lambda_i^{\binom{C_i}{3}}\eta_i^{\binom{C_i}{2}}\nu_{i+1}^{C_i},\\
\lambda^{\binom{{C_i}}{i+3}}\mu_1^{\binom{{C_i}}{i+2}}
                 \mu_2^{\binom{{C_i}}{i+1}}\ldots\mu_{i}^{\binom{C_i}{3}}&=                 
\lambda_i^{\binom{C_i}{4}}\eta_i^{\binom{C_i}{3}}\nu_{i+1}^{\binom
{C_i}{2}}\nu_{i+2}^{C_i},\\
&~\vdots\\
\lambda^{\binom{{C_i}}{n}}\mu_1^{\binom{{C_i}}{n-1}}
                 \mu_2^{\binom{{C_i}}{n-2}}\ldots\mu_{i}^{\binom{C_i}{n-i}}
&=\lambda_i^{\binom{C_i}{n-i+1}}\eta_i^{\binom{C_i}{n-i}}
\nu_{i+1}^{\binom{C_i}{n-i-1}}\nu_{i+2}^{\binom{C_i}{n-i-2}}\ldots\nu_{n-1}^{C_i},
\end{align*}
i.e., choose $\eta_i$ to be a $C_i$-th root of
$$\lambda^{\binom{{C_i}}{i+1}}\mu_1^{\binom{{C_i}}{i}}                 
\mu_2^{\binom{{C_i}}{i-1}}\ldots\mu_{i}^{C_i}\overline{\lambda_i}                
^{\binom{C_i}{2}}=(-1)^{C_i+1}\lambda^{\binom{{C_i}}{i+1}}
\mu_1^{\binom{{C_i}}{i}}
                 \mu_2^{\binom{{C_i}}{i-1}}\ldots\mu_{i}^{C_i}$$
($\overline{\lambda_i}^{\binom{C_i}{2}}=(-1)^{C_i+1}$ since $\lambda_i$ is a 
primitive
$C_i$-th root of unity), and then choose $\nu_{i+1}$ to be a $C_i$-th root of
$$\lambda^{\binom{{C_i}}{i+2}}\mu_1^{\binom{{C_i}}{i+1}}
                 \mu_2^{\binom{{C_i}}{i}}\ldots\mu_{i}^{\binom{C_i}{2}}                 
\overline{\lambda_i}^{\binom{C_i}{3}}\overline{\eta_i}^{\binom{C_i}{2}},$$
and then choose $\nu_{i+2}$ to be a $C_i$-th root of
$$\lambda^{\binom{{C_i}}{i+3}}\mu_1^{\binom{{C_i}}{i+2}}
                 \mu_2^{\binom{{C_i}}{i+1}}\ldots\mu_{i}^{\binom{C_i}{3}}                 
\overline{\lambda_i}^{\binom{C_i}{4}}\overline{\eta_i}^{\binom{C_i}{3}}
\overline{\nu_{i+1}}^{\binom{C_i}{2}},$$
and so on choose $\nu_{i+3},\ldots,\nu_{n-1}$ recursively to satisfy the above
equations.\\

To see that $\tau_i$ commutes with the actions of $\Bbb{Z}$, take a point
$P \in \mathbf{Y}_i \times \Bbb{T}^{n-i}$. Then
$P=\phi_i^r(\overbrace{1,\ldots,1}^{i},v_{i+1},\ldots,v_n)$ for some $0 \leq 
r < C_i$
and $v_{i+1},\ldots,v_n \in \Bbb{T}$, and
\begin{align*}
\tau_i \circ \phi_i(P)&=\tau_i \circ
\phi_i^{r+1}(\overbrace{1,\ldots,1}^{i},v_{i+1},\ldots,v_n)\\&=\psi_i^{r+1}\circ 
\tau_i
(\overbrace{1,\ldots,1}^{i},v_{i+1},\ldots,v_n)\\&=\psi_i \circ \tau_i \circ 
\phi_i^r
(\overbrace{1,\ldots,1}^{i},v_{i+1},\ldots,v_n)\\&=\psi_i \circ \tau_i(P),
\end{align*}
as required.
  \end{proof}
\begin{coro}
$A^{(n)}_i 
\cong\mathcal{C}(\Bbb{Z}_{C_i}\times\Bbb{T}^{n-i})\rtimes_{\psi_i}\Bbb{Z}$,
where $\psi_i$ was introduced in the preceding lemma. In particular, the 
algebra $A^{(n)}_i$
is also a simple infinite dimensional quotient of 
$C^*(\mathfrak{D}_{n-i+1})$.
\end{coro}
\begin{proof}
The first part is clear, in view of the preceding lemma. For the second 
part, note that
$\mathcal{C}(\Bbb{Z}_{C_i}\times\Bbb{T}^{n-i})\rtimes_{\psi_i}\Bbb{Z}$ is 
generated
by unitaries satisfying \textnormal{${\bahram}_{n-i+1,1}$} (Theorem 
\ref{theo:3}).
\end{proof}
  \begin{theo}\label{theo:4}
  The algebra $A^{(n)}_i=\mathcal{C}(\mathbf{Y}_{i} 
\times\Bbb{T}^{n-i})\rtimes_{\phi_{i}}\mathbb{Z}$
  above is isomorphic to $M_{C_i}(B^{(n)}_i)$, where $B^{(n)}_i$ is the 
  $C^*$-algebra generated
  by \textnormal{${\widetilde{\bahram}}_{n,i}$} for 
$\zeta_i=(-1)^{C_i+1}\eta_i^{C_i}$ and
  $i=1,\ldots,n-1$.
  \end{theo}
  \begin{proof}
  By Lemma \ref{lemm:2},
  $A^{(n)}_i 
\cong\mathcal{C}(\Bbb{Z}_{C_i}\times\Bbb{T}^{n-i})\rtimes_{\psi_i}\Bbb{Z}$, 
hence
  this crossed product is simple too.
  For convenience, put $q:=C_i$ and $m:=n-i$ and let 
$D=\mathcal{C}(\Bbb{Z}_q \times
  \Bbb{T}^m)\rtimes_\psi\Bbb{Z}$, where $\psi(v,v_1,v_2,\ldots,v_m)=
  (\lambda v,\eta v v_1,v_1 v_2,\ldots,v_{m-1}v_m)$, in which 
  $\lambda$ is a primitive $q$-th root of unity
  and $\eta$ not a root of unity. Since $D$ is simple, it is the unique
  $C^*$-algebra generated by unitaries $U,V,V_1,\ldots,V_m$ such that 
$V^q=I$ and
  \begin{equation}
  [U,V]=\lambda,~[U,V_1]=\eta V,~[U,V_2]=V_1,~\ldots,~[U,V_m]=V_{m-1}
  \end{equation}
  (all other pairs of unitaries from $U,V,V_1,\ldots,V_m$ commute).\\
  Let $\zeta:=(-1)^{q+1}\eta^q$ and $B$ be the unique (simple) 
  $C^*$-algebra
  generated by unitaries $\tilde{U},\tilde{V}_1,\ldots,\tilde{V}_m$
  such that
  \begin{equation}
  [\tilde{U},\tilde {V}_1]=\zeta,~
[\tilde{U},\tilde{V}_2]=\tilde {V}_1^{q},~
\ldots,~
[\tilde{U},\tilde{V}_m]=\tilde {V}_1^{\binom{q}{m-1}}\tilde {V}_2^
{\binom{q}{m-2}}\ldots\tilde {V}_{m-1}^{q}
  \end{equation}
  (all other pairs of unitaries from 
$\tilde{U},\tilde{V}_1,\ldots,\tilde{V}_m$ commute).
  We prove that $D \cong M_q(B)$.
  Define unitaries $U',V',V'_1,\ldots,V'_m$ in $M_q(B)$ as follows
(all unspecified entries being 0).
\\
\linebreak
$U'$ has $\tilde{U}$ in the upper right-hand corner and 1's on the 
subdiagonal.\\
$V'=\mathrm{diag}(1,\bar{\lambda},\bar{\lambda}^2,\ldots,
\bar{\lambda}^{q-1})$.\\
$V'_1=\mathrm{diag}(b \tilde{V}_1,b \bar{\eta}\lambda\tilde{V}_1,
b \bar{\eta}^2 \lambda^3 \tilde{V}_1,\ldots,
b \bar{\eta}^{q-1}\lambda^{\binom{q}{2}}\tilde{V}_1)$.\\
$V'_2=\mathrm{diag}(c_{21} \tilde{V}_2,c_{22}\tilde{V}_1^{-1}\tilde{V}_2,
c_{23}\tilde{V}_1^{-2}\tilde{V}_2,
\ldots,c_{2q} \tilde{V}_1^{-(q-1)}\tilde{V}_2)$.\\
$V'_3=\mathrm{diag}(c_{31} \tilde{V}_3,c_{32}\tilde{V}_1 
\tilde{V}_2^{-1}\tilde{V}_3,
c_{33}\tilde{V}_1^{3}\tilde{V}_2^{-2}\tilde{V}_3,
\ldots,c_{3q} \tilde{V}_1^{\binom{q}{2}}\tilde{V}_2^{-(q-1)}\tilde{V}_3)$.\\
\vdots\\
$V'_m=\mathrm{diag}(c_{m1} \tilde{V}_m,c_{m2}\tilde{V}_1^{(-1)^{m-1}}
\tilde{V}_2^{(-1)^{m-2}}\ldots\tilde{V}_m,
\ldots,\\ ~~~~~~~~~~~~~~~~~~~~~~~~~~~~~~~~~~~~~~
c_{mq} \tilde{V}_1^{(-1)^{m-1}\binom{q+m-3}{m-1}}\tilde{V}_2^{(-1)^{m-2}
\binom{q+m-4}{m-2}}\ldots \tilde{V}_{m-1}^{-(q-1)}\tilde{V}_m).$\\
\linebreak
Thus the $j$-th entry of $V'_i$ equals
$c_{ij}\prod_{r=1}^i \tilde{V}_r^{(-1)^{i-r}\binom{j+i-r-2}{i-r}}$, and the 
constants
$b_i$, $c_{ij}\in\Bbb{T}$ must be chosen so that the unitaries
$U',V',V'_1,\ldots,V'_m$ satisfy ($4$). Note that $V'^q=I_q$ and $U'$ has 
the property
that
$$[U',\mathrm{diag}(\tilde{W}_1,\ldots\tilde{W}_q)]=\mathrm{diag}
(\tilde{U}\tilde{W}_q \tilde{U}^{-1}\tilde{W}_1^{-1},\tilde{W}_1 
\tilde{W}_2^{-1}
,\ldots,\tilde{W}_{q-1}\tilde{W}_q^{-1})$$
for arbitrary invertibles $\tilde{W}_1,\ldots\tilde{W}_q \in B$. Therefore\\

$~[U',V']=\mathrm{diag}(\bar{\lambda}^{q-1},\lambda,\ldots,\lambda)=\lambda 
I_q.$\\
\begin{align*}
[U',V'_1]&=\mathrm{diag}(\tilde{U}(b 
\bar{\eta}^{q-1}\lambda^{\binom{q}{2}}\tilde{V}_1)
\tilde{U}^{-1}(b \tilde{V}_1)^{-1},(b \tilde{V}_1) (b 
\bar{\eta}\lambda\tilde{V}_1)^{-1}
,\ldots,\\&~~~~~~~~~~~~~~~~~~~~~~~~~~~~~~~~~~~~~~~~(b 
\bar{\eta}^{q-2}\lambda^{\binom{q-1}{2}}
\tilde{V}_1)(b \bar{\eta}^{q-1}\lambda^{\binom{q}{2}}\tilde{V}_1)^{-1})\\
&=\mathrm{diag}(\bar{\eta}^{q-1}\lambda^{\binom{q}{2}}[\tilde{U},\tilde{V}_1],
\eta\bar{\lambda},\ldots,\eta\bar{\lambda}^{\binom{q}{2}-\binom{q-1}{2}})\\
&=\mathrm{diag}(\bar{\eta}^{q-1}\lambda^{\binom{q}{2}}\zeta,
\eta\bar{\lambda},\ldots,\eta\bar{\lambda}^{q-1})\\
&=\mathrm{diag}(\bar{\eta}^{q-1}\lambda^{\binom{q}{2}}\zeta,
\eta\bar{\lambda},\ldots,\eta\bar{\lambda}^{q-1})\\
&=\eta V'
\end{align*}
(note that since $\lambda$ is a primitive $q$-th root of unity, one can 
easily see
that $\lambda^{\binom{q}{2}}=(-1)^{q+1}$).
\begin{align*}
[U',V'_2]&=\mathrm{diag}(\tilde{U}(c_{2q} \tilde{V}_1^{-(q-1)}\tilde{V}_2)
\tilde{U}^{-1}(c_{21} \tilde{V}_2)^{-1},(c_{21} \tilde{V}_2)
(c_{22}\tilde{V}_1^{-1}\tilde{V}_2)^{-1}
,\ldots,\\&~~~~~~~~~~~~~~~~~~~~~~~~~~~~~~~~~~~~~~~~(c_{2,q-1} 
\tilde{V}_1^{-(q-2)}
\tilde{V}_2)(c_{2q} \tilde{V}_1^{-(q-1)}\tilde{V}_2)^{-1})\\
&=\mathrm{diag}(\frac{c_{2q}}{c_{21}}
(\tilde{U}\tilde{V}_1^{-(q-1)}\tilde{U}^{-1})
(\tilde{U}\tilde{V}_2\tilde{U}^{-1})\tilde{V}_2^{-1},
\frac{c_{21}}{c_{22}}\tilde{V}_1,
\ldots,\frac{c_{2,q-1}}{c_{2q}}\tilde{V}_1)\\
&=\mathrm{diag}(\frac{c_{2q}}{c_{21}}(\zeta^{-(q-1)}\tilde{V}_1^{-(q-1)})
(\tilde{V}_1^q 
\tilde{V}_2)\tilde{V}_2^{-1},\frac{c_{21}}{c_{22}}\tilde{V}_1,
\ldots,\frac{c_{2,q-1}}{c_{2q}}\tilde{V}_1)\\
&=\mathrm{diag}(\frac{c_{2q}}{c_{21}}\zeta^{-(q-1)}\tilde{V}_1,
\frac{c_{21}}{c_{22}}\tilde{V}_1,\ldots,
\frac{c_{2,q-1}}{c_{2q}}\tilde{V}_1)\\
&=V'_1.
\end{align*}
The last equality holds if, and only if,
\begin{equation}
\frac{c_{2q}}{c_{21}}\zeta^{-(q-1)}=b,~\frac{c_{21}}{c_{22}}=b 
\bar{\eta}\lambda,~\ldots,~
\frac{c_{2,q-1}}{c_{2q}}=b \bar{\eta}^{q-1}\lambda^{\binom{q}{2}}.
\end{equation}
By multiplying these equations, one obtains $b^q 
\bar{\eta}^{\binom{q}{2}}\lambda
^{\binom{q+1}{3}}=\bar{\zeta}^{q-1}$, so one must choose $b$ to be a $q$-th 
root of
$$\bar{\lambda}^{\binom{q+1}{3}}\eta^{\binom{q}{2}}\bar{\zeta}^{q-1}.$$
\begin{align*}
[U',V'_3]&=\mathrm{diag}(\tilde{U}(c_{3q} \tilde{V}_1^{\binom{q}{2}}
\tilde{V}_2^{-(q-1)}\tilde{V}_3)
\tilde{U}^{-1}(c_{31} \tilde{V}_3)^{-1},(c_{31} \tilde{V}_3)
(c_{32}\tilde{V}_1 \tilde{V}_2^{-1}\tilde{V}_3)^{-1}
,\\\ldots,&~~~~~~~~~~~~~~~~~~~~~~~~~~
(c_{3,q-1} \tilde{V}_1^{\binom{q-1}{2}}\tilde{V}_2^{-(q-2)}\tilde{V}_3)
(c_{3q} \tilde{V}_1^{\binom{q}{2}}\tilde{V}_2^{-(q-1)}\tilde{V}_3)^{-1})
\end{align*}
\begin{align*}
&=\mathrm{diag}(\frac{c_{3q}}{c_{31}}\tilde{U} \tilde{V}_1^{\binom{q}{2}}
\tilde{V}_2^{-(q-1)}\tilde{V}_3
\tilde{U}^{-1}\tilde{V}_3^{-1},\frac{c_{31}}{c_{32}}\tilde{V}_1^{-1}
\tilde{V}_2,\ldots,
\frac{c_{3,q-1}}{c_{3q}}\tilde{V}_1^{\binom{q-1}{2}-\binom{q}{2}}
\tilde{V}_2)\\
&=\mathrm{diag}(\frac{c_{3q}}{c_{31}}(\tilde{U} 
\tilde{V}_1^{\binom{q}{2}}
\tilde{U}^{-1})(\tilde{U}
\tilde{V}_2^{-(q-1)}\tilde{U}^{-1})(\tilde{U}\tilde{V}_3
\tilde{U}^{-1})\tilde{V}_3^{-1},\frac{c_{31}}{c_{32}}\tilde{V}_1^{-1}
\tilde{V}_2,\ldots,\\
&~~~~~~~~~~~~~~~~~~~~~~~~~~~~~~~~~~~~~~~~~~~~~~~~~~~~~~~~~~~~~~~~~~~~~~~~
\frac{c_{3,q-1}}{c_{3q}}\tilde{V}_1^{-(q-1)}\tilde{V}_2)
\end{align*}
\begin{align*}
&=\mathrm{diag}(\frac{c_{3q}}{c_{31}}(\zeta^{\binom{q}{2}}
\tilde{V}_1^{\binom{q}{2}})
(\tilde{V}_1^{-q(q-1)}\tilde{V}_2^{-(q-1)})
(\tilde{V}_1^{\binom{q}{2}}\tilde{V}_2^q \tilde{V}_3)\tilde{V}_3^{-1},
\frac{c_{31}}{c_{32}}\tilde{V}_1^{-1}\tilde{V}_2,\ldots,\\
&~~~~~~~~~~~~~~~~~~~~~~~~~~~~~~~~~~~~~~~~~~~~~~~~~~~~~~~~~~~~~~~~~~~~~~~~
\frac{c_{3,q-1}}{c_{3q}}\tilde{V}_1^{-(q-1)}\tilde{V}_2)\\
&=\mathrm{diag}(\frac{c_{3q}}{c_{31}}\zeta^{\binom{q}{2}}\tilde{V}_2,
\frac{c_{31}}{c_{32}}\tilde{V}_1^{-1}\tilde{V}_2,\ldots,
\frac{c_{3,q-1}}{c_{3q}}\tilde{V}_1^{-(q-1)}\tilde{V}_2)\\
&=V'_2.
\end{align*}
The last equality holds if, and only if,
\begin{equation}
\frac{c_{3q}}{c_{31}}\zeta^{\binom{q}{2}}=c_{21},
~\frac{c_{31}}{c_{32}}=c_{22}
,~\ldots,~\frac{c_{3,q-1}}{c_{3q}}=c_{2q}.
\end{equation}
By multiplying the above equalities, one obtains
$c_{21}\ldots c_{2q}=\zeta^{\binom{q}{2}}$ and combining with (6)
one must choose $c_{21}$ to be a $q$-th root of
$$\lambda^{\binom{q+2}{4}}\bar{\eta}^{\binom{q+1}{3}}b^{\binom{q}{2}}
\zeta^{\binom{q}{2}}$$
and then
$$c_{2j}=c_{21}\bar{\lambda}^{\binom{j+1}{3}}
\eta^{\binom{j}{2}}\bar{b}^{j-1}.$$\\

One can continue this procedure down to
\begin{align*}
[U',V'_m]&=\mathrm{diag}(\tilde{U}(c_{mq} 
\tilde{V}_1^{(-1)^{m-1}\binom{q+m-3}{m-1}}
\ldots \tilde{V}_{m-1}^{-(q-1)}\tilde{V}_m)
\tilde{U}^{-1}(c_{m1} \tilde{V}_m)^{-1},\\&~~~(c_{m1} 
\tilde{V}_m)
(c_{m2}\tilde{V}_1^{(-1)^{m-1}}
\tilde{V}_2^{(-1)^{m-2}}\ldots\tilde{V}_m)^{-1}
,\ldots,\\&~~~
(c_{m,q-1} \tilde{V}_1^{(-1)^{m-1}\binom{q+m-4}{m-1}}
\ldots \tilde{V}_{m-1}^{-(q-2)}\tilde{V}_m)\\
&~~~~~~~~~~~~~~~~~~~~~~~~~~~~~~~~~~~~~~~~~
(c_{mq} \tilde{V}_1^{(-1)^{m-1}\binom{q+m-3}{m-1}}
\ldots \tilde{V}_{m-1}^{-(q-1)}\tilde{V}_m)^{-1})
\end{align*}
\begin{align*}
&=\mathrm{diag}(\frac{c_{mq}}{c_{m1}}\tilde{U}
(\tilde{V}_1^{(-1)^{m-1}\binom{q+m-3}{m-1}}
\ldots \tilde{V}_{m-1}^{-(q-1)}\tilde{V}_m)
\tilde{U}^{-1}\tilde{V}_m^{-1},\\&~~~
\frac{c_{m1}}{c_{m2}}\tilde{V}_1^{(-1)^{m-2}}
\tilde{V}_2^{(-1)^{m-3}}\ldots\tilde{V}_{m-1},\ldots,
\end{align*}
\begin{align*}
&~~~\frac{c_{m,q-1}}{c_{mq}}\tilde{V}_1^{(-1)^{m-1}
\binom{q+m-4}{m-1}-\binom{q+m-3}{m-1}}
\ldots \tilde{V}_{m-2}^{\binom{q-1}{2}-\binom{q}{2}}
\tilde{V}_{m-1})
\end{align*}
\begin{align*}
&=\mathrm{diag}(\frac{c_{mq}}{c_{m1}}(\tilde{U}
\tilde{V}_1^{(-1)^{m-1}\binom{q+m-3}{m-1}}
\tilde{U}^{-1})\ldots (\tilde{U}
\tilde{V}_{m-1}^{-(q-1)}\tilde{U}^{-1})
(\tilde{U}\tilde{V}_m\tilde{U}^{-1})\tilde{V}_m^{-1},
\end{align*}
\begin{align*}
&~~~
\frac{c_{m1}}{c_{m2}}\tilde{V}_1^{(-1)^{m-2}}
\tilde{V}_2^{(-1)^{m-3}}\ldots\tilde{V}_{m-1},\ldots,\\
&~~~
\frac{c_{m,q-1}}{c_{mq}}
\tilde{V}_1^{(-1)^{m-2}\binom{q+m-4}{m-2}}
\ldots \tilde{V}_{m-2}^{-(q-1)}\tilde{V}_{m-1})
\end{align*}
\begin{align*}
&=\mathrm{diag}(\frac{c_{mq}}{c_{m1}}
(\zeta\tilde{V}_1)^{(-1)^{m-1}\binom{q+m-3}{m-1}}
\ldots (\tilde{V}_1^{\binom{q}{m-2}}
\ldots\tilde{V}_{m-1})^{-(q-1)}\\
&~~~~~~~~~~
(\tilde{V}_1^{\binom{q}{m-1}}\ldots\tilde{V}_m)\tilde{V}_m^{-1},
\frac{c_{m1}}{c_{m2}}\tilde{V}_1^{(-1)^{m-2}}
\tilde{V}_2^{(-1)^{m-3}}\ldots\tilde{V}_{m-1},\ldots,
\\&~~~~~~~~~~~~~~~~~~~~~~~~~~~~~~~~~~~~~~~
\frac{c_{m,q-1}}{c_{mq}}\tilde{V}_1^{(-1)^{m-2}
\binom{q+m-4}{m-2}}
\ldots \tilde{V}_{m-2}^{-(q-1)}\tilde{V}_{m-1})\\
&=\mathrm{diag}(\frac{c_{mq}}{c_{m1}}
\zeta^{(-1)^{m-1}\binom{q+m-3}{m-1}}
\prod_{s=1}^{m-1}\tilde{V}_s^{\sum_{r=1}^m 
(-1)^{m-r}\binom{q+m-r-2}{m-r}\binom{q}{r-s}},
\\&~~~~~~~~~~
\frac{c_{m1}}{c_{m2}}\tilde{V}_1^{(-1)^{m-2}}
\tilde{V}_2^{(-1)^{m-3}}\ldots\tilde{V}_{m-1},\ldots,
\\&~~~~~~~~~~~~~~~~~~~~~~~~~~~~~~~~~~~~~~~
\frac{c_{m,q-1}}{c_{mq}}\tilde{V}_1^{(-1)^{m-2}
\binom{q+m-4}{m-2}}
\ldots \tilde{V}_{m-2}^{-(q-1)}\tilde{V}_{m-1})\\
&=V'_{m-1}.
\end{align*}
Using the next lemma and the paragraph following it, we have
\begin{equation}
\sum_{r=1}^m (-1)^{m-r}\binom{q+m-r-2}{m-r}
\binom{q}{r-s}=\delta_{s,m-1}~~~
(1 \leq s \leq m-1)
\end{equation}
where $\delta$ denotes the delta function. 
Therefore the last equality
holds if, and only if,
\begin{equation}
\frac{c_{mq}}{c_{m1}}\zeta^{(-1)^{m-1}
\binom{q+m-3}{m-1}}=c_{m-1,1},~
\frac{c_{m1}}{c_{m2}}=c_{m-1,2}~,\ldots,~
\frac{c_{m,q-1}}{c_{mq}}=c_{m-1,q}.
\end{equation}
By multiplying the above equalities, one obtains $$\prod_{j=1}^q c_{m-1,j}=
\zeta^{(-1)^{m-1}\binom{q+m-3}{m-1}}.$$
One can see that $c_{m-1,1}$ must be chosen as a $q$-th root of
$$\zeta^{(-1)^{m-1}\binom{q+m-3}{m-1}}
\prod_{k=-1}^{m-2}c_{k1}^{(-1)^{m-k}\binom{q+m-k-2}{m-k}},$$
where $c_{-1,1}:=\lambda,~c_{01}:=\eta$ and $c_{11}:=b$. Then
$$c_{m-1,j}=\prod_{k=-1}^{m-1}c_{k1}^{(-1)^{m-k-1}\binom{q+m-k-3}{m-k-1}}.$$
Then one can show that
$$c_{m,j}:=\prod_{k=-1}^{m}c_{k1}^{(-1)^{m-k}\binom{q+m-k-2}{m-k}}$$
is a solution for (8), where $c_{k1}$ was chosen 
recursively in the previous steps for $k=1,\ldots,m-1$ and 
$c_{m1}\in\Bbb{T}$ is chosen
arbitrary. In accordance with the previous steps, we let $c_{m1}$ be a 
$q$-th root
of
\begin{equation*}\zeta^{(-1)^{m}\binom{q+m-2}{m}}
\prod_{k=-1}^{m-1}c_{k1}^{(-1)^{m-k+1}\binom{q+m-k-1}{m-k+1}}.
\end{equation*}
So the unitaries $U',V',V'_1,\ldots,V'_m$ satisfy (5). It just
remains to prove that they generate $M_q(B)=M_q(\Bbb{C})\otimes
B$. Let $E_{ij}$ denote the $q \times q$ matrix with all zero
entries except for a $1$ in the $i,j$ entry ($1 \leq i,j \leq q$).
These form a set of matrix units for $M_q(\Bbb{C})$. We first
show that all $E_{ij}$'s are generated by $U',V'$. Take
$E=\frac{1}{q}\sum_{k=0}^{q-1}V'^k$. Then
$E=\mathrm{diag}(1,0,0,\ldots,0)=E_{11}$ and one can check that
$E_{ij}=U'^{(i-1)}EV'^{(1-j)}$. Now that we have a copy of
$M_q(\Bbb{C})$ generated by $U'$ and $V'$, it is enough to prove
that the elements $E_{11}\otimes\tilde{U}$ and
$E_{11}\otimes\tilde{V}_i$ (for $i=1,\ldots,m$) belong to
$C^*(U',V',V'_1,\ldots,V'_m)$ (for these elements
generate $E_{11}\otimes B$ and moving this around with the matrix
units will generate $M_q(\Bbb{C})\otimes B$). But according to the
definition of $U',V',V'_1,\ldots,V'_m$ one can easily verify that
$E_{11}\otimes\tilde{U}=U'E_{q1}$ and
$E_{11}\otimes\tilde{V}_i=\bar{c}_{i1}E_{11}V'_i$ (for
$i=1,\ldots,m$), in which the constants $c_{i1}\in\Bbb{T}$ for
$i>1$ were introduced above and $c_{11}:=b$. This completes the
proof.
  \end{proof}

The following combinatorial lemma presents
  a generalized form of the identity 
$\binom{m-1}{k}=\binom{m}{k}-\binom{m-1}{k-1}$.

  \begin{lemm}
  Using Notation \textnormal{\ref{nota:1}}, we have
    $$\binom{m-q}{k}=\sum_{j \geq 0} (-1)^j \binom{m-j}{k-j}\binom{q}{j}$$
    for all $m,k \in\Bbb{Z}$ and $q \in \Bbb{N}$.
  \end{lemm}
\begin{proof}
First one checks that for $q=1$, the identity
$\binom{m-1}{k}=\binom{m}{k}-\binom{m-1}{k-1}$
holds for all $m,k \in\Bbb{Z}$. Then using
induction on $q$, one has
\begin{align*}
\binom{m-q-1}{k}=&\binom{m-q}{k}-\binom{m-q-1}{k-1}\\
=&\sum_{j \geq 0} (-1)^j \binom{m-j}{k-j}\binom{q}{j}-
\sum_{j \geq 0} (-1)^j \binom{m-1-j}{k-1-j}\binom{q}{j}\\
=&\sum_{j \geq 0} (-1)^j \binom{m-j}{k-j}\binom{q}{j}-
\sum_{j \geq 1} (-1)^{j-1} \binom{m-j}{k-j}\binom{q}{j-1}\\
=&\binom{m}{k}+\sum_{j \geq 1} (-1)^{j} 
\binom{m-j}{k-j}\{\binom{q}{j}+\binom{q}{j-1}\}\\
=&\sum_{j \geq 0} (-1)^j \binom{m-j}{k-j}\binom{q}{j}.
\end{align*}
\end{proof}

Now, to prove (9) in the preceding theorem, let $j:=r-s$ and $k:=m-s$. Then 
we have
\begin{align*}
\sum_{r=1}^m (-1)^{m-r}\binom{q+m-r-2}{m-r}\binom{q}{r-s}\\=&
\sum_{j=1-s}^k (-1)^{k+j}\binom{(q+k-2)-j}{k-j}\binom{q}{j}\\
=&(-1)^k \sum_{j \geq 0} (-1)^j \binom{(q+k-2)-j}{k-j}\binom{q}{j}\\
=&(-1)^k \binom{k-2}{k}=\delta_{k1}=\delta_{m-s,1}\\
=&\delta_{s,m-1}.
\end{align*}

  \section{$K$-theory for $\mathcal{A}_{n,\theta}$}

  In this section we study the $K$-theory of $\mathcal{A}_{n,\theta}$. 
First, we
  describe the method of computation for the $K$-groups of 
$\mathcal{C}(\Bbb{T}^n)\rtimes_\alpha\Bbb{Z}$,
  where $\alpha$ is an arbitrary homeomorphism of $\Bbb{T}^n$.
  To do this, we will pay attention to the algebraic structure of 
$K^{*}(\Bbb{T}^n)$. Note that it
  is sufficient to consider the special case of ``linear" homeomorphisms 
since as
  stated in the introduction, we know that every continuous function from 
$\Bbb{T}^n$ to $\Bbb{T}$
  is homotopic to a ``linear" function $f(v_1,\ldots,v_n)=
v_1^{a_1} v_2^{a_2}\ldots v_n^{a_n}$ ($a_1,a_2,\ldots,a_n \in\mathbb{Z}$), 
and that
the $K$-groups of $\mathcal{C}(\Bbb{T}^n)\rtimes_\alpha\Bbb{Z}$ depend up to
isomorphism only on the homotopy class of $\alpha$. \\

  It is well known that $K^{*}(\Bbb{T}^n)$
  is a $\mathbb{Z}_2$-graded ring and by the K\"{u}nneth
  formula \cite{mA67}, it is an exterior algebra (over $\mathbb{Z}$) on $n$ 
generators, where the
  elements of even degree are in $K^{0}(\Bbb{T}^n)$ and those of odd degree 
are in $K^{1}(\Bbb{T}^n)$.
  The generators of this exterior algebra correspond to the generators of 
the dual group $\mathbb{Z}^n$ of
  $\Bbb{T}^n$ \cite[p. 185]{jlT75}. Indeed, in this case the Chern character 
$$\mathrm{ch}:K^{*}(\Bbb{T}^n)\longrightarrow
  \Check{H}^{*}(\Bbb{T}^n,\mathbb{Q})$$ is integral and gives the Chern 
isomorphisms
  $$\mathrm{ch}_{0}:K^{0}(\Bbb{T}^n)\longrightarrow 
\Check{H}^{\text{even}}(\Bbb{T}^n,\mathbb{Z}),$$
  $$\mathrm{ch}_{1}:K^{1}(\Bbb{T}^n)\longrightarrow 
\Check{H}^{\text{odd}}(\Bbb{T}^n,\mathbb{Z}),$$
  where 
$\Check{H}^{*}(\Bbb{T}^n,\mathbb{Z})\cong
\Lambda^{*}_{\mathbb{Z}}(e_1,\ldots,e_n)$
is the (\v{C}ech) cohomology ring of $\Bbb{T}^n$ under the cup product, and
$\Check{H}^{k}(\Bbb{T}^n,\mathbb{Z})\cong
\Lambda^{k}_{\mathbb{Z}}(e_1,\ldots,e_n)$. 
On the other hand, $K^{*}(\Bbb{T}^n)\cong
K_{*}(\mathcal{C}(\Bbb{T}^n))$. So, by introducing $e_i:=[v_i]_1$
(i.e. the class in $K_{1}(\mathcal{C}(\Bbb{T}^n))$
of the coordinate function $v_i:\Bbb{T}^n \rightarrow\Bbb{T}$ as an element 
of
  $\mathcal{U}(\mathcal{C}(\Bbb{T}^n))~)$ for $i=1,\ldots,n$, we have the 
isomorphisms
$K_{*}(\mathcal{C}(\Bbb{T}^n))\cong
\Lambda^{*}_{\mathbb{Z}}(e_1,\ldots,e_n)\cong\Lambda^{*}\mathbb{Z}^n$,
and it is in a way that respects the canonical embedding of $\mathbb{Z}^n$. 
Moreover, such an isomorphism is
unique since only the identity automorphism of the ring 
$\Lambda^{*}\mathbb{Z}^n$ fixes each
element of $\mathbb{Z}^n$.

Now we use the Pimsner-Voiculescu six term exact sequence \cite{mP80} as the 
main tool for computing the $K$-groups
of $\mathcal{C}(\Bbb{T}^n)\rtimes_\alpha\Bbb{Z}$. Let 
$\alpha_*(=K_*(\alpha))$ be the
  ring automorphism of $K_{*}(\mathcal{C}(\Bbb{T}^n))$ induced by $\alpha$ 
and let $\alpha_i$ be the restriction
  of $\alpha_*$ on $K_{i}(\mathcal{C}(\Bbb{T}^n)); (i=0,1)$.
  Let $A:=\mathcal{C}(\Bbb{T}^n)\rtimes_\alpha\Bbb{Z}$.
  Then we have the following exact sequence.  
\[
\begin{CD}
      K_{0}(\mathcal{C}(\Bbb{T}^n)) 
@>\alpha_{0}-\mathrm{id}>>K_{0}(\mathcal{C}(\Bbb{T}^n)) @>\jmath_{0}>>K_0(A)\\
      @A \delta_{1}AA  @. @VV \delta_{0}V \\
      K_{1}(A) @<\jmath_{1}<< K_{1}(\mathcal{C}(\Bbb{T}^n)) 
@<\alpha_{1}-\mathrm{id}<<K_1(\mathcal{C}(\Bbb{T}^n))
      \end{CD}
      \]
     Here, $\jmath:\mathcal{C}(\Bbb{T}^n)\rightarrow A$ is the canonical 
embedding of
     $\mathcal{C}(\Bbb{T}^n)$ in $A$, $\jmath_0:=K_0(\jmath)$ and 
$\jmath_1:=K_1(\jmath)$. Also from now on
     $\mathrm{id}$ denotes the identity function on each underlying set.
      As a result, we have the following short exact sequences
       $$0 \longrightarrow \mathrm{coker}
      (\alpha_0-\mathrm{id})\longrightarrow 
K_0(\mathcal{C}(\Bbb{T}^n)\rtimes_\alpha\Bbb{Z})\longrightarrow
      \ker(\alpha_{1}-\mathrm{id})\longrightarrow 0,$$
       $$0 \longrightarrow \mathrm{coker}
      (\alpha_1-\mathrm{id})\longrightarrow 
K_1(\mathcal{C}(\Bbb{T}^n)\rtimes_\alpha\Bbb{Z})\longrightarrow
      \ker(\alpha_{0}-\mathrm{id})\longrightarrow 0.$$
      Since all the groups involved are abelian and $\ker(\alpha_{i}-\mathrm{id})$ is 
torsion-free $(i=0,1)$,
      these short exact sequences split and we have
      \begin{equation}
      \boxed{K_{0}(\mathcal{C}(\Bbb{T}^n)\rtimes_\alpha\Bbb{Z})\cong
      \mathrm{coker}(\alpha_{0}-\mathrm{id})\oplus \ker(\alpha_{1}-\mathrm{id}),}
      \end{equation}
      \begin{equation}
      \boxed{K_{1}(\mathcal{C}(\Bbb{T}^n)\rtimes_\alpha\Bbb{Z})\cong
      \mathrm{coker}(\alpha_{1}-\mathrm{id})\oplus \ker(\alpha_{0}-\mathrm{id}).}
      \end{equation}
      So, it suffices to determine the kernel and cokernel of 
$(\alpha_0-\mathrm{id})$ and $(\alpha_1-\mathrm{id})$ acting as
      endomorphisms on the finitely generated abelian groups \linebreak
      $\Lambda^{\text{even}}_{\mathbb{Z}}(e_1,\ldots,e_n)
      (\cong\Bbb{Z}^{2^{n-1}})$
      and $\Lambda^{\text 
{odd}}_{\mathbb{Z}}(e_1,\ldots,e_n)(\cong\Bbb{Z}^{2^{n-1}})$, respectively.
      Note that from the isomorphisms (11) and (12), the $K$-groups of 
$\mathcal{C}(\Bbb{T}^n)\rtimes_\alpha\Bbb{Z}$ are finitely generated abelian 
groups.
      Now since $\alpha_*$ becomes a ring homomorphism, it suffices
      to know the action of $\alpha_*$ on $e_1,\ldots,e_n$. In fact for a 
general
      basis element $e_{i_1}\wedge e_{i_2}\wedge\ldots\wedge e_{i_r}$ of 
$K_{*}(\mathcal{C}(\Bbb{T}^n))\cong
      \Lambda^{*}_{\mathbb{Z}}(e_1,\ldots,e_n )$ we have
      $$\alpha_*(e_{i_1}\wedge e_{i_2}\wedge\ldots\wedge
      e_{i_r})=\alpha_*(e_{i_1})\wedge\alpha_* 
(e_{i_2})\wedge\ldots\wedge\alpha_*(e_{i_r}).$$
      Thus if we consider $\{e_1,\ldots,e_n \}$ as the canonical basis of 
$\Bbb{Z}^n$ and
      take $\hat{\alpha}=\left.\alpha_*\right|_{\Bbb{Z}^n}$, we have 
$\alpha_*=\wedge^*\hat{\alpha}
      =\oplus_{r=1}^n \wedge^r \hat{\alpha}$,
      $\alpha_0=\wedge^{\text{even}}\hat{\alpha}=\oplus_{r \ge 0}
\wedge^{2r}\hat{\alpha}$ and
$\alpha_1=\wedge^{\text{odd}}\hat{\alpha}=\oplus_{r \ge 0}
\wedge^{2r+1}\hat{\alpha}$, where $\wedge^{i}\hat{\alpha}$ is the $i$-th 
exterior
power of $\hat{\alpha}$, which acts on 
$\Lambda^{i}{\mathbb{Z}^n}~(i=0,1,\ldots,n)$.
       Now let
      $\alpha=(f_1,\ldots,f_n)$ and $b_{ij}:=A_i[f_j]$ or in other words, 
assume that
      $f_i$ is homotopic to $(v_1,\ldots,v_n)\mapsto v_1^{b_{1i}} 
v_2^{b_{2i}}\ldots v_n^{b_{ni}}$.
      So we can write
      \begin{align*}
      \alpha_*(e_i)=\alpha_*[v_i]_1=[\alpha(v_i)]_1 &=[f_i( 
v_1,\ldots,v_n)]_1 =
      [v_1^{b_{1i}} v_2^{b_{2i}}\ldots v_n^{b_{ni}}]_1\\
      &=\sum_{j=1}^n b_{ji}[v_j]_1=\sum_{j=1}^n b_{ji}e_j.
      \end{align*}
    Therefore $\hat{\alpha}$
      acts on $\Bbb{Z}^n$ via the corresponding integer matrix
      $\mathsf{A}:=[b_{ij}]_{n \times n}\in\mathrm{GL}(n,\Bbb{Z})$
      and $\alpha_*$ acts on $\Lambda^*\Bbb{Z}^n$ via $\wedge^*\mathsf{A}$, 
and we have the following isomorphisms
\begin{align*}
K_0(\mathcal{C}(\Bbb{T}^n)\rtimes_\alpha\Bbb{Z})&\cong 
\mathrm{coker}(\alpha_0-\mathrm{id})\oplus \ker(\alpha_1-\mathrm{id})\\
                             &=\mathrm{coker}(\oplus_{r \ge 
0}\wedge^{2r}\hat{\alpha}-\mathrm{id})
\oplus \ker(\oplus_{r \ge 0}\wedge^{2r+1}\hat{\alpha}-\mathrm{id}),
\end{align*}
so we can write
$$\boxed{K_0(\mathcal{C}(\Bbb{T}^n)\rtimes_\alpha\Bbb{Z})
\cong \bigoplus_{r \ge 0}[\mathrm{coker}(\wedge^{2r}\hat{\alpha}-\mathrm{id})\oplus
\ker(\wedge^{2r+1}\hat{\alpha}-\mathrm{id})],}$$
and similarly
\begin{align*}
\boxed{K_1(\mathcal{C}(\Bbb{T}^n)\rtimes_\alpha\Bbb{Z})
\cong \bigoplus_{r \ge 0}[\mathrm{coker}(\wedge^{2r+1}\hat{\alpha}-\mathrm{id})
\oplus\ker(\wedge^{2r}\hat{\alpha}-\mathrm{id})].}
\end{align*}
Therefore, in order to compute the $K$-groups of
$\mathcal{C}(\Bbb{T}^n)\rtimes_\alpha\Bbb{Z}$, we must find the kernel and 
cokernel of
$\wedge^{r}\hat{\alpha}-\mathrm{id}$ as an endomorphism of $\Lambda^{r}{\mathbb{Z}^n}
~(r=0,1,\ldots,n)$. Note that the matrix of $\wedge^{r}\hat{\alpha}-\mathrm{id}$ with 
respect to the canonical basis
$\{e_{i_1}\wedge\ldots\wedge e_{i_r}|1 \le i_1 <  \ldots < i_r \le n \}$ 
with lexicographic
order is $\mathsf{A}_{n,r}:=\wedge^r \mathsf{A}- \mathsf{I}_{\binom{n}{r}}$, 
which is an integer
matrix of order $\binom{n}{r}$ ($\mathsf{I}_k$ is the identity matrix of 
order $k$, we often
omit $k$ whenever it is clear). So by computing
the kernel and cokenel of $\mathsf{A}_{n,r}$ for $r=0,1,\ldots,n$ with 
appropriate tools
(such as the Smith normal form for integer matrices \cite[p. 26]{mN72}), 
one can 
determine the
$K$-groups of $\mathcal{C}(\Bbb{T}^n)\rtimes_\alpha\Bbb{Z}$. The authors 
have
written a program in Maple 6.0 to do these computations.\\

We summarize the discussion above in the following proposition.

\begin{prop}\label{prop:1}
Let $\alpha$ be a homeomorphism of $\Bbb{T}^n$ and 
$\hat{\alpha}\in\mathrm{Aut}(\Bbb{Z}^n)$
be the restriction of $\alpha_*$ to $\Bbb{Z}^n $ (as above). Then
$\alpha_*=\wedge^*\hat{\alpha}=\oplus_{r=1}^n \wedge^r \hat{\alpha}$
on $K^*(\Bbb{T}^n)=\Lambda^*\Bbb{Z}^n$ and
\begin{align*}
K_0(\mathcal{C}(\Bbb{T}^n)\rtimes_\alpha\Bbb{Z})
&\cong \bigoplus_{r \ge 0}[\mathrm{coker}(\wedge^{2r}\hat{\alpha}-\mathrm{id})\oplus
\ker(\wedge^{2r+1}\hat{\alpha}-\mathrm{id})],\\
K_1(\mathcal{C}(\Bbb{T}^n)\rtimes_\alpha\Bbb{Z})&\cong \bigoplus_{r \ge 
0}[\mathrm{coker}(\wedge^{2r+1}\hat{\alpha}-\mathrm{id})\oplus
\ker(\wedge^{2r}\hat{\alpha}-\mathrm{id})].
\end{align*}
\end{prop}
\begin{coro}\label{coro:3}
The $K$-groups of $\mathcal{C}(\Bbb{T}^n)\rtimes_\alpha\Bbb{Z}$ are finitely 
generated
abelian groups with the same rank. Moreover, this common rank equals
$$\mathrm{rank}\ker(\wedge^* \hat{\alpha}-\mathrm{id})=
\sum_{r=0}^n \mathrm{rank}\ker(\wedge^r \hat{\alpha}-\mathrm{id}).$$
\end{coro}
\begin{proof}
Use the proposition and note that for any
$\varphi\in\mathrm{End}(\Bbb{Z}^n)$ one has 
$\mathrm{rank}\ker\varphi$
$=\mathrm{rank}\hspace{2pt}\mathrm{coker}\hspace{2pt}\varphi$.
\end{proof}

Now for the $K$-groups of $\mathcal{A}_{n,\theta}\cong 
\mathcal{C}(\Bbb{T}^n)\rtimes_\sigma\Bbb{Z}$,
the ``linearized" form of the corresponding affine homeomorphism $\sigma$ is 
as follows
$$(v_1,v_2,\ldots,v_n)\mapsto(v_1,v_1 v_2,\ldots,v_{n-1}v_n).$$
      So $\hat{\sigma}(e_i)=e_{i-1}+e_i$ for $i=1,\ldots,n$ $(e_0:=0)$.
    The matrix  with respect to the canonical basis $\{e_1,\ldots,e_n \}$ of
    $\mathbb{Z}^n$ that corresponds to $\hat{\sigma}$ is

      \[\mathsf{S}_n := \left[\begin{array}{ccccc}
1 & 1 & 0 & \cdots& 0\\
0 & 1 & 1 & &\vdots\\
0 &\ddots&\ddots&\ddots&0\\
\vdots&   & 0 & 1 & 1\\
0 &\cdots  & 0 & 0 & 1
\end{array}
\right]_{n \times n}.
\]\\

\begin{nota} 
We let $a_n:=\textnormal{rank}\hspace{2pt}K_0(\mathcal{A}_{n,\theta})=
\textnormal{rank}\hspace{2pt}K_1(\mathcal{A}_{n,\theta})$ and $a_{n,r}:=
\textnormal{rank}\ker(\wedge^r \mathsf{S}_n-\mathsf{I})$ for 
$r=0,1,\ldots,n$. From the preceding corollary
we have $$a_n=\textnormal{rank}\ker(\wedge^* 
\mathsf{S}_n-\mathsf{I})=\sum_{r=0}^n a_{n,r}.$$
\end{nota}\label{nota:2}
\begin{coro}
$K_i(C^*(\mathfrak{D}_n))\cong K_i(\mathcal{A}_{n+1,\theta})$ for $i=0,1$. 
In particular
$$\mathrm{rank}\hspace{2pt}K_0(C^*(\mathfrak{D}_n))
=\mathrm{rank}\hspace{2pt}K_1(C^*(\mathfrak{D}_n))=a_{n+1}.$$
\end{coro}
\begin{proof}
Since $\mathfrak{D}_n \cong \Bbb{Z}^{n+1}\rtimes_\eta \Bbb{Z}$, so
$C^*(\mathfrak{D}_n)\cong C^*(\Bbb{Z}^{n+1})\rtimes_{\tilde{\eta}} \Bbb{Z}
\cong\mathcal{C}(\Bbb{T}^{n+1})\rtimes_{\tilde{\eta}} \Bbb{Z}$ and the
integer matrix corresponding to $\tilde{\eta}$ is the $(n+1)\times(n+1)$ 
matrix
$\mathsf{M}_n$ introduced in Section 1, which is precisely 
$\mathsf{S}_{n+1}$. The rest of
proof follows from the last proposition.
\end{proof}

Some examples will illustrate the methods described.
\begin{example}
We compute the $K$-groups of $A^{5,5}_\theta$, which have been computed in 
\cite{sW02} by another method.
In fact, the Chern character and noncommutative geometry were used in 
\cite{sW02} to compute the
kernel and cokernel of $\sigma_i -\mathrm{id}$
($i$=0,1). As $A^{5,5}_\theta=\mathcal{A}_{3,\theta}$, we compute the kernel 
and cokernel of $\mathsf{S}_{3,r}:=\wedge^r \mathsf{S}_3-
\mathsf{I}_{\binom{3}{r}}$ for $r=0,1,2,3$, in which
\[\mathsf{S}_3 = \left[\begin{array}{ccccc}
1 & 1 & 0 \\
0 & 1 & 1 \\
0 & 0 & 1\\
\end{array}
\right].
\]\\
\begin{itemize}
\item $r=0$, $\mathsf{S}_{3,0}=\wedge^0 \mathsf{S}_3-\mathsf{I}_1=[0]$. So, 
$\ker \mathsf{S}_{3,0}=\mathbb{Z}$ and
$\mathrm{coker}\hspace{2pt} \mathsf{S}_{3,0}=\mathbb{Z}/\langle 0 
\rangle\cong
\mathbb{Z}.$
\item $r=1$, $\mathsf{S}_{3,1}=\wedge^1 
\mathsf{S}_3-\mathsf{I}_3=\left[\begin{smallmatrix}
1 & 1 & 0 \\
0 & 1 & 1 \\
0 & 0 & 1
\end{smallmatrix}
\right]
-\left[\begin{smallmatrix}
1 & 0 & 0 \\
0 & 1 & 0 \\
0 & 0 & 1\\
\end{smallmatrix}
\right]=\left[\begin{smallmatrix}
0 & 1 & 0 \\
0 & 0 & 1 \\
0 & 0 & 0\\
\end{smallmatrix}
\right].$
So,
\begin{align*}
\ker \mathsf{S}_{3,1}&=\{(x,y,z)\in
\mathbb{Z}^3|\left[\begin{smallmatrix}
0 & 1 & 0 \\
0 & 0 & 1 \\
0 & 0 & 0\\
\end{smallmatrix}
\right]\left[\begin{smallmatrix}
x \\
y \\
z \\
\end{smallmatrix}
\right]=\left[\begin{smallmatrix}
0 \\
0 \\
0 \\
\end{smallmatrix}
\right]\}=(\mathbb{Z},0,0)\cong\mathbb{Z}\\
\mathrm{coker}\hspace{2pt}\mathsf{S}_{3,1}&=
\mathbb{Z}^3/\mathsf{S}_{3,1}(\mathbb{Z}^3)
=\mathbb{Z}^3/\langle e_1,e_2 \rangle \cong \mathbb{Z}.
\end{align*}

\item $r=2$, $\mathsf{S}_{3,2}=\tiny\wedge^2 
\mathsf{S}_3-\mathsf{I}_3=\left[\begin{smallmatrix}
1 & 1 & 1 \\
0 & 1 & 1 \\
0 & 0 & 1\\
\end{smallmatrix}
\right]-\left[\begin{smallmatrix}
1 & 0 & 0 \\
0 & 1 & 0 \\
0 & 0 & 1\\
\end{smallmatrix}
\right]=\left[\begin{smallmatrix}
0 & 1 & 1 \\
0 & 0 & 1 \\
0 & 0 & 0\\
\end{smallmatrix}
\right].$
\begin{align*}
\ker \mathsf{S}_{3,2}&=\{(x,y,z)\in\mathbb{Z}^3|~y+z=z=0 
\}=(\mathbb{Z},0,0)\cong\mathbb{Z}\\
\mathrm{coker}\hspace{2pt}\mathsf{S}_{3,2}&=
\mathbb{Z}^3/\mathsf{S}_{3,2}(\mathbb{Z}^3)
=\mathbb{Z}^3/\langle e_1,e_2 \rangle \cong \mathbb{Z}.
\end{align*}
\item $r=3$, $\mathsf{S}_{3,3}=\wedge^3 
\mathsf{S}_3-\mathsf{I}_1=[0]$. So, 
$\ker \mathsf{S}_{3,3}
=\mathbb{Z}$ and $\mathrm{coker}\hspace{2pt} 
\mathsf{S}_{3,3}=\mathbb{Z}/\langle 0 \rangle\cong
\mathbb{Z}.$
\end{itemize}
Therefore
\begin{align*}
K_0(A^{5,5}_\theta)&=K_0({\mathcal{A}_{3,\theta}})\cong 
(\mathrm{coker}\hspace{2pt}
\mathsf{S}_{3,0}\oplus \mathrm{coker}\hspace{2pt}\mathsf{S}_{3,2})\oplus
(\ker \mathsf{S}_{3,1}\oplus \ker \mathsf{S}_{3,3})\\
&\cong\mathbb{Z}\oplus\mathbb{Z}\oplus\mathbb{Z}\oplus\mathbb{Z}=\mathbb{Z}^4,\\
K_1(A^{5,5}_\theta)&=K_1({\mathcal{A}_{3,\theta}})
\cong (\mathrm{coker}\hspace{2pt}\mathsf{S}_{3,1}\oplus 
\mathrm{coker}\hspace{2pt}\mathsf{S}_{3,3})
\oplus
(\ker \mathsf{S}_{3,0}\oplus \ker 
\mathsf{S}_{3,2})\\&\cong\mathbb{Z}\oplus\mathbb{Z}
\oplus\mathbb{Z}\oplus\mathbb{Z}=\mathbb{Z}^4.
\end{align*}
\end{example}
\begin{example}
Using our method, we have obtained the $K$-groups of
$\mathcal{A}_{n,\theta}$ by computer for $4 \leq n \leq 11$. We find the 
kernels
and cokernels of
$\mathsf{S}_{n,r}:=\wedge^r \mathsf{S}_n-\mathsf{I}_{\binom{n}{r}}$ for 
$r=0,1,\ldots,n$
using the Smith normal form theorem \cite[p. 26]{mN72}. The results  obtained are
stated in Table 1. Because of computational limitations, we do not have
any results yet for $n > 11$, except for $\{a_n\}$ for which we have 
obtained
the generating functions (see Subsection 6.1).

\begin{table}[h]
\begin{center}
\caption{\it The $K$-groups of $\mathcal{A}_{n,\theta}$ for $1 \leq n \leq 
11$}\vspace{4pt}
\begin{tabular}{|c|c|c|c|}
\hline
$n$  & $K_0(\mathcal{A}_{n,\theta})$ & $K_1(\mathcal{A}_{n,\theta})$ & $a_n$ 
\\ \hline
1 & $\Bbb{Z}^2$ & $\Bbb{Z}^2$ & 2 \\ 
2 & $\Bbb{Z}^3$ & $\Bbb{Z}^3$ & 3  \\ 
3 & $\Bbb{Z}^4$ & $\Bbb{Z}^4$ & 4  \\ 
4 & $\Bbb{Z}^6$ & $\Bbb{Z}^6$ & 6  \\ 
5 & $\Bbb{Z}^8$ & $\Bbb{Z}^8$ & 8  \\ 
6 & $\Bbb{Z}^{13}$ & $\Bbb{Z}^{13} \oplus \Bbb{Z}_2$ & 13  \\ 
7 & $\Bbb{Z}^{20}$ & $\Bbb{Z}^{20}$ & 20  \\ 
8 & $\Bbb{Z}^{32} \oplus \Bbb{Z}_8^{(2)} $ & $\Bbb{Z}^{32} \oplus 
\Bbb{Z}_{18}^{(2)}$ & 32  \\ 
9 & $\Bbb{Z}^{52} \oplus \Bbb{Z}_3^{(2)} \oplus \Bbb{Z}_9^{(2)}$ & 
$\Bbb{Z}^{52} \oplus \Bbb{Z}_3^{(2)} \oplus \Bbb{Z}_9^{(2)}$ & 52  \\ 
10 & $\Bbb{Z}^{90} \oplus \Bbb{Z}_{55}^{(4)}$ & $\Bbb{Z}^{90} \oplus 
\Bbb{Z}_{11}^{(2)} \oplus \Bbb{Z}_{99} \oplus \Bbb{Z}_{198} \oplus 
\Bbb{Z}_{2574}$ & 90  \\ 
11 & $\Bbb{Z}^{152} \oplus \Bbb{Z}_{11}^{(12)} \oplus\Bbb{Z}_{143}^{(4)} 
\oplus \Bbb{Z}_{286}^{(2)}$ & $\Bbb{Z}^{152} \oplus \Bbb{Z}_{11}^{(12)} 
\oplus \Bbb{Z}_{143}^{(4)} \oplus \Bbb{Z}_{286}^{(2)}$ & 152  \\ \hline
\end{tabular}
\end{center}
\end{table}
In this table, the group $\mathbb{Z}_{k}^{(m)}$ is the direct product of $m$ 
copies of the cyclic group $\mathbb{Z}_{k}=\mathbb{Z}/k \mathbb{Z}$, and 
it seems that the \text{$K$-groups} of $\mathcal{A}_{n,\theta}$ generally 
have torsion. The first example is $K_1(\mathcal{A}_{6,\theta})$;
 this is in fact because 
$\mathrm{coker}\hspace{2pt}\mathsf{S}_{6,3}=
\mathrm{coker}(\wedge^3 \mathsf{S}_6-\mathsf{I}_{20})\cong
\Bbb{Z}^{3} \oplus \Bbb{Z}_2$.
One of the things that we are interested in studying is 
the behavior of the sequence $\{a_n \}$. We will show below the importance 
of this sequence, namely $a_n$ is the common rank of the \text{$K$-groups} 
of a certain set of $C^*$-algebras including the Furstenberg transformation
group $C^*$-algebras $A_{F_{f,\theta}}$ of $\Bbb{T}^n$ introduced in 
\cite{rJ86}.
\end{example}

Note that by Proposition \ref{prop:1}, the $K$-groups of 
$\mathcal{C}(\Bbb{T}^n)\rtimes_\alpha
\Bbb{Z}$ are completely determined by the corresponding homomorphism 
$\hat\alpha\in\mathrm
{Aut}(\Bbb{Z}^n)$ and its exterior powers. From a computational point of 
view, we only
need the cokernels of the maps involved, since we know that for any 
endomorphism
$\varrho$ on $\Bbb{Z}^m$,
$\ker\varrho\cong\mathrm{coker}\varrho/(\mathrm{coker}\varrho)_{\mathrm{tor}}$.
When $\det\hat\alpha=1$, we don't even need to compute
all the cokernels. In other words, we have the next proposition, for which 
we recall a definition.
\begin{defi}
Let $\hat\alpha,\hat\beta\in\mathrm{End}(\Bbb{Z}^m)$. We say that 
$\hat\alpha$ is equivalent to
$\hat\beta$ over $\mathbb{Z}$ (and write $\hat\alpha$ \textnormal{equiv} 
$\hat\beta$)
if there exist $\hat{u},\hat{v} \in \mathrm{Aut}(\Bbb{Z}^m)$ such that
$\hat{u} \circ \hat\alpha\circ \hat{v}=\hat\beta$. Similarly, if 
$\mathsf{A}$ and $\mathsf{B}$ are integer $m \times m$ matrices,
$\mathsf{A}$ is equivalent to $\mathsf{B}$ if there exist ${\sf U,V} 
\in\mathrm{GL}(m,\Bbb{Z})$ such that
$\sf UAV=B$.
\end{defi}
Recall that $\hat\alpha$ \textnormal{equiv} $\hat\beta$ if, and only if, 
$\mathrm{coker}\hat\alpha\cong
\mathrm{coker}\hat\beta$ if, and only if, $\hat\alpha$ and $\hat\beta$ have the same 
Smith normal form.
Also, $\mathsf{A}$ equiv $\mathsf{B}$ if, and only if, $\mathsf{B}$ is obtainable from 
$\mathsf{A}$ by a finite number of elementary
operations. An elementary operation on an integer matrix is one of the 
following types:
interchanging two rows (or two columns), adding an integer multiple of one 
row (or column) to
another,  and multiplying a row (or column) by $-1$.
\begin{prop}
Let $\hat\alpha\in\mathrm{SL}(n,\Bbb{Z})$ (i.e. $\det\hat\alpha=1$). Then
$\wedge^r \hat\alpha-\mathrm{id}$ and $\wedge^{n-r} \hat\alpha-\mathrm{id}$ 
are equivalent as endomorphisms of $\Lambda^r 
\Bbb{Z}^n=\Lambda^{n-r}\Bbb{Z}^n=\Bbb{Z}^{\binom{n}{r}}$.
Equivalently, $\mathrm{coker}(\wedge^r \hat\alpha-\mathrm{id})\cong
\mathrm{coker}(\wedge^{n-r} \hat\alpha-\mathrm{id})$ for $r=0,1,\ldots,n$.
\end{prop}
\begin{proof}
We prove the equivalence of the endomorphisms as matrices with respect to
some basis. Let $\mathcal{E}=\{e_1,\ldots,e_n \}$
be a basis for $\Bbb{Z}^n$ and put $S=\{1,2,\ldots,n\}$. For $I 
=\{i_1,\ldots,i_r \}\subset S$, where
$1 \leq i_1<\ldots<i_r \leq n$, put $e_I=e_{i_1}\wedge\ldots\wedge 
e_{i_r}\in \Lambda^r
\Bbb{Z}^n$. Then $\mathcal{E}_r:=\{e_I \mid I \subset S~,~|I|=r \}$ is a 
basis
for $\Lambda^r \Bbb{Z}^n$. Let $\omega:=e_1 \wedge\ldots\wedge e_n$, which 
generates
$\Lambda^n \Bbb{Z}^n$. When $r=0$, $\wedge^0 \hat\alpha-\mathrm{id}=0$ and
$\wedge^n 
\hat\alpha(\omega)=\hat\alpha(e_1)\wedge\ldots\wedge\bar\alpha(e_n)
=(\det\hat\alpha)(e_1 \wedge\ldots\wedge e_n)=\omega$, so $\wedge^n 
\hat\alpha-\mathrm{id}=0$.
Now, fix an $r \in\{1,\ldots,n-1\}$. For an arbitrary subset $I \subset S$ 
with
$|I|=r$, take $J=\mathcal{E}\setminus I=\{j_1,\ldots,j_{n-r} \}$, so 
$|J|=n-r$.
Then $e_I \wedge e_J=(\mathrm{sgn}\hspace{2pt}\mu)\omega$, in which $\mu\in 
S_n$ is the permutation
that converts $(1,2,\ldots,n)$ to $(i_1,\ldots,i_r,j_1,\ldots,j_{n-r})$. It 
is easily seen
that $\mu=\mu_1 \ldots \mu_r$, where $\mu_k$ is the permutation that takes 
$i_k$ from
its position in $(1,2,\ldots,n)$ to its new position in 
$(i_1,\ldots,i_r,j_1,\ldots,j_{n-r})$.
One can see that $\mu_k$ is the combination of $i_{k}-(r-k+1)$ 
transpositions
($k=1,\ldots,r$). Thus
$$\mathrm{sgn}\hspace{2pt}\mu=\prod_{k=1}^r(-1)^{i_{k}-(r-k+1)}=
(-1)^{\ell(I)-\frac{r(r+1)}{2}},$$ where $\ell(I):=\sum_{k=1}^r
i_k$. Now take $m=\binom{n}{r}=\binom{n}{n-r}$ and let
$\mathcal{E}_r=\{e_{I_1},\ldots,\linebreak e_{I_m}\}$ be a basis
for $\Lambda^r \Bbb{Z}^n$. Now write $\mathcal{E}_{n-r}=
\{e_{J_1},\ldots,e_{J_m}\}$ as the basis for
$\Lambda^{n-r} \Bbb{Z}^n$ such that $J_k=\mathcal{E}\setminus
I_k$ for $k=1,\ldots,m$. From the above argument one can write
$$e_{I_i}\wedge 
e_{J_j}=(-1)^{\ell(I_i)-\frac{r(r+1)}{2}}\delta_{ij}\omega,$$
since if $i \ne j$ then $I_i \cap J_j \ne \emptyset$ thus $e_{I_i}\wedge 
e_{J_j}=0$.
Now let $\mathsf{A}=[a_{ij}]_{m \times m}$ and $\mathsf{B}=[b_{ij}]_{m 
\times m}$ be the corresponding
integer matrices of $\wedge^r \hat\alpha$ and $\wedge^{n-r} \hat\alpha$ with 
respect to
$\mathcal{E}_r$ and $\mathcal{E}_{n-r}$, respectively. So
$\wedge^r \hat\alpha(e_{I_i})=\sum_{p=1}^m a_{pi}e_{I_P}$ and
$\wedge^{n-r} \hat\alpha(e_{J_j})=\sum_{q=1}^m b_{qj}e_{J_q}$. What we want 
to show is that
$\mathsf{A}-\mathsf{I}$ is equivalent to $\mathsf{B}-\mathsf{I}$. We have
\begin{multline*}
\wedge^n \hat\alpha(e_{I_i}\wedge e_{J_j})
=(-1)^{\ell(I_i)-\frac{r(r+1)}{2}}\delta_{ij}\omega=\\
\wedge^r
\hat\alpha(e_{I_i})\bigwedge\wedge^{n-r}\hat\alpha(e_{J_j})=
\sum_{p,q=1}^m
a_{pi}b_{qj}(-1)^{\ell(I_p)-\frac{r(r+1)}{2}}\delta_{pq}\omega;
\end{multline*}
thus one obtains
\begin{equation}
\sum_{k=1}^m (-1)^{\ell(I_k)-\ell(I_i)}a_{ki}b_{kj}=\delta_{ij}.
\end{equation}
Therefore if we take $c_{ij}=(-1)^{\ell(I_j)-\ell(I_i)}a_{ji}$ and
$\mathsf{C}:=[c_{ij}]_{m \times m}$, then
$c_{ij}-\delta_{ij}=(-1)^{\ell(I_j)-\ell(I_i)}(a_{ji}-\delta_{ji})$;
therefore $\mathsf{C}-\mathsf{I}$ is obtained from $\mathsf{A}-\mathsf{I}$ 
by changing rows (and columns)
and occasionally multiplying some rows (and columns) by -1. So 
$\mathsf{C}-\mathsf{I}$ is equivalent to $\mathsf{A}-\mathsf{I}$.
On the other hand, (13) means that ${\sf CB=I}$. Thus 
$\mathsf{C}-\mathsf{I}=\mathsf{C}(\mathsf{B}-\mathsf{I})(-\mathsf{I})$, 
hence $\mathsf{B}-\mathsf{I}$
is equivalent to $\mathsf{C}-\mathsf{I}$. So $\mathsf{A}-\mathsf{I}$ is 
equivalent to $\mathsf{B}-\mathsf{I}$.
\end{proof}

\begin{coro}
If  $\det\hat\alpha=1$, then $\mathrm{rank}\ker(\wedge^{r}
\hat\alpha-\mathrm{id})= \mathrm{rank}\ker(\wedge^{n-r} \hat\alpha-\mathrm{id})$.
In particular, $a_{n,r}=a_{n,n-r}$ for $r=0,1,\ldots,n$.
\end{coro}

\begin{coro}
Let $A:=\mathcal{C}(\Bbb{T}^{2m-1})\rtimes_\alpha\Bbb{Z}$ for which the 
corresponding
homomorphism $\hat\alpha$ satisfies $\det\hat\alpha=1$. Then $K_0(A)\cong 
K_1(A)$ and
the (common) rank of the $K$-groups of $A$ is an even number.
In particular, for every Furstenberg transformation group $C^*$-algebra 
$A_{F_{f,\theta}}$
of an \text{odd-dimensional} torus (including $\mathcal{A}_{2m-1,\theta}$), 
one has
$K_0(A_{F_{f,\theta}})\cong K_1(A_{F_{f,\theta}})$.
\end{coro}
\begin{proof}
Combining the preceding proposition and Proposition \ref{prop:1}, one 
obtains
$$K_0(A)\cong K_1(A)\cong\bigoplus_{k=0}^{m-1}[\mathrm{coker}(\wedge^k 
\hat\alpha-\mathrm{id})
\oplus \ker(\wedge^k \hat\alpha-\mathrm{id})].$$
As a result, the rank of the $K$-groups of $A$ is an even number since the 
ranks of
the cokernel and kernel of an endomorphism coincide. Note that for 
$A_{F_{f,\theta}}$
the corresponding integer matrix of $\hat\alpha$ is an upper triangular 
matrix with
1's on diagonal. Thus $\det\hat\alpha=1.$
\end{proof}

\subsection{The rank $a_n$ of the $K$-groups of $\mathcal{A}_{n,\theta}$}
In this part we study some general properties of $a_n$. We specify some
\text{$C^*$-algebras} whose ranks of $K$-groups are related to the sequence 
$\{a_n \}$.
As an application, we
characterize the rank of the $K$-groups of Furstenberg transformation group 
\text{$C^*$-algebras}
$A_{F_{f,\theta}}$ and simple infinite dimensional quotients of
$C^*(\mathfrak{D}_n)$, which were studied in Sections 4 and 5. To this end, 
we remind
the reader of some linear algebraic properties of nilpotent and unipotent 
matrices.

\begin{defi}
Let $V$ be a (complex) vector space. An ${\hat\epsilon}\in 
\mathrm{End}_{\mathbb{C}}V$
is called nilpotent (respectively, unipotent) if ${\hat\epsilon}^k=0$ 
(respectively,
$({\hat\epsilon}-\mathrm{id})^k=0$) for some positive integer $k$. The minimum value 
of $k$ with this property is
called the degree of ${\hat\epsilon}$, denoted $\deg({\hat\epsilon})$.
\end{defi}

As an example, every upper (respectively, lower) triangular matrix with zeros on the 
diagonal is nilpotent.
Also, $\mathsf{S}_n$ as defined above is a unipotent matrix of degree $n$. 
Note that all eigenvalues of
a nilpotent (respectively, unipotent) matrix are zero (respectively, one). In particular, 
every unipotent matrix
is invertible. Thus every unipotent endomorphism is an automorphism.

\begin{coro}
Let $V$ be a finite dimensional complex vector space and ${\hat\epsilon}$ be 
a
nilpotent (respectively, unipotent) endomorphism of $V$.
Then $\deg({\hat\epsilon})$ is equal to the maximum order of its Jordan 
blocks.
\end{coro}
\begin{proof}
It suffices to prove the statement for the nilpotent case. Since all the 
eigenvalues of ${\hat\epsilon}$ are zero, each Jordan
block is a zero matrix of order one or in the form
$$\left[\begin{array}{ccccc}
0 & 1 & 0 & \cdots& 0\\
0 & 0 & 1 & &\vdots\\
0 &\ddots&\ddots&\ddots&0\\
\vdots&   & 0 & 0 & 1\\
0 &\cdots  & 0 & 0 & 0
\end{array}
\right],$$
which is a nilpotent matrix and its degree is the same as its order,
which is greater than 1. Now the rest of proof is clear.
\end{proof}

\begin{defi}
Let $V$ be a finite dimensional complex vector space and
${\hat\epsilon}\in \mathrm{End}_{\mathbb{C}}V$ be nilpotent (respectively, 
unipotent).
We say that ${\hat\epsilon}$ is of maximal degree if 
$\deg({\hat\epsilon})=\dim V$.
\end{defi}

\begin{coro}
Let $V$ be an \text{$n$-dimensional} complex vector space and
${\hat\epsilon}\in \mathrm{End}_{\mathbb{C}}V$ be nilpotent (respectively, 
unipotent).
Then $\deg({\hat\epsilon})\leq n$. If $\deg({\hat\epsilon})=n$,
the Jordan normal form of ${\hat\epsilon}$ is the following matrix of order 
$n$
$$\left[\begin{array}{ccccc}
0 & 1 & 0 & \cdots& 0\\
0 & 0 & 1 & &\vdots\\
0 &\ddots&\ddots&\ddots&0\\
\vdots&   & 0 & 0 & 1\\
0 &\cdots  & 0 & 0 & 0
\end{array}
\right].$$
In particular, all nilpotent (respectively, unipotent) matrices of maximal degree 
acting on $V$ are similar.
\end{coro}

\begin{proof}
Use the preceding corollary and the Jordan normal form theorem.
\end{proof}
\begin{example}
Let $\mathsf{B}=[b_{ij}]_{n \times n}$ be any upper triangular matrix whose 
diagonal elements
are all zero (respectively, one) and whose entries $b_{i,i+1}$ for $i=1,\ldots,n-1$ 
are all nonzero.
Then $\mathsf{B}$ is nilpotent (respectively, unipotent) of maximal degree. In fact, 
let $n$ be the order of $\mathsf{B}$
and let $b:=\prod_{i=1}^{n-1} b_{i,i+1}$, which is a nonzero number. Then 
one can
easily see that $\mathsf{B}^n=0$ (respectively, $(\mathsf{B}-\mathsf{I})^n=0$) and 
$\mathsf{B}^{n-1}$ (respectively, $(\mathsf{B}-\mathsf{I})^{n-1}$)
is a matrix with $b$ appearing in the northeast corner and zeros elsewhere.
So $\deg(\mathsf{B})=n$ i.e., $\mathsf{B}$ is of maximal degree.
\end{example}

As a useful conclusion, we compare the ranks of the $K$-groups of 
$C^*$-algebras of
the form $\mathcal{C}(\Bbb{T}^n)\rtimes_\alpha\Bbb{Z}$ in the following 
theorem,
which shows the importance of the rank $a_n$ of the $K$-groups of 
$\mathcal{A}_{n,\theta}$.
\begin{theo}\label{theo:5}
Let $A:=\mathcal{C}(\Bbb{T}^n)\rtimes_\alpha\Bbb{Z}$, in which $\alpha$ is a 
homeomorphism
of $\Bbb{T}^n$ whose corresponding integer matrix 
$\mathsf{A}\in\mathrm{GL}(n,\Bbb{Z})$ is unipotent of
maximal degree (i.e. $\deg(\mathsf{A})=n$). Then
$$\mathrm{rank}\hspace{2pt}K_0(A)=
\mathrm{rank}\hspace{2pt}K_1(A)=a_n
=\mathrm{rank}\hspace{2pt}K_0(\mathcal{A}_{n,\theta})
=\mathrm{rank}\hspace{2pt}K_1(\mathcal{A}_{n,\theta}).$$
In particular, the rank of the $K$-groups of any Furstenberg transformation 
group \text{$C^*$-algebra} $A_{F_{f,\theta}}=
\mathcal{C}(\Bbb{T}^n)\rtimes_{F_{f,\theta}}\Bbb{Z}$ is
equal to the rank of the $K$-groups of $\mathcal{A}_{n,\theta}$, i.e. to 
$a_n$.
\end{theo}
\begin{proof}
Let $\hat{\alpha}$ denote the restriction of $\alpha_*$ to $\Bbb{Z}^n$ and 
$\mathsf{A}$ be
the corresponding matrix of $\hat\alpha$ acting on $\Bbb{Z}^n$. Also let 
$\mathsf{S}_n$
be the corresponding matrix for $\mathcal{A}_{n,\theta}$. Now $\mathsf{A}$ 
is assumed to
be unipotent of maximal degree, and $\mathsf{S}_n$ is unipotent of maximal 
degree too. So
the last corollary yields that $\mathsf{A}$ is similar to $\mathsf{S}_n$. In 
fact the
Jordan normal form of $\mathsf{A}-\mathsf{I}$ is precisely 
$\mathsf{S}_n-\mathsf{I}$.
On the other hand, from Corollary \ref{coro:3}, we know that the rank of the 
$K$-groups of $A$ is equal to
$\mathrm{rank}\ker(\wedge^*\mathsf{A}-\mathsf{I})$. But the similarity 
of $\mathsf{A}$ and
$\mathsf{S}_n$ yields the similarity of $\wedge^*\mathsf{A}-\mathsf{I}$ and 
$\wedge^*\mathsf{S}_n-\mathsf{I}$ as
matrices acting on $\Lambda^*\Bbb{C}^n$. So 
$\dim_{\Bbb{C}}\ker(\wedge^*\mathsf{A}-\mathsf{I})
=\dim_{\Bbb{C}}\ker(\wedge^*\mathsf{S}_n-\mathsf{I})$, which yields the 
result.

For the second part, note that the corresponding integer matrix of a 
Furstenberg
transformation $F_{f,\theta}$ on $\Bbb{T}^n$ is of the form

\[\left[\begin{array}{ccccc}
1 & b_{12} & b_{13} & \cdots& b_{1n}\\
0 & 1 & b_{23} & &\vdots\\
0 &\ddots&\ddots&\ddots& b_{n-2,n}\\ \tag{$\heartsuit$}
\vdots&   & 0 & 1 & b_{n-1,n}\\
0 &\cdots  & 0 & 0 & 1
\end{array}
\right]_{n \times n},
\]
which, following the preceding example, is unipotent of maximal degree 
since $b_{i,i+1}\ne 0$ for $i=1,\ldots,n-1$. Now the proof of
the first part yields the result.
\end{proof}

\begin{rema}\label{rema:3}
In the preceding theorem, the basis for $\Bbb{Z}^n$ for the 
matrices involved
is $\{e_1,\ldots,e_n \}$, where $e_i:=[v_i]_1$
as introduced at the beginning of Section 6. But it is interesting to 
observe that
if $\hat{\alpha}$ is an arbitrary unipotent automorphism of $\Bbb{Z}^n$, 
then there is a basis
for $\Bbb{Z}^n$ with respect to which the integer matrix $\mathsf{A}$ of 
$\hat{\alpha}$ is of the form
$(\heartsuit)$ above (but not necessarily with $b_{i,i+1}\ne 0$ for 
$i=1,\ldots,n-1$,
unless $\hat{\alpha}$ is of maximal degree) \cite[Theorems 16, 18]{fjH63}.
The unipotence of $\hat{\alpha}$ also has important effects on the dynamics 
of the generated
flow. For example, if $\alpha$ is an affine transformation on $\Bbb{T}^n$ 
and $\hat{\alpha}$
is unipotent, then the dynamical system $(\Bbb{T}^n,\alpha)$ has 
quasi-discrete
spectrum \cite[Theorem 19]{fjH63}. More generally, let
$\alpha=(\boldsymbol{t},\mathsf{A})$ be an affine transformation on 
$\Bbb{T}^n$
and take $Z_p(\mathsf{A})=
\ker(\mathsf{A}^p-\mathrm{id})\subset\Bbb{Z}^n$ for $p \in \Bbb{N}$  and consider the 
following conditions\\

\begin{itemize}
\item [\textnormal{(1)}] $Z_1(\mathsf{A})=Z_p(\mathsf{A})$, $\forall p 
\in\Bbb{N},$
\item [\textnormal{(2)}] $\boldsymbol{t}$ is rationally independent over
$Z_1(\mathsf{A}),$ i.e., if $\boldsymbol{k}=(k_1,\ldots,k_n)\in 
Z_1(\mathsf{A})$
is such that $\langle \boldsymbol{t},\boldsymbol{k}\rangle:=\sum_{j=1}^n t_j 
k_j$ is a rational
number, then $\boldsymbol{k}=\boldsymbol 0$.
\item [\textnormal{(3)}] $Z_1(\mathsf{A})\ne \{0\},$
\item [\textnormal{(4)}] $\mathsf{A}$ is unipotent.
\end{itemize}~\\
\linebreak
Then $(\Bbb{T}^n,\alpha)$ is ergodic with respect to Haar measure if and 
only if $\alpha$
satisfies (1) and (2) \cite{fjH63}. Moreover, if $\alpha$ satisfies (1) 
through (4),
$(\Bbb{T}^n,\alpha)$ is minimal, uniquely ergodic with respect to Haar 
measure, and has
quasi-discrete spectrum. Conversely, any minimal transformation of 
$\Bbb{T}^n$ with
topologically quasi-discrete spectrum is conjugate to an affine 
transformation which
must satisfy (1) through (4) \cite{fjH65}. The $C^*$-algebras corresponding 
to such
actions are therefore simple and have a unique tracial state.
\end{rema}

\begin{coro}
Let $A$ be a simple infinite dimensional quotient of $C^*(\mathfrak{D}_n)$.
Then 
$\mathrm{rank}\hspace{2pt}K_0(A)=
\mathrm{rank}\hspace{2pt}K_1(A)=a_{n-i}$
for some $~i \in \{0,1,\ldots,n-1 \}$, determined by the
isomorphism $A \cong \mathcal{C}(\mathbf{Y}_{i}
\times\Bbb{T}^{n-i})\rtimes_{\phi_{i}}\mathbb{Z}$ in
Theorem \textnormal{\ref{theo:2}}.
\end{coro}
\begin{proof}
It was proved in Section 5 that $A$ is isomorphic to a matrix algebra over
a Furstenberg transformation group $C^*$-algebra $B^{(n)}_i$ on $\Bbb{T}^{n-i}$ 
for some
suitable $i$. So $K_j(A) \cong K_j(B^{(n)}_i)$ for $j=0,1$. The rest of proof
is clear from the preceding theorem.
  \end{proof}
  
  We will see in Proposition \ref{prop:2} below that $\{a_n \}$ is a 
strictly increasing
  sequence. Therefore the preceding corollary is a first step towards the 
classification
of the simple infinite dimensional
quotients of $C^*(\mathfrak{D}_n)$ by means of $K$-theory. But as is seen, 
the rank of
the $K$-groups can not distinguish the algebras in the same ``level" (i.e. 
those algebras that
are included in the same case (i) in Section 5, but with different values of 
the parameters).
The other powerful $K$-theoretical object that helps us do this is the trace 
invariant, i.e. the
range of the unique tracial state acting on the $K_0$-group.
\begin{prop}
Suppose $A \cong \mathcal{C}(\mathbf{Y}_{i} 
\times\Bbb{T}^{n-i})\rtimes_{\phi_{i}}\mathbb{Z}$
is a simple infinite dimensional quotient of $C^*(\mathfrak{D}_n)$ as in 
Theorem
\textnormal{\ref{theo:2}}. Then $A$ has a unique tracial state $\tilde\tau$
and $\tilde\tau_{*}K_0(A)=\frac{1}{C_i}(\Bbb{Z}+\Bbb{Z}\vartheta_i)$, where 
$C_i=|\mathbf{Y}_{i}|$
and $e^{2 \pi i \vartheta_i}=\zeta_i=(-1)^{C_i+1}\eta_i^{C_i}$ as in Lemma 
\textnormal{\ref{lemm:2}} and Theorem
\textnormal{\ref{theo:4}}.
\end{prop}
\begin{proof}
Following Theorem \ref{theo:4}, $A$ is isomorphic to 
$M_{C_i}(B^{(n)}_i)=M_{C_i}(\Bbb{C})\otimes B^{(n)}_i$, where 
$B^{(n)}_i$ is the simple
$C^*$-algebra generated
  by \textnormal{${\widetilde{\bahram}}_{n,i}$} for 
$\zeta_i=(-1)^{C_i+1}\eta_i^{C_i}$.
  By Corollary \ref{coro:1}, $B^{(n)}_i$ has a unique tracial state $\tau$ and
  $\tau_*K_0(B^{(n)}_i)=\Bbb{Z}+\Bbb{Z}\vartheta_i$, where $e^{2 \pi i 
\vartheta_i}=\zeta_i$
  \cite[Theorem 2.23]{rJ86}. Thus $A$ has the unique tracial state
  $\tilde\tau=(\frac{1}{C_i}\mathrm{Tr})\otimes\tau$, in which $\mathrm{Tr}$ 
is the usual
  trace on $M_{C_i}(\Bbb{C})$, and so
  $\tilde\tau_{*}K_0(A)=\frac{1}{C_i}(\Bbb{Z}+\Bbb{Z}\vartheta_i)$ 
\cite[Lemma 3.5]{rJ86}.
\end{proof}

\begin{coro}\label{coro:4}
$\mathcal{A}_{n,\theta}\cong\mathcal{A}_{n',\theta'}$ if, and only if, $n=n'$ and there 
exists
an integer $k$ such that $\theta=k \pm\theta'$. More generally, let 
$A_i^{(n)}\cong
\mathcal{C}(\mathbf{Y}_{i} \times\Bbb{T}^{n-i})\rtimes_{\phi_{i}}\mathbb{Z}$
be a simple infinite dimensional quotient of $C^*(\mathfrak{D}_n)$ with the
structure constants $\lambda,\mu_1,\ldots,\mu_i$ as in case (i)
of Section \textnormal{5} and let $A_{i'}^{(n')}\cong
\mathcal{C}(\mathbf{Y'}_{i'} 
\times\Bbb{T}^{n'-i'})\rtimes_{\phi'_{i}}\mathbb{Z}$
be a simple infinite dimensional quotient of $C^*(\mathfrak{D}_{n'})$ with 
the
structure constants $\lambda',\mu'_1,\ldots,\mu'_{i'}$.
Suppose that $C_i=|\mathbf{Y}_{i}|$ and $C'_{i'}=|\mathbf{Y'}_{i'}|$. Then 
$A_i^{(n)}\cong A_{i'}^{(n')}$ if, and only if,
$n-i=n'-i'$, $C_i=C'_{i'}$ and
$$\lambda^{\binom{{C_i}}{i+1}}\mu_1^{\binom{{C_i}}{i}}
                 \mu_2^{\binom{{C_i}}{i-1}}\ldots\mu_{i}^{C_i}=                 
\lambda'^{\binom{{C'_{i'}}}{i'+1}}{\mu'_1}^{\binom{{C'_{i'}}}{i'}}                 
{\mu'_2}^{\binom{{C'_{i'}}}{i'-1}}\ldots{\mu'_{i'}}^{C'_{i'}}$$
                 or
                 $$\lambda^{\binom{{C_i}}{i+1}}\mu_1^{\binom{{C_i}}{i}}
                 \mu_2^{\binom{{C_i}}{i-1}}\ldots\mu_{i}^{C_i}=                 
(\lambda'^{\binom{{C'_{i'}}}{i'+1}}{\mu'_1}^{\binom{{C'_{i'}}}{i'}}                 
{\mu'_2}^{\binom{{C'_{i'}}}{i'-1}}\ldots{\mu'_{i'}}^{C'_{i'}})^{-1}.$$
\end{coro}
\begin{proof}
Use the proposition and the fact that $\{a_n \}$ is a strictly increasing 
sequence
(see Proposition \ref{prop:2} below). Note that 
$$\zeta_i=(-1)^{C_i+1}\eta_i^{C_i}=
\lambda^{\binom{{C_i}}{i+1}}\mu_1^{\binom{{C_i}}{i}}
                 \mu_2^{\binom{{C_i}}{i-1}}\ldots\mu_{i}^{C_i}.$$
\end{proof}
\begin{rema}
Note that $C_i=|\mathbf{Y}_{i}|$ is completely determined by the structural 
constants
$\lambda,\mu_1,\ldots,\mu_{i-1}$ (which are roots of unity). More precisely
$$C_i=\min\{r \in\Bbb{N}\mid\lambda^r=\lambda^{\binom{r}{2}}\mu_1^r=\ldots=
\lambda^{\binom{r}{i}}\mu_1^{\binom{r}{i-1}}
                 \mu_2^{\binom{r}{i-2}}\ldots\mu_{i-1}^r=1\}.$$
				 As an example, see \cite[Lemma 5.4]{pM97}.
\end{rema}~\\
\subsection{Some combinatorial properties of $a_{n,r}$}
As mentioned before, our main goal is to describe $a_n$ as the rank of the 
$K$-groups of
$\mathcal{A}_{n,\theta}$. Since $a_n=\sum_{r=0}^n a_{n,r}$, we study 
$a_{n,r}$ first.
So this part is devoted to some combinatorial properties of $a_{n,r}$ as the 
rank of
$\ker(\wedge^r \hat\alpha-\mathrm{id})$ for $r=0,1,\ldots,n$. In other words we show 
that $a_{n,r}$
equals the number of partitions of $[r(n+1)/2]$ to $r$ distinct positive
integers not greater than $n$.  To do this, we will use properties of the
irreducible representations of the semisimple Lie algebra 
$\mathfrak{sl}_2(\Bbb{C})$.
First, we need some more elementary properties of nilpotent linear mappings.
\begin{lemm}\label{lemm:3}
Let $V$ be a complex vector space and ${\hat\epsilon}\in 
\mathrm{End}_{\mathbb{C}} V$
be nilpotent of degree $k$. Then
$\exp({\hat\epsilon})$ is unipotent of degree $k$. Moreover, 
$\exp({\hat\epsilon})-\mathrm{id}$ is similar to ${\hat\epsilon}$.
\end{lemm}
\begin{proof}
For the first part, we know that
$\exp(\hat\epsilon)-\mathrm{id}={\hat\epsilon}+{\hat\epsilon}^2/2!+\ldots+
{\hat\epsilon}^{k-1}/(k-1)!={\hat\epsilon}\omega$,
where $\omega:=\mathrm{id}+{\hat\epsilon}/2!+\ldots
+{\hat\epsilon}^{k-2}/(k-1)!$ commutes with ${\hat\epsilon}$ and is 
invertible since it is unipotent. So, $(\exp(\hat\epsilon)-\mathrm{id})^r
=({\hat\epsilon}\omega)^r={\hat\epsilon}^r \omega^r$ for each positive 
integer $r$. Thus $\exp({\hat\epsilon})-\mathrm{id}$ is unipotent with
the same degree of ${\hat\epsilon}$. For the second part, using the Jordan 
normal form of ${\hat\epsilon}$, it is sufficient
to prove the statement for the special case when ${\hat\epsilon}$ is a 
Jordan block with zeros on the diagonal. Since in this
case ${\hat\epsilon}$ is of maximal degree, by the first part, 
$\exp({\hat\epsilon})-\mathrm{id}$ is also of maximal degree. Therefore they are
similar by Corollary $6.3$.
\end{proof}

Let $V$ be an arbitrary (complex) vector space and ${\hat\phi}:V \rightarrow 
V$ be a
linear mapping. Then ${\hat\phi}$ can be extended in a unique way to a 
homomorphism
$\wedge^{*}{\hat\phi}:\Lambda^{*}V \rightarrow \Lambda^{*}V$ such that 
$\wedge^{*}{\hat\phi}(1)=1$, yielding
$$\wedge^{*}{\hat\phi}(x_1 \wedge\ldots\wedge x_p)={\hat\phi}(x_1) 
\wedge\ldots\wedge {\hat\phi}(x_p)~~~~~~(x_i \in V).$$
Also, ${\hat\phi}$ can be extended in a unique way to a derivation 
$\mathcal{D}^{\hspace{1pt}*}{\hat\phi}:\Lambda^{*}V \rightarrow
\Lambda^{*}V$, yielding $$\mathcal{D}^{\hspace{1pt}*}{\hat\phi}(x_1 
\wedge\ldots\wedge x_p)=\sum_{i=1}^{p}x_1 \wedge\ldots
\wedge{\hat\phi}(x_i)\wedge\ldots\wedge x_p~~~~~(~p \geq 2 ,x_i \in V).$$
Let's define $\wedge^{r}{\hat\phi}:=\wedge^{*}{\hat\phi}|_{\Lambda^{r}V}$ 
and $\mathcal{D}^{\hspace{1pt}r} 
{\hat\phi}:=\mathcal{D}^{\hspace{1pt}*}{\hat\phi}|_{\Lambda^{r}V}$
as induced linear mappings on the $r$-th exterior power of $V$ $(r \geq 0)$. 
Then we have
$$\wedge^{*}{\hat\phi}=\bigoplus_{r \geq 0}\wedge^{r}{\hat\phi}~,
~~~~\mathcal{D}^{\hspace{1pt}*}{\hat\phi}
=\bigoplus_{r \geq 0}\mathcal{D}^{\hspace{1pt}r}{\hat\phi}~.$$
\begin{lemm}
With the above notation, if ${\hat\phi}:V \rightarrow V$ is nilpotent, 
$\wedge^{r}{\hat\phi}$ and
$\mathcal{D}^{\hspace{1pt}r}{\hat\phi}$ are also nilpotent for $r \geq 1$. 
If $V$ is finite dimensional,
$\mathcal{D}^{\hspace{1pt}*}{\hat\phi}$ is nilpotent.
\end{lemm}
\begin{proof}
Assume that ${\hat\phi}^t=0$ for some $t \in \mathbb{N}$. We have
$(\wedge^r {\hat\phi})^t(x_1 \wedge\ldots\wedge 
x_r)={\hat\phi}^t(x_1)\wedge\ldots\wedge {\hat\phi}^t(x_r)=0$
so $(\wedge^r {\hat\phi})^t=0$ and $\wedge^r {\hat\phi}$ is nilpotent. For 
$\mathcal{D}^{\hspace{1pt}r}{\hat\phi}$ we know that
$\mathcal{D}^{\hspace{1pt}r}{\hat\phi}(x_1 \wedge\ldots\wedge 
x_r)=\sum_{i=1}^{r}x_1 \wedge\ldots
\wedge{\hat\phi}(x_i)\wedge\ldots\wedge x_r$ and one can deduce that
$$\mathcal{D}^{\hspace{1pt}r}{\hat\phi}^p(x_1 \wedge\ldots\wedge 
x_r)=\underset{(i_j \geq 0)}{\sum_{i_1+\ldots+i_r=p}}
\frac{p!}{(i_1)!\ldots(i_r)!}~{\hat\phi}^{i_1}x_1 
\wedge\ldots\wedge{\hat\phi}^{i_r}x_r.$$
Now since $i_1+\ldots+i_r=p$, there exists an $i_j$ with $i_j \geq p/r$. So 
if
$p \geq r t$ then ${\hat\phi}^{i_j}=0$ and 
$\mathcal{D}^{\hspace{1pt}r}{\hat\phi}^p=0$. Thus 
$\mathcal{D}^{\hspace{1pt}r}{\hat\phi}$ is nilpotent.\\

For the next part, let $m:=\dim V$. Since 
$\mathcal{D}^{\hspace{1pt}*}{\hat\phi}=\bigoplus_{r \geq 0}^{m}
\mathcal{D}^{\hspace{1pt}r}{\hat\phi}$ and ${\hat\phi}_0=0$, from the first 
part we have $(\mathcal{D}^{\hspace{1pt}*}{\hat\phi})^{mt}=0$
and $\mathcal{D}^{\hspace{1pt}*}{\hat\phi}$ is nilpotent too.
  \end{proof}

  \begin{lemm}
  Let ${\hat\phi}:V \rightarrow V$ be a nilpotent linear mapping. Then
  $$\exp(\mathcal{D}^{\hspace{1pt}*}{\hat\phi})=\wedge^{*} 
\exp({\hat\phi})$$
  on $\Lambda^{*}V$.
  \end{lemm} 
  \begin{proof}
  We have
  {\allowdisplaybreaks
  \begin{align*}
\exp(\mathcal{D}^{\hspace{1pt}*}{\hat\phi})(x_1 \wedge\ldots\wedge 
x_r)&=\sum_{p \geq 0}
  \frac{1}{p!}\mathcal{D}^{\hspace{1pt}r}{\hat\phi}^p(x_1 \wedge\ldots\wedge 
x_r)\\
  &=\sum_{p \geq0}\frac{1}{p!}(\underset{(i_j \geq 
0)}{\sum_{i_1+\ldots+i_r=p}}
\frac{p!}{(i_1)!\ldots(i_r)!}~{\hat\phi}^{i_1}x_1 
\wedge\ldots\wedge{\hat\phi}^{i_r}x_r)\\
&=\sum_{i_j \geq 0}\frac{1}{(i_1)!\ldots(i_r)!}~{\hat\phi}^{i_1}x_1 
\wedge\ldots\wedge{\hat\phi}^{i_r}x_r\\
&=(\sum_{i_1 \geq 0}\frac{1}{(i_1)!}{\hat\phi}^{i_1}x_1)\wedge\ldots\wedge
(\sum_{i_r \geq 0}\frac{1}{(i_r)!}{\hat\phi}^{i_r}x_r)\\
&=(\exp({\hat\phi}) x_1)\wedge\ldots\wedge(\exp({\hat\phi}) x_r)\\
&=\wedge^{*} \exp({\hat\phi})(x_1 \wedge\ldots\wedge x_r),
   \end{align*}}
   \linebreak
  which yields the result. Note that all sums in the above equalities are 
finite according to the
  previous lemma.
   \end{proof}
   \begin{rema}
   The nilpotence of ${\hat\phi}$ is not necessary in the preceding lemma. 
In fact, one may
   use the definition of
   $\exp:\mathfrak{gl}(\Lambda^*V)\rightarrow\mathrm{GL}(\Lambda^*V)$. More 
precisely,
   define $s:\mathbb{R}\rightarrow\mathrm{GL}(\Lambda^*V)$ by
   $s(t)=\wedge^* \exp({t {\hat\phi}})$. Then one may check that $s$ is the 
1-parameter
   subgroup generated by $\mathcal{D}^{\hspace{1pt}*}{\hat\phi}$ (i.e. 
$\dot{s}(0)=\mathcal{D}^{\hspace{1pt}*}{\hat\phi}$)
   and $s(1)=\wedge^{*} \exp({\hat\phi})$.
   \end{rema}
   \begin{coro}
Let ${\hat\phi}:V \rightarrow V$ be a nilpotent linear mapping. Then for
  $r \geq 0$, $\exp(\mathcal{D}^{\hspace{1pt}r}{\hat\phi})=\wedge^{r} 
\exp({\hat\phi})$ on $\Lambda^{r}V$. In particular,
   if ${\hat\epsilon}:={\hat\phi}+\mathrm{id}$, then $\wedge^{r} {\hat\epsilon}-\mathrm{id}$ 
is similar to
    $\mathcal{D}^{\hspace{1pt}r}{\hat\phi}$.
	\end{coro}
   \begin{proof}
The first part follows immediately from the lemma. For the second part, we 
know from
   Lemma \ref{lemm:3} that $\exp({\hat\phi})-\mathrm{id}$ is similar to ${\hat\phi}$, 
hence $\exp({\hat\phi})$ is similar to
   ${\hat\phi}+\mathrm{id}={\hat\epsilon}$. So
   \begin{align*}
   \wedge^{r} {\hat\epsilon}-\mathrm{id} &\sim \wedge^{r} \exp({\hat\phi})-\mathrm{id}\\
   &=\exp({\mathcal{D}^{\hspace{1pt}r}{\hat\phi}})-\mathrm{id}\\
   & \sim \mathcal{D}^{\hspace{1pt}r}{\hat\phi} .
      \end{align*}
   \end{proof}
   
   \begin{prop}
   Let $\sigma$ be the Anzai transformation on $\Bbb{T}^n$ and $\sigma_{*}$
   be the corresponding induced homomorphism on 
$K_*(\mathcal{C}(\mathbb{T}^n))=\Lambda^*\Bbb
   {Z}^n$. Let $\hat\sigma$ be the restriction of $\sigma_{*}$ to 
$\mathbb{Z}^n$ and consider
   the linear mapping $\hat\sigma\otimes{1}$ on $V:=\Bbb{Z}^n 
\otimes\Bbb{C}$. Take
   ${\hat\varphi}=\hat\sigma\otimes{1}-\mathrm{id}$ and 
$\mathcal{D}^{\hspace{1pt}r}{\hat\varphi}=
\mathcal{D}^{\hspace{1pt}*}{\hat\varphi}|_{\Lambda^{r}V}$ 
as above. Then
     
$$a_{n,r}=\mathrm{rank}\ker(\wedge^{r}\hat\sigma-\mathrm{id})=
\dim\ker\mathcal{D}^{\hspace{1pt}r}{\hat\varphi}.$$
   \end{prop}
\begin{proof}
Using the preceding corollary one has $\wedge^{r} {(\hat\sigma\otimes 1)}-\mathrm{id} 
\sim
\mathcal{D}^{\hspace{1pt}r}{\hat\varphi}$. Therefore
\begin{equation*}
\mathrm{rank}\ker(\wedge^{r}\hat\sigma-\mathrm{id})=\dim\ker(\wedge^{r}
(\hat\sigma\otimes 1)-\mathrm{id})=
\dim\ker\mathcal{D}^{\hspace{1pt}r}{\hat\varphi}.
\end{equation*}
\end{proof}

Now, to compute $a_{n,r}$ as explicitly as possible, we find a connection to 
representation
theory of $\mathfrak{sl}_2(\mathbb{C})$. First, we recall some definitions 
and properties.\\

Let $\mathfrak{sl}_2(\mathbb{C})$ denote the special linear algebra over 
$\mathbb{C}^2$ $$\mathfrak{sl}_2(\mathbb{C}):=
\{a \in M_2(\mathbb{C})\mid tr(a)=0\}.$$ We know that 
$\mathfrak{sl}_2(\mathbb{C})$ is a $3$-dimensional
semisimple complex Lie algebra. One can check that
$$
\mathfrak{B}:=\{h:=\left[\begin{array}{cc}
1 & 0\\
0 & -1\\
\end{array}
\right],e:=\left[\begin{array}{cc}
0 & 1\\
0 & 0\\
\end{array}
\right],f:=\left[\begin{array}{cc}
0 & 0\\
1 & 0\\
\end{array}
\right]\}$$
is a basis for this Lie algebra. 
Let $V$ be an \text{$n$-dimensional} complex vector space with a basis 
$\{e_1,\ldots,e_n \}$.
The following equalities for $i=1,\ldots,n$ define a representation
$\pi_n:\mathfrak{sl}_2(\mathbb{C})\rightarrow \mathfrak{gl}(V)$

\begin{itemize}
\item [(a)]$\pi_n(h)e_i=(2i-n-1)e_i$
\item [(b)]$\pi_n(e)e_i=i(n-i)e_{i+1}$; $(e_{n+1}:=0)$
\item [(c)]$\pi_n(f)e_i=e_{i-1}$; $(e_0:=0)$.
\end{itemize}
So, the relationship between our topics becomes clear. In fact, we know that
$\hat\sigma(e_i)=e_i+e_{i-1}$. Thus $\hat\varphi(e_i)=e_{i-1}$
, and hence $\hat\varphi=\pi_n(f)$.\\

Now recall the following theorem \cite[p. 33]{jeH72}.
\begin{theo}
Let $\pi_n$ be the representation described above. Then
\begin{itemize}
\item [\textnormal{(i)}]
$\pi_n$ is an irreducible representation of
$\mathfrak{sl}_2(\mathbb{C})$.
\item [\textnormal{(ii)}]
Any \text{$n$-dimensional} representation of $\mathfrak{sl}_2(\mathbb{C})$
is equivalent to $\pi_n$.
\item [\textnormal{(iii)}]
Suppose $V$ is any finite dimensional $\mathfrak{sl}_2(\mathbb{C})$-module
and define $$V_\alpha=\{v\in V \mid h.v=\alpha~v\}$$
for $\alpha\in\mathbb{C}$. Then $V$ decomposes into a direct sum of 
irreducible submodules
and in any such decomposition, the number of summands is precisely
$\dim V_0 + \dim V_1$.
\end{itemize}
\end{theo}

For convenience, put $\pi=\pi_n$ and extend $\pi$
to $\pi^r:\mathfrak{sl}_2(\mathbb{C})\rightarrow \mathfrak{gl}(\Lambda^r V)$ 
with
$\pi^1=\pi$. More precisely, for every $X \in \mathfrak{sl}_2(\mathbb{C})$, 
define
$$\pi^r(X)(v_1 \wedge \ldots\wedge v_r)=(\pi(X) v_1)\wedge v_2 
\wedge\ldots\wedge v_r
+\ldots +v_1 \wedge\ldots\wedge v_{r-1}\wedge (\pi(X)v_r).$$
Thus $\mathcal{D}^{\hspace{1pt}r}{\hat\varphi}=\pi^r(f)$.\\

Now we are ready to compute $a_{n,r}=\dim\ker 
\mathcal{D}^{\hspace{1pt}r}{\hat\varphi}$,
which is essential for studying the rank $a_n$ of the $K$-groups of 
$\mathcal{A}_{n,\theta}$.
We need some notation first.\\

\begin{nota}\label{nota:3}
Let $n,k,r$ be positive integers. $P(n,r,k)$ denotes the number
of partitions of $k$ to $r$ distinct positive integers not greater than $n$. 
In other
words
$$P(n,r,k)=\textnormal{card}\{(i_1,\ldots,i_r)\mid i_1+\ldots+i_r=k, 
1 \leq i_1<\ldots<i_r
\leq n \}.$$
We conventionally take $P(n,0,0)=1$ and $P(n,r,0)=P(n,0,k)=0$ for $r,k \geq 
1$.
\end{nota}
\begin{prop}
With the above notation, $\dim \ker 
\mathcal{D}^{\hspace{1pt}r}{\hat\varphi}=P(n,r,[r(n+1)/2])$, in which $[x]$ 
denotes
the greatest integer not greater than $x$.
\end{prop}
\begin{proof}
Following Weyl's theorem \cite[p. 28]{jeH72}, $\pi^r$ is completely irreducible, 
which means that $\Lambda^r V=
\oplus_{p=1}^N W_p$, where the $W_p$'s are $\pi^r$-invariant irreducible 
subspaces of
$\Lambda^r V$ and $N$ is the number of such subspaces, which following the 
preceding theorem is
equal to $\dim E_0+\dim E_1$, where $E_j=\{v \in \Lambda^r V \mid 
\pi^r(h)v=jv\}$. On
the other hand, $\ker 
\mathcal{D}^{\hspace{1pt}r}{\hat\varphi}=\oplus_{p=1}^N 
\ker({\mathcal{D}^{\hspace{1pt}r}{\hat\varphi}}|_{\tiny{W_p}})$, hence
$\dim\ker \mathcal{D}^{\hspace{1pt}r}{\hat\varphi}=\sum_{p=1}^N 
\dim\ker(\mathcal{D}^{\hspace{1pt}r}{\hat\varphi}|_{W_p})=N$, since 
$\dim\ker(\mathcal{D}^{\hspace{1pt}r}{\hat\varphi}|_{W_p})=1$
, so $\dim\ker \mathcal{D}^{\hspace{1pt}r}{\hat\varphi}=\dim E_0+\dim E_1$. 
To compute the last term, note that from the
preceding theorem we have
$$\pi^r(h)(e_{i_1}\wedge\ldots\wedge e_{i_r})=(2(i_1+\ldots+i_r)-r(n+1))
e_{i_1}\wedge\ldots\wedge e_{i_r}$$
so for even $r(n+1)$, $E_1=\{0\}$ and $\dim E_0=P(n,r,r(n+1)/2)$ and for odd
$r(n+1)$, $E_0=\{0\}$ and $\dim E_1=P(n,r,r(n+1)/2-1)$. To summarize,
$\dim \ker\mathcal{D}^{\hspace{1pt}r}{\hat\varphi}=N=\dim E_0+\dim E_1 
=P(n,r,[r(n+1)/2])$.
\end{proof}

Therefore, as the main result of this part, we have the following theorem.
\begin{theo}
$a_{n,r}=P(n,r,[r(n+1)/2])$ for $r=0,1,\ldots,n$.
\end{theo}
\begin{proof}
Use Propositions 6.4 and 6.5.
\end{proof}

As a result, we can prove that $\{a_n \}$ is a strictly increasing sequence.
We need a lemma first.
\begin{lemm}
$P(n+1,r,k+s)\geq P(n,r,k)$ for $s=0,1,\ldots,r$.
\end{lemm}
\begin{proof}
For $s=0$, the proof is clear. Now let $1 \leq s \leq r$ and suppose that
$(j_1,\ldots,j_r)$ is a partition of $k$ such that $1 \leq j_1 <\ldots<j_r 
\leq n$.
Now define $i_q:=j_q$ for $1 \leq q \leq r-s$ and  $i_q:=j_q+1$
for $r-s+1 \leq q \leq r$. Then $(i_1,\ldots,i_r)$ is a partition of $k+s$ 
and
$1 \leq i_1<\ldots<i_r \leq n+1$. Thus $P(n+1,r,k+s)\geq P(n,r,k)$.
\end{proof}

\begin{prop}\label{prop:2}
$\{a_n \}$ is a strictly increasing sequence.
\end{prop}
\begin{proof}
First note that $a_{n,0}=a_{n,n}=P(n,0,0)=P(n,n,n(n+1)/2)=1$ and from the 
preceding theorem
we have $a_n=\sum_{r=0}^n P(n,r,[r(n+1)/2])$. Now we prove that
for every $m \in\Bbb{N}$, $a_{2m+1}>a_{2m}>a_{2m-1}$. Applying the lemma to 
the terms
of the following equalities yields the result.
\begin{align*}
a_{2m+1}&=1+\sum_{r=0}^m P(2m+1,2r,2rm+2r)+\\
&~~~~~~~~~~~~~~~~~~~~~~~~~~~~~~~\sum_{r=0}^{m-1} P(2m+1,2r+1,2rm+2r+m+1),
\end{align*}
\begin{align*}
a_{2m}&=\sum_{r=0}^m P(2m,2r,2rm+r)+\sum_{r=0}^{m-1} P(2m,2r+1,2rm+m+r)\\
&=1+\sum_{r=0}^{m-1} P(2m,2r,2rm+r)+\sum_{r=0}^{m-1} P(2m,2r+1,2rm+m+r),
\end{align*}
\begin{align*}
a_{2m-1}&=\sum_{r=0}^{m-1} P(2m-1,2r,2rm)+\sum_{r=0}^{m-1} 
P(2m-1,2r+1,2rm+m).
\end{align*}
\end{proof}

\subsection{Generating functions for $a_n$}
In this part, we express the rank of the $K$-groups of 
$\mathcal{A}_{n,\theta}$ as explicitly as possible. In fact, 
we present them as the constant terms in the polynomial 
expansions of certain functions. First of all, we need the 
following basic lemma.

\begin{lemm}
Let $P(n,r,k)$ denote the number of partitions of $k$ to $r$ distinct 
positive
integers not greater than $n$. Then $P(n,r,k)$ is the coefficient of $u^r 
t^k$
in the polynomial expansion of $F_n(u,t):=\prod_{i=1}^n (1+u t^i)$. In other 
words
$$\sum_{r,k \geq 0}P(n,r,k)u^r t^k=\prod_{i=1}^n (1+u t^i).$$
\end{lemm}
\begin{proof}
\begin{align*}
\prod_{i=1}^n (1+u t^i)&=1+\sum_{r=1}^n
\underset{1 \leq i_1<\ldots<i_r \leq n}
{\sum_{\small{(i_1,\ldots,i_r)}}}
{(ut^{i_1})\ldots (ut^{i_r})}\\
&=1+\sum_{r=1}^n \sum_{k \geq 1}P(n,r,k)u^r t^k\\
&=\sum_{r,k \geq 0}P(n,r,k)u^r t^k
\end{align*}
\end{proof}

Now, we have the following expressions for the rank $a_n$ of the $K$-groups
of $\mathcal{A}_{n,\theta}$.

\begin{theo}\label{theo:6}
Let $a_n=\mathrm{rank}\hspace{2pt}K_0(\mathcal{A}_{n,\theta})=
\mathrm{rank}\hspace{2pt}K_1(\mathcal{A}_{n,\theta})$. Then \\

\begin{itemize}
\item for odd $n$, $a_n$ is the constant term in the polynomial expansion of
$$\prod_{i=1}^n (1+t^{i-\frac{n+1}{2}})$$
\item for even $n$, $a_n$ is the constant term in the polynomial expansion 
of
$$(1+t^{\frac{1}{2}})\prod_{i=1}^n (1+t^{i-\frac{n+1}{2}}).$$
\end{itemize}
\end{theo}
\begin{proof}
We know that $a_n=\sum _{r=0}^n a_{n,r}$ and $a_{n,r}=P(n,r,[r(n+1)/2])$.
Now let $n=2m-1$ be an odd number. So, $a_n=\sum_{r=0}^{2m-1}P(2m-1,r,rm)$.
Take $y=u t^m$. From the preceding lemma, we have
$$F_n(u,t)=F_{2 m-1}(y t^{-m},t)=\prod_{i=1}^{2m-1}(1+y t^{i-m})=
\sum_{r,k \geq 0}P(2m-1,r,k)y^r t^{k-rm}$$
so for $y=1$ we have
$$\prod_{i=1}^n (1+t^{i-\frac{n+1}{2}})=\prod_{i=1}^{2m-1}(1+t^{i-m})=
\sum_{r,k \geq 0}P(2m-1,r,k)t^{k-rm},$$
which yields the result.\\

For even $n$, say $n=2m$, we have
\begin{align*}
a_n=a_{2m}&=\sum_{r=0}^{2m}P(2m,r,[r(m+\frac{1}{2})])\\
          &=\sum_{r=0}^m P(2m,2r,r(2m+1))+\sum_{r=0}^{m-1} 
P(2m,2r+1,2rm+m+r)\\
		  &=:A_m + B_m.
\end{align*}
First, let's determine $A_m$. Note that from the lemma we have
\begin{align*}
\frac{1}{2}\{\prod_{i=1}^{2m}(1+ut^i)+\prod_{i=1}^{2m}(1-ut^i)\}&=
\sum_{r,k \geq 0}P(2m,r,k)\{\frac{1+(-1)^r}{2}\}u^r t^k\\
&=\sum_{r,k \geq 0}P(2m,2r,k)u^{2r} t^k.
\end{align*}
If we define $y:=u^2 t^{2m+1}$, we have the following identity
\begin{multline*}\frac{1}{2}
\{\prod_{i=1}^{2m}(1+y^{\frac{1}{2}}t^{i-(m+\frac{1}{2})})
+\prod_{i=1}^{2m}(1-y^{\frac{1}{2}}t^{i-(m+\frac{1}{2})})\}=\\
\sum_{r,k \geq 0}P(2m,2r,k)y^r t^{k-r(2m+1)},
\end{multline*}
which for $y=1$ yields
$$\frac{1}{2}\{\prod_{i=1}^{2m}(1+t^{i-(m+\frac{1}{2})})
+\prod_{i=1}^{2m}(1-t^{i-(m+\frac{1}{2})})\}=
\sum_{r,k \geq 0}P(2m,2r,k)t^{k-r(2m+1)},$$
hence $A_m$ is the constant term in polynomial expansion of
$$\frac{1}{2}\{\prod_{i=1}^{2m}(1+t^{i-(m+\frac{1}{2})})
+\prod_{i=1}^{2m}(1-t^{i-(m+\frac{1}{2})})\}.$$
Similarly, for $B_m$ we have
\begin{align*}
\frac{1}{2}\{\prod_{i=1}^{2m}(1+ut^i)-\prod_{i=1}^{2m}(1-ut^i)\}&=
\sum_{r,k \geq 0}P(2m,r,k)\{\frac{1-(-1)^r}{2}\}u^r t^k\\
&=\sum_{r,k \geq 0}P(2m,2r+1,k)u^{2r+1} t^k.
\end{align*}
If we define $y^2:=u^2 t^{2m+1}$, we have the following identities
\begin{multline*}
\frac{1}{2}\{\prod_{i=1}^{2m}(1+y^{\frac{1}{2}}t^{i-(m+\frac{1}{2})})
-\prod_{i=1}^{2m}(1-y^{\frac{1}{2}}t^{i-(m+\frac{1}{2})})\}\\
=\sum_{r,k \geq 0}P(2m,2r+1,k)y^{2r+1} t^{k-(2rm+r+m)-\frac{1}{2}}\\
=t^{-\frac{1}{2}}\sum_{r,k \geq 0}P(2m,2r+1,k)y^{2r+1} t^{k-(2rm+r+m)},
\end{multline*}
which for $y=1$ yields
\begin{multline*}
\frac{t^{\frac{1}{2}}}{2}\{\prod_{i=1}^{2m}(1+t^{i-(m+\frac{1}{2})})
-\prod_{i=1}^{2m}(1-t^{i-(m+\frac{1}{2})})\}=\\
\sum_{r,k \geq 0}P(2m,2r+1,k)t^{k-(2rm+r+m)},
\end{multline*}
hence $B_m$ is the constant term in the polynomial expansion of
$$\frac{t^{\frac{1}{2}}}{2}\{\prod_{i=1}^{2m}(1+t^{i-(m+\frac{1}{2})})
-\prod_{i=1}^{2m}(1-t^{i-(m+\frac{1}{2})})\}.$$
Therefore $a_n=a_{2m}=A_m+B_m$ is the constant term in the polynomial
expansion of
\begin{multline*}
\frac{1}{2}\{\prod_{i=1}^{2m}(1+t^{i-(m+\frac{1}{2})})
+\prod_{i=1}^{2m}(1-t^{i-(m+\frac{1}{2})})\\+\frac{t^{\frac{1}{2}}}{2}
\prod_{i=1}^{2m}(1+t^{i-(m+\frac{1}{2})})-\frac{t^{\frac{1}{2}}}{2}
\prod_{i=1}^{2m}(1-t^{i-(m+\frac{1}{2})})\},
\end{multline*}
or equivalently, the constant term in the polynomial
expansion of
$$\frac{1}{2}\{(1+z)\prod_{i=1}^{2m}(1+z^{2i-(2m+1)})
+(1-z)\prod_{i=1}^{2m}(1-z^{2i-(2m+1)})\},$$
or equivalently, the constant term in the polynomial
expansion of
$$(1+z)\prod_{i=1}^{2m}(1+z^{2i-(2m+1)}),$$
or equivalently, the constant term in the polynomial
expansion of
$$(1+t^{\frac{1}{2}})\prod_{i=1}^{2m}(1+t^{i-(m+\frac{1}{2})}),$$
which yields the result.
\end{proof}

One can use this theorem to determine the asymptotic 
behavior of $a_n$.
\begin{coro}
$a_n \thicksim (\frac{24}{\pi})^{\frac{1}{2}}2^n n^{-\frac{3}{2}}$ when $n 
\rightarrow \infty$.
\end{coro}

  \section{Concluding remarks}
  \begin{enumerate}
\item The torsion parts of the $K$-groups of $\mathcal{A}_{n,\theta}$ seem 
much more
difficult to describe explicitly in terms of $n$. Nevertheless, it is an 
interesting
problem to find such descriptions. It is of interest to compute the 
\text{$K$-groups} of
$\mathcal{A}_{n,\theta}$ (or more general algebras $A_{F_{f,\theta}}$), 
since in the class
of $C^*$-algebras generated by uniquely ergodic minimal diffeomorphisms on a 
compact
manifold, $K$-theory is a complete invariant. More precisely, suppose that 
$M$
is a connected compact smooth manifold with $\dim(M)>0$ and \text{$h:M 
\rightarrow M$} is a
uniquely ergodic minimal diffeomorphism, and put $A:=\mathcal{C}(M)\rtimes_h 
\Bbb{Z}$.
Let $\tau$ be the trace induced by the unique invariant probability measure, 
and assume
that $\tau_*K_0(A)$ is dense in $\Bbb{R}$. Then the $4$-tuple
$$(K_0(A),K_0(A)_{+},[1_A],K_1(A))$$
is a complete algebraic invariant (called the \emph{Elliott invariant} of 
$A$) \cite{qL02}.
In this case, $A$ has stable rank one, real rank zero and tracial 
topological rank
zero in the sense of H. Lin \cite{hL00}. The order on $K_0(A)$ is also 
determined
by the unique trace $\tau$.
\item The method used in section $6$ above for computing the $K$-groups of 
the
transformation group $C^*$-algebras of the tori may be extended to more 
general settings. Let $G$ be a
  compact connected Lie group with $\pi_1(G)$ torsion-free.
  Then $K^*(G)$ is torsion-free and can be given the structure of a
  $\mathbb{Z}_2$-graded Hopf algebra over the integers \cite{lH67}. 
Moreover, regarded as
  a Hopf algebra,
  $K^*(G)$ is the exterior algebra on the module of the primitive elements, 
which are
  of degree 1. The module of the primitive elements of $K^*(G)$ may also be 
described
  as follows. Let $U(n)$ denote the group of unitary matrices of order $n$ 
and let
  $U:=\cup_{n=1}^\infty U(n)$ be the stable unitary group.  Any unitary 
representation $\rho:G \rightarrow U(n)$,
  by composition with the inclusion $U(n)\subset U$, defines a homotopy 
class
  $\beta(\rho)$ in $[G,U]=K^1(G)$. The module of the primitive elements in 
$K^1(G)$
  is exactly the module generated by all classes $\beta(\rho)$ of this type.
  If in addition, $G$ is semisimple and simply connected of rank $l$, there 
are
  $l$ basic irreducible representations $\rho_1,\ldots,\rho_l$, whose 
maximum weights
  $\lambda_1,\ldots,\lambda_l$ form a basis for the character group 
$\mathbf{\hat{T}}$
  of the maximal torus $\textbf{T}$ of $G$ and the classes 
$\beta(\rho_1),\ldots,
  \beta(\rho_l)$ form a basis for the module of the primitive elements in 
$K^1(G)$ and
  $K^*(G)=\Lambda^*(\beta(\rho_1),\ldots,\beta(\rho_l))$. In any case, to 
compute
  $K_*(\mathcal{C}(G)\rtimes_{\alpha}\mathbb{Z})$, it is sufficient to 
determine
  the homotopy classes of $\alpha\circ\rho$ for irreducible representations 
$\rho$
  of $G$ in terms of $\beta(\rho)$'s.

  \item There is a relation between the $K$-theory of transformation group
  \text{$C^*$-algebras} of the tori and the topological $K$-theory of compact 
nilmanifolds. In fact
  let $\alpha=(\boldsymbol{t},\mathsf{A})$ be an affine transformation on 
$\Bbb{T}^n$
  satisfying the conditions (1) through (4) in Remark \ref{rema:3}. Then it 
has been shown
  in \cite{fjH63} that $\alpha$ is conjugate (in the group of affine 
transformations of
   $\Bbb{T}^n$) to the transformation 
$\alpha'=(\boldsymbol{t'},\mathsf{A}')$,
   where $\mathsf{A}'$ has an upper triangular matrix whose bottom right 
   $k \times k$
   corner is the identity matrix $\mathsf{I}_k$ and 
$\boldsymbol{t}'=(0,\ldots,0,t'_1,
   \ldots,t'_k)$. $\alpha'$ is called a \emph{standard form} for $\alpha$ 
\cite{jP86}.
   Assume that $\alpha$ is given in standard form. J. Packer associates
   an \emph{induced flow} $(\Bbb{R}, N/\Gamma)$ to the flow 
$(\Bbb{Z},\Bbb{T}^n)$
   generated by $\alpha$, where $N$ is a simply connected nilpotent Lie 
group of dimension $n+1$,
  $\Gamma$ is a cocompact subgroup of $N$, and the action of $\Bbb{R}$ is 
given by
  translation on the left by $\exp sX$ for $s \in\mathbb{R}$ and some $X 
\in\mathfrak{n}$, the Lie
  algebra of $N$. One of the most important facts is
  that the \text{$C^*$-algebra} $\mathcal{C}(N/\Gamma)\rtimes_\beta\Bbb{R}$ 
corresponding to the
  induced flow is strongly Morita equivalent to
  $\mathcal{C}(\Bbb{T}^n)\rtimes_\alpha\Bbb{Z}$ 
  \cite[Proposition 3.1]{jP86}. Consequently, one has
  \begin{equation}
  K_i(\mathcal{C}(\Bbb{T}^n)\rtimes_\alpha\Bbb{Z})\cong
  K_i(\mathcal{C}(N/\Gamma)\rtimes_\beta\Bbb{R})\cong 
K^{1-i}(N/\Gamma);~~~i=0,1.
  \end{equation}
  The second isomorphism here, is the Connes' Thom
  isomorphism. So the $K$-theory of 
$\mathcal{C}(\Bbb{T}^n)\rtimes_\alpha\Bbb{Z}$
  is converted to the topological $K$-theory of the compact nilmanifold 
$N/\Gamma$.
  Following the proof of Proposition 3.1 in \cite{jP86}, one can conclude 
that for the special case of the Anzai flows, $N=F_{n-1}$ and 
$\Gamma=\mathfrak{D}_{n-1}$
  which were defined in Section 2.\\
  On the other hand, following \cite[Theorem 3.6]{jR84}, one has the 
following isomorphism
  \begin{equation}
  K_i(C^*(\Gamma))\cong K^{i+n+1}(N/\Gamma);~~~i=0,1.
  \end{equation}
  Combining (14) and (15) one gets
  \begin{equation}
  K_i(\mathcal{C}(\Bbb{T}^n)\rtimes_\alpha\Bbb{Z})\cong 
K^{i+1}(N/\Gamma)\cong K_{i+n}(C^*(\Gamma));~~~i=0,1.
  \end{equation}
Using the above isomorphisms, one can relate the algebraic invariants of
the involved $C^*$-algebras and topological information about the 
corresponding
nilmanifold. For example, since $N/\Gamma$ is a classifying space for 
$\Gamma$,
one has the following isomorphisms
\begin{equation}
H^*_{\mathrm{dR}}(N/\Gamma)\cong \Check{H}^*(N/\Gamma,\Bbb{R})\cong 
H^*(\Gamma,\Bbb{R})
\cong H^*(N,\Bbb{R})\cong H^*(\mathfrak{n},\Bbb{R}),
\end{equation}
where $H^*_{\mathrm{dR}}(N/\Gamma)$ denotes the de Rham cohomology of the 
manifold
$N/\Gamma$, $H^*(N/\Gamma,\Bbb{R})$ denotes the \v{C}ech cohomology of 
$N/\Gamma$ with
coefficients in $\Bbb{R}$, $H^*(\Gamma,\Bbb{R})$ denotes the group 
cohomology of $\Gamma$ with
coefficients in the trivial $\Gamma$-module $\Bbb{R}$, $H^*(N,\Bbb{R})$ 
denotes the
Moore cohomology group of $N$ (as a locally compact group) with coefficients 
in the trivial
Polish $N$-module $\Bbb{R}$,
and $H^*(\mathfrak{n},\Bbb{R})$ denotes the cohomology of the Lie algebra 
$\mathfrak{n}$
with coefficients in the trivial $\mathfrak{n}$-module $\Bbb{R}$. Now using 
the Chern
isomorphisms 
$\mathrm{ch}_0:K^0(N/\Gamma)\otimes\Bbb{Q}
\rightarrow\Check{H}^{\mathrm{even}}(N/\Gamma,\Bbb{Q})$ and 
$\mathrm{ch}_1:K^1(N/\Gamma)\otimes\Bbb{Q}\rightarrow
\Check{H}^{\mathrm{odd}}(N/\Gamma,\Bbb{Q})$, one concludes that
the even and odd cohomology groups stated in (17) are all isomorphic to 
$\mathbb{R}^{k}$,
where $k$ is the (common) rank of the $K$-groups of
$\mathcal{C}(\Bbb{T}^n)\rtimes_\alpha\Bbb{Z}$ as in Corollary \ref{coro:3}.
As an example, if $N=F_{n-1}$, $\Gamma=\mathfrak{D}_{n-1}$, and
$\mathfrak{n}=\mathfrak{f}_{n-1}$, then the even and odd cohomology groups
stated in (17) are all isomorphic to $\mathbb{R}^{a_{n}}$,
where $a_{n}$ is the rank of the $K$-groups of $\mathcal{A}_{n,\theta}$ that
was studied in detail in Section 6. \\
Conversely, one may use the topological
tools for $N/\Gamma$ to get some information about
$\mathcal{C}(\Bbb{T}^n)\rtimes_\alpha\Bbb{Z}$ and $C^*(\Gamma)$.
For example, we know that $N/\Gamma$ as a compact nilmanifold
can be constructed as a principal $\Bbb{T}$-bundle over a lower
dimensional compact nilmanifold \cite{vG97}. Then we can compute the 
topological
$K$-groups of $N/\Gamma$ using the six term Gysin exact sequence
\cite[IV.1.13, p. 187]{mK78}. As an example, one can see that 
$F_{n-1}/\mathfrak{D}_{n-1}$
is a principal \text{$\Bbb{T}$-bundle} over $F_{n-2}/\mathfrak{D}_{n-2}$, 
and the corresponding
Gysin exact sequence is in fact the topological version of the 
Pimsner-Voiculescu exact
sequence for the crossed product
$\mathcal{A}_{n,\theta}\cong 
\mathcal{A}_{n-1,\theta}\rtimes_\alpha\mathbb{Z}$ as in
Theorem \ref{theo:1}.
\end{enumerate}

\textbf{Acknowledgement.}
This paper is part of the first author's
 Ph.D. thesis in Tarbiat Modarres University under the
 supervision of Professor Ali Reza Medghalchi. The first author would
 like to thank Professor Medghalchi for his moral support and 
scientific guidance during his studies.
 He would also like to thank Graham Denham, George
 Elliott, Herve Oyono-Oyono, N. C. Phillips, and Tim
 Steger for helpful discussions.

\end{document}